\documentstyle[amssymb,amsfonts]{amsart}

\newcommand\crit{{\operatorname{crit}}}

\newcommand\lin{{\operatorname{lin}}}
\newcommand\loc{{\operatorname{loc}}}
\newcommand\lo{{\operatorname{lo}}}
\newcommand\hi{{\operatorname{hi}}}
\newcommand\stress{{\operatorname{T}}}
\newcommand\R{{{\mathbf R}}}
\renewcommand\H{{{\mathbf H}}}
\newcommand\Z{{{\mathbf Z}}}
\newcommand\C{{{\mathbf C}}}
\renewcommand\S{{{\mathcal S}}}
\newcommand{\nabb}{\mbox{$\nabla \mkern-13mu /$\,}}
\newcommand\N{{{\mathcal N}}}
\newcommand\g{{{\mathfrak g}}}
\newcommand\eps{{\varepsilon}}
\newenvironment{proof}{\noindent {\bf Proof} }{\endprf\par}
\def \endprf{\hfill  {\vrule height6pt width6pt depth0pt}\medskip}
\def\emph#1{{\it #1}}
\def\textbf#1{{\bf #1}}

\theoremstyle{plain}
  \newtheorem{theorem}[subsection]{Theorem}
  
  \newtheorem{proposition}[subsection]{Proposition}
  \newtheorem{lemma}[subsection]{Lemma}
  \newtheorem{corollary}[subsection]{Corollary}

\theoremstyle{remark}
  \newtheorem{remark}[subsection]{Remark}

\theoremstyle{definition}
  \newtheorem{definition}[subsection]{Definition}

\include{psfig}

\begin{document}

\title[Nonlinear evolution equations]{Global behaviour of nonlinear dispersive and wave equations}
\author{Terence Tao}
\address{Department of Mathematics, UCLA, Los Angeles CA 90095-1555}
\email{ tao@@math.ucla.edu}

\vspace{-0.3in}
\begin{abstract}
We survey recent advances in the analysis of the large data global (and asymptotic) behaviour of nonlinear dispersive equations such as the nonlinear wave (NLW), nonlinear Schr\"odinger (NLS), wave maps (WM), Schr\"odinger maps (SM), generalised Korteweg-de Vries (gKdV), Maxwell-Klein-Gordon (MKG), and Yang-Mills (YM) equations.  The classification of the nonlinearity as \emph{subcritical} (weaker than the linear dispersion at high frequencies), \emph{critical} (comparable to the linear dispersion at all frequencies), or \emph{supercritical} (stronger than the linear dispersion at high frequencies) is fundamental to this analysis, and much of the recent progress has pivoted on the case when there is a critical conservation law.  We discuss
how one synthesises a satisfactory critical (scale-invariant) global theory, starting the basic building blocks of perturbative analysis, conservation laws, and monotonicity formulae, but also incorporating more advanced (and recent) tools such as gauge transforms, concentration-compactness, and induction on energy. 
\end{abstract}

\maketitle

\section{Introduction}

The purpose of this survey is to discuss recent progress in understanding the global and asymptotic behaviour of various model nonlinear evolution equations of dispersive or wave type (as opposed to parabolic, transport, or kinetic equations) on Euclidean spacetimes $\R \times \R^d$ for various dimensions $d$.  These equations are \emph{semilinear}, meaning that they are perturbations of a linear dispersive or wave equation by a nonlinearity of lower order (i.e. using fewer derivatives than the linear part of the equation); the evolution is then a competition between the linear part of the equation (which tends to disperse the solution) and the nonlinear part (which can either focus or defocus the solution, depending on the sign of the nonlinearity).  They are also Hamiltonian (and hence time-reversible), in contrast to parabolic equations (such as the heat equation, or Navier-Stokes) which are dissipative and non-time-reversible.  The evolution can be expected to broadly be 
a combination of one of three forms:

\begin{itemize}

\item \textbf{Linearly dominated behaviour.}  In some cases the linear effects dominate the nonlinear effects, and the solution exists globally and converges asymptotically to a linear solution (which itself should disperse to zero).  In such cases one tends to have very good spacetime bounds (basically, the nonlinear solution should obey almost the same bounds as the linear solution) and a complete scattering theory for the equation. This scenario tends to occur for small data, high regularities, short times, low dimensions, and weak (low-power) nonlinearities.

\item \textbf{Nonlinearly dominated behaviour.} In opposition to the previous case, it is possible for the nonlinear effects to dominate the linear effects.  In ``focusing'' cases, this typically causes the solution to become very unstable, and singularities develop in finite time or even instantaneously.  In ``defocusing'' cases, the solution is still rather unstable for medium times, but typically the nonlinearity acts to disperse the solution, at which point the evolution switches over to linearly dominated behaviour.
This scenario tends to occur for large data, low regularities, long times, high dimensions, and strong (high-power) nonlinearities.

\item \textbf{Intermediate behaviour.}  A third regime of behaviour emerges when the nonlinear and linear effects are roughly in balance.  The most notable example of this are the \emph{soliton} solutions in focusing (or at least non-defocusing) equations, which are typically stationary or traveling wave solutions in which the dispersive effect of the linear equations 

\end{itemize}

A large part of the analytical theory of these equations revolves around how to rigorously classify, based on the equation and on the class of initial data involved, whether the global evolution of the equation exhibits linear behaviour or nonlinear behaviour.  In doing so, two basic features of these equation have proven to be of vital importance.  The first are the \emph{conserved quantities} (and to a lesser extent, the \emph{monotone quantities}) of the evolution, and more precisely those quantities which are \emph{coercive} (in that they provide non-trivial upper bounds on the size of the solution) or at least positive semi-definite to top order.  In the large data theory, the conserved and monotone quantities determine what control one can retain on the solution after long times.  The second is the natural scale-invariance (or approximate scale invariance) of the equation, which provides an identification between the fine-scale and coarse-scale behaviour of the evolution.  Using this invariance, one can classify the conservation laws as being either \emph{subcritical} (strong at fine scales, weak at coarse scales), \emph{critical} (scale-invariant), or \emph{supercritical} (strong at coarse scales, but weak at fine scales).
One can similarly classify regularity classes (such as Sobolev spaces $H^s_x(\R^d)$) as being subcritical, critical, or supercritical for a certain equation.
The equations with critical conservation laws provide a context where the nonlinear and linear parts are roughly comparable in strength, and represent the frontier of current technology for analysing large data global behaviour of evolution equations.

After the classification of equations and their conservation laws as being subcritical, critical, or supercritical, the next most important distinction is whether the equation is \emph{defocusing}, \emph{focusing}, or \emph{neither}.  These terms do not have a fully precise meaning, but roughly speaking in a defocusing equation the nonlinear component of the equation is typically aligned to have the ``same sign'' as the linear component, thus (hopefully) amplifying the dispersive effects of the linear equation, whereas in the focusing case the opposite is true, and the dispersive effects can be attenuated, halted (to cause stationary or travelling wave solutions such as soliton solutions) or even reversed (to cause blowup).  In some cases (e.g. for the Korteweg-de Vries, Maxwell-Klein-Gordon and Yang-Mills equations) the nonlinearity does not have a preferred sign, and 
can act either to aid the dispersion or to counteract it.

One can seek to understand the behaviour of these solutions either at high regularities (smooth solutions) or low regularities (rough solutions).  In many applications, it is the smooth solutions which are of importance; but even if one is only ultimately interested in high regularity solutions, it is often worthwhile to fully develop the low regularity theory, as the estimates obtained as a consequence of that theory are often extremely useful in controlling the global and asymptotic behaviour of
smooth solutions, and in particular in obtaining precise criteria as to whether blowup or other bad behaviour will occur from smooth initial data.  In any event, in cases where the key conserved quantity is critical, the smooth theory and the critical-regularity theory are often very closely related, and many of the deepest results concerning smooth solutions to these equations arose directly from, or were at least inspired by, the critical-regularity theory\footnote{This does not necessarily mean however, that one has to abandon the classical concept of solution for weaker notions of solution (such as distributional solutions); in many cases, one can proceed by working entirely in the category of smooth solutions, so long as one is always seeking estimates which are \emph{scale-invariant} in nature, and in particular not reliant on the high regularity norms of the solution, except to justify certain formal computations or to run qualitative arguments such as continuity arguments.}.

There are a large number of interesting nonlinear evolution equations of dispersive or wave type.
In contrast to other fields of mathematics, it is not always profitable to try to treat all of these at once by working with an abstract class of PDE; while a limited amount of generalisation is possible, each individual PDE typically has its own ``personality'' which requires separate treatment, especially when working with the particularly delicate issue of global large data theory at critical regularities.  On the other hand, the \emph{techniques}, heuristics, and principles for analysing these PDE are remarkably constant from one equation to the next.  Furthermore, we shall see that there are several connections and analogies, both formal and heuristic, between different equations.  Thus it is important to study these equations both individually (at the rigorous level) and collectively (at the informal level).

With a few notable exceptions (KdV, mKdV, 1-dimensional cubic NLS, 1-dimensional wave maps, restricted classes of Yang-Mills), the majority of equations discussed here are not completely integrable, and almost certainly not reducible via algebraic transformations to a linear evolution; thus there is essentially no hope of finding exact solutions to these equations from \emph{general} initial data via some algebraic formula, although there are certainly many important \emph{special} exact solutions (e.g. solitons, highly symmetric solutions, or the trivial vacuum solution $0$) which play major roles in the subject and provide important examples and intuition.  In the absence of exact formulae for general solutions, the analytical theory instead revolves around \emph{qualitative} and \emph{quantitative} properties of the solutions.  Qualitative properties include the fundamental question of \emph{wellposedness} (existence, uniqueness, and continuous dependence of the solution on the initial data in some prescribed data class), as well as regularity, approximation by smooth solutions, justification of formal algebraic manipulations (e.g. conservation laws), and asymptotics at infinity.  At very low regularities, even the utterly basic (but surprisingly subtle) question of what it even \emph{means} to be a solution has to be properly addressed.  Quantitative properties typically involve estimating various spatial or spacetime norms of the solution (e.g. Sobolev or Lebesgue norms) in terms of various norms of the initial data (such as the mass and energy).  The two types of properties are often closely intertwined; one needs quantitative estimates in order to conclude enough convergence or continuity to justify a qualitative argument, and conversely qualitative results are often needed to justify quantitative computations; in many cases one needs a bootstrap, continuity, or iteration argument to produce both the quantitative and qualitative results simultaneously.  Our focus here 
shall be more on the ``hard'' quantitative components of recent results; the ``soft'' qualitative arguments are also a necessary component of these results, but these tend to be relatively routine once the quantitative estimates are obtained. In particular we shall often assume that a solution has been \emph{a priori} given to us, and is already sufficiently regular to justify all formal computations, but lacks strong quantitative estimates; we shall then work hard to establish such quantitative estimates 
(known as \emph{a priori estimates}).  Once these estimates are obtained, there are a number of ``soft'' 
techniques (approximation, penalisation, iterative methods, continuity methods, use of 
higher-regularity wellposedness theory) to remove the a priori restriction and show existence and uniqueness of 
solutions with the desired bounds from all data in a given class.  While these arguments are necessary and sometimes subtle, the 
technical issues they raise tend to distract from the physical intuition underlying the dynamics of these equations, and so
we will not dwell on them here.

\section{The model equations}

In this section we describe several model equations which we will discuss in this survey.  
There are many model nonlinear equations of dispersive or wave type which are of importance, but we shall select only some
particularly symmetric ones, in particular those which enjoy an exact translation-invariance and scaling-invariance, as these are slightly simpler to study analytically and already exhibit many of the key phenomena that one wishes to understand in this field.  Also, the presence of symmetries naturally leads one to special self-symmetric sub-classes of solutions (e.g. travelling wave solutions, self-similar solutions, spherically symmetric solutions) of interest.  We shall also focus attention on those equations for which our current level of understanding is at or very close to the critical regularity level; there are other equations (e.g. Benjamin-Ono, Einstein, Zakharov, Kadomtsev-Petviashvili, etc.) for which there are additional obstructions which seem to prevent us from getting close to a critical theory, and we will not discuss these here.

The analytic theory associated to each of the equations is extensive, and we will not be able to even begin to survey all of the developments for each of the equations in this paper, focusing instead only on some representative recent results.  In this particular section we shall concentrate instead on the more algebraic features of these equations, such as the conservation laws, symmetries (especially scaling symmetry), soliton-like solutions, and exact embeddings (or asymptotic embeddings) from one equation to another.

\subsection{Spacetime geometry}

The model equations are intimately tied to the geometry of the underlying spacetime domain, and in some cases also to
the geometry of the target (which is a manifold for the wave maps equation, or a vector bundle for the Maxwell-Klein-Gordon or Yang-Mills equations).  For simplicitly we are considering spacetimes which are completely flat and scale-invariant, but it is still important to note some key geometric features of these spacetimes.

\begin{definition}[Spacetime conventions]  We use $\R^{1+d}$ to denote \emph{Minkowski spacetime}, i.e. the points
$(t,x_1,\ldots,x_d)$ endowed with the Minkowski metric $ds^2 = -c^2 dt^2 + dx_1^2 + \ldots + dx_d^2$, where $c > 0$ is the speed of light (we shall usually normalise $c=1$).  We also write $x_0$ for $x_d$ and $\partial_\alpha$ for $\frac{\partial}{\partial x_\alpha}$.  We use Roman indices $i,j,k$ to sum from $1,\ldots,d$, and Greek indices $\alpha,\beta,\gamma$ to sum from $0,\ldots,d$.  We use $\nabla = \nabla_x$ to denote the spatial gradient, and $\nabla_{t,x}$ for the spacetime gradient.
We raise and lower Greek indices using the Minkowski metric, thus for instance $\partial^i = \frac{\partial}{\partial x_i}$ but $\partial^0 = - c^2 \frac{\partial}{\partial t}$.  Repeated indices will be implicitly summed as per usual, thus for instance the d'Lambertian operator $\Box := \partial^\alpha \partial_\alpha$ can be written in co-ordinates as
$$ \Box = \partial^\alpha \partial_\alpha = \partial^0 \partial_0 + \ldots + \partial^d \partial_d = - c^2 \partial_t^2 + \Delta$$
where $\Delta = \partial_i \partial_i$ is the spatial Laplacian.  We use $\R \times \R^d$ to denote \emph{Galilean spacetime}, which as a set is identical to Minkowski spacetime, but without the Minkowski metric\footnote{There is a natural pseudometric that one should place on Galilean spacetime, which in some sense is the limit of the Minkowski metrics $-c^2 dt^2 + dx_1^2 + \ldots + dx_d^2$ as $c \to \infty$, but defining the pseudometric structure rigorously is somewhat tedious.  Since Galilean spacetime is the only pseudometric space which we will ever consider here, we shall not detail this structure here, though we do remark that this pseudometric can be used to justify the terminology ``pseudoconformal'' which appears later.  Much later on we will also encounter \emph{parabolic spacetime} $\R^+ \times \R^d$, which is the natural spacetime for handling parabolic equations.}; thus with these spacetimes we do not use Greek indices or raising and lowering operations.
\end{definition}

Both Minkowski and Galilean spacetimes enjoy the symmetries of spatial (Euclidean) rotations and reflections, 
spatial translation, time translation, and time reversal.  Minkowski space also enjoys the additional scaling symmetry
$(t,x) \mapsto (\lambda t, \lambda x)$ and the Lorentz boosts
$$ (t,x) \mapsto (\frac{t + v \cdot x/c^2}{\sqrt{1-|v|^2/c^2}}, x_{v^\perp} + \frac{x_v + v t}{\sqrt{1-|v|^2/c^2}} )$$
for any velocity vector $v \in \R^d$ with $|v| < c$, where $x_v$ is the orthogonal projection to the space spanned by $v$, and $x_{v^\perp} := x - x_v$ is the projection to the space orthogonal to $v$.  Meanwhile, Galilean spacetime enjoys a two-parameter scaling symmetry $(t,x) \mapsto (\lambda' t, \lambda x)$ and a Galilean invariance
$$ (t,x) \mapsto (t, x + vt)$$
which is the limit of the Lorentz invariance in the \emph{nonrelativistic limit} $c \to \infty$.  Many of these symmetries will 
be reflected in the model equations; one reason for this is that many of these equations have Lagrangian formulations where the Lagrangian can be defined purely in terms of the geometry of the domain and range and so are automatically invariant (or covariant, in the case of non-scalar equations) under all the symmetries of the underlying geometry.

\subsection{The nonlinear wave equation}  

Let $d \geq 1$, and consider the nonlinear wave equation (NLW)
\begin{equation}\label{nlw}
\Box u = \mu |u|^{p-1} u
\end{equation}
where $u: \R^{1+d} \to \C$ is a complex scalar field, $p > 1$ is the power of the nonlinearity, and $\mu = \pm 1$ is the sign of the nonlinearity (the case $\mu=+1$ is \emph{defocusing}, while the case $\mu=-1$ is \emph{focusing}).  One often restricts attention to the case when $u$ is real-valued, though most of the analysis extends without difficulty to the complex case also.
This equation is also the Euler-Lagrange equation for the functional
$$ \int_{\R^{1+d}} \frac{1}{2} \partial^\alpha \overline{u} \partial_\alpha u + \mu |u|^{p+1}\ dx dt$$
and is thus one of the simplest nonlinear Lagrangian perturbations of the free wave equation (which has the same Lagrangian but with $\mu=0$).  They also appear as special cases of more geometric equations such as wave maps (see below).

Henceforth we normalise $c=1$.  The equation \eqref{nlw} has a conserved \emph{energy}
\begin{equation}\label{nlw-energy}
 E(u) = E(u[t]) := \int_{\R^d} \frac{1}{2} |\partial_t u(t,x)|^2 + \frac{1}{2} |\nabla u(t,x)|^2 + \mu \frac{1}{p+1} |u(t,x)|^{p+1}\ dx.
\end{equation}
Here we adopt the useful convention that $u[t] := (u(t),\partial_t u(t))$ denotes the instantaneous state (both position and velocity) of the field $u$ at time $t$.
Indeed, one can easily verify from differentiating under the integral sign that $E(u[t])$ is independent of $t$ if $u$ is a sufficiently smooth and rapidly decreasing solution to \eqref{nlw}; one can also observe that this energy is the Hamiltonian for \eqref{nlw} using the symplectic structure $\{ (u,u_t), (v,v_t) \} := \Re \int_{\R^d} \overline{u} v_t - v \overline{u}_t\ dx$.  
Observe that in the defocusing case $\mu=+1$ the nonlinear component $\mu \frac{1}{p+1} |u|^{p+1}$ of the energy density has the same sign as the linear component $\frac{1}{2} |u_t|^2 + \frac{1}{2} |\nabla u|^2$, whereas in the focusing case these components have opposing signs.  Thus in the defocusing case we heuristically expect the nonlinearity to amplify the dispersive effects of the linear equation, while in the focusing case we expect the nonlinearity to oppose this dispersion.

The equation \eqref{nlw} also enjoys the scaling invariance
\begin{equation}\label{nlw-scaling}
u(t,x) \mapsto \frac{1}{\lambda^{2/(p-1)}} u(\frac{t}{\lambda}, \frac{x}{\lambda}).
\end{equation}

In the \emph{energy-critical} case $d \geq 3$, $p = 1+\frac{4}{d-2}$, the scaling \eqref{nlw-scaling} preserves the energy \eqref{nlw-energy}.  Note also that in this case the exponent $\frac{2d}{d-2}$ appearing in the nonlinear component of the energy \eqref{nlw-energy} is precisely the exponent appearing in the endpoint Sobolev inequality
$$ \| f \|_{L^{2d/(d-2)}_x(\R^d)} \leq C_d \| \nabla f \|_{L^2_x(\R^d)}.$$
Historically, the energy-critical wave equation was one of the first critical nonlinear evolution equations to have a satisfactory global theory.  This is due to a number of factors, including the finite speed of propagation property (which allows one to analyse blowup by localisation in space), as well as the fact that the conserved momentum
$$ \vec p(u) = \vec p(u[t]) := -\Re \int_{\R^d} \overline{u}_t(t,x) \nabla u(t,x)\ dx$$
(which will ultimately be the source for a key monotonicity formula in the defocusing case) has the same scaling as the conserved energy.

In the focusing case $\mu=-1$ we have the \emph{stationary solutions} $u(t,x) = Q_\omega(x) e^{i\omega t}$, where $\omega > 0$ is a time-frequency and $Q_\omega$ solves the elliptic equation
$$ \Delta Q_\omega + |Q_\omega|^{p-1} Q_\omega = \omega^2 Q_\omega.$$
One can also create travelling wave solutions by applying Lorentz transforms to the stationary solution.  When $Q_\omega$ is a ground state (i.e. it is positive), then these solutions are believed to mark the transition between linear behaviour (such as decay in time) and nonlinear behaviour (such as blowup, or at least lack of decay in time); very recently there has been some progress in making this behaviour rigorous.  One also expects these stationary solutions to play a prominent role in analysis of singularities (blowup) of solutions, though the precise relationship here is presently rather poorly understood.

When $d \leq 2$, or when $d \geq 3$ and $p < 1+\frac{4}{d-2}$, the equation \eqref{nlw} is \emph{energy-subcritical}, because
the scaling \eqref{nlw-scaling} for $\lambda > 1$ will decrease the energy rather than preserve it.  Thus a bounded amount of energy at fine scales is equivalent (after scaling) to a small amount of energy at unit scales, and so we therefore expect the fine-scale behaviour of bounded-energy solutions to be close to linear.  Because of this, the local theory of subcritical equations is very well understood, though the global asymptotic behavior remains a mystery.  

There are a number of other important exponents $p$, such as the \emph{conformal power} $p = 1 + \frac{4}{d-1}$, which makes the equation \eqref{nlw} invariant under conformal transformations of spacetime, and in particular under the Kelvin inversion
$$ u(t,x) \mapsto (c^2 t^2-|x|^2)^{-(d-1)/2} u( \frac{t}{c^2 t^2-|x|^2}, \frac{x}{c^2 t^2-|x|^2} ).$$
With this power the equation is energy-subcritical, though the symplectic structure is now critical.  We will however not discuss this equation in this survey (focusing instead on equations with a 
critical conserved quantity which is positive definite to top order).

\subsection{The nonlinear Schr\"odinger equation}\label{nls-eq}

Take $d \geq 1$ and consider the energy-critical nonlinear Schr\"odinger equation (NLS)\footnote{It is sometimes convenient to replace the linear part $i\partial_t + \Delta$ of this operator with $-i\partial_t + \Delta$, $i \partial_t + \frac{1}{2} \Delta$, or $-i \partial_t + \frac{1}{2} \Delta$ to make certain formulae slightly prettier, however it is a trivial matter to transform one equation to the other (by conjugating, dilating, or stretching the solution $u$ in space or time) and so all choices of operator here are essentially equivalent.}
\begin{equation}\label{nls}
i u_t + \Delta u = \mu |u|^{p-1} u
\end{equation}
where $u: \R \times \R^d \to \C$ is a complex scalar field, and $\mu = \pm 1$ is the sign of the nonlinearity (again, $\mu=+1$ is \emph{defocusing}, while the case $\mu=-1$ is \emph{focusing}).   These equations arise naturally as models describing various forms of weakly dispersive behaviour; see \cite{sulem} (as well as the discussion on the gKdV equation below).  
The case $d=1,p=3$ happens to be completely integrable, but in general the equations are merely Hamiltonian (though they do enjoy a large, but finite, number of conserved quantities).

The scaling symmetry is now given by
\begin{equation}\label{nls-scaling}
u(t,x) \mapsto \frac{1}{\lambda^{2/(p-1)}} u(\frac{t}{\lambda^2}, \frac{x}{\lambda})
\end{equation}
while the conserved energy is now
\begin{equation}\label{nls-energy}
 E(u) = E(u(t)) := \int_{\R^d} \frac{1}{2} |\nabla u(t,x)|^2 + \mu \frac{1}{p+1} |u|^{p+1}(t,x)\ dx.
\end{equation}
Again, this energy can be interpreted as a Hamiltonian for \eqref{nls}, using the symplectic form $\{u,v\} = \int_{\R^d} \Im( \overline{u} v )\ dx$.  The NLS also has an additional phase rotation symmetry $u(t,x) \mapsto e^{i\theta} u(t,x)$, which leads (via Noether's theorem) to a second important conserved quantity\footnote{The analogue of this quantity for NLW would be the charge $\int \Im( \overline{u} u_t )\ dx$, but this quantity vanishes for the most important case of real scalar fields $u$ and so has not been of major importance in the analysis.}, the \emph{mass} (or \emph{charge})
\begin{equation}\label{nls-mass}
M(u) = M(u(t)) = \int_{\R^d} |u(t,x)|^2\ dx.
\end{equation}
The translation symmetry $u(t,x) \mapsto u(t-x_0)$ also leads to a third conserved quantity, the \emph{momentum}
\begin{equation}\label{nls-momentum}
 \vec p(u) = \vec p(u(t)) := 2\int_{\R^d} \Im( \overline{u(t,x)} \nabla u(t,x))\ dx.
\end{equation}

When $d \geq 3$ and $p = 1+\frac{4}{d-2}$, the equation \eqref{nls} is energy-critical but mass-supercritical and momentum-supercritical; conversely, in the \emph{pseudoconformal} case $p = \frac{4}{d}$ the equation \eqref{nls} is mass-critical but energy-subcritical and momentum-subcritical.   Thus in both cases, the momentum (which supplies a crucial monotonicity formula in the large data theory) is not scale-invariant, which causes significant technical difficulties in the analysis.

Of the two critical equations, the mass-critical equation is considered harder to analyse.  This is because in this case
the NLS equation enjoys two less obvious symmetries, namely the \emph{Galilean invariance}
$$ u(t,x) \mapsto e^{-i t |v|^2/4} e^{i v \cdot x/2} u(t, x-vt)$$
where $v \in \R^d$ is arbitrary\footnote{Indeed, this invariance holds for all powers $p$, being the analogue of the Lorentz invariance for the NLW.  The pseudoconformal symmetry however is restricted to the pseudoconformal exponent $p = 1 + \frac{4}{d}$.}, as well as the \emph{pseudoconformal symmetry}
\begin{equation}\label{pc-conf}
 u(t,x) \mapsto \frac{1}{|t|^{d/2}} e^{i|x|^2/4t} u( \frac{1}{t}, \frac{x}{t} ) 
 \end{equation}
for $t \neq 0$.  These two symmetries (as well as spatial translation symmetry) also preserve the mass \eqref{nls-mass}, thus the mass is in fact critical with respect to quite a large group of symmetries.  This wealth of symmetries complicates the analysis, because it implies quite a serious breakdown of compactness for the ``essential'' part of the dynamics.  (The Galilean invariance is not a serious issue for the energy-critical equation, basically because it does not leave the energy invariant.)

As with NLW, the focusing NLS ($\mu=-1$) also enjoys stationary solutions (or \emph{solitons}) $u(t,x) = Q_\omega(x) e^{i\omega t}$, where $\omega > 0$ is a time-frequency and $Q_\omega$ solves the elliptic equation
$$ \Delta Q_\omega + |Q_\omega|^{p-1} Q_\omega = \omega Q_\omega.$$
One can apply Galilean invariance to also obtain travelling soliton solutions.  As with NLW, the ground state solitons are expected to demarcate the transition between linear and nonlinear behaviour, and to dominate the dynamics of blowup (at least in certain cases), and there are now several rigorous results that demonstrate this fact.  

There is an algebraic embedding of NLS into NLW: if $u: \R \times \R^d \to \C$ solves \eqref{nls} in $d$ spatial dimensions, then the complex field $\tilde u: \R^{1+(d+1)} \to \C$ defined by
$$ \tilde u(t,x_1,\ldots,x_{d+1}) := e^{i(t+x_{d+1})} u(t-x_{d+1}, x_1,\ldots,x_d)$$
solves \eqref{nlw} in $d+1$ spatial dimensions (with $c=1$); in Fourier space, this fact becomes the geometric observation that a $d$-dimensional paraboloid can be viewed as a section of a $d+1$-dimensional cone.  This allows one to deduce many 
algebraic identities for the $d$-dimensional NLS from the corresponding identities for the $d+1$-dimensional NLW (the ``method of descent'').  However, this embedding of NLS into NLW, while exact, is not very useful analytically as it maps finite-energy solutions to infinite-energy ones.  There is a more profitable \emph{asymptotic} embedding from NLS to a variant of NLW, the \emph{nonlinear Klein-Gordon equation} (NLKG)
$$ \Box u = c^4 u + \mu |u|^{p-1} u,$$
namely that if $u: \R \times \R^d$ solves NLS, then the complex field $\tilde u: \R^{1+d} \to \C$ defined by
$$ \tilde u(t,x) := e^{-ic^2 t} u( t / 2c^2, x )$$
solves NLKG up to errors which are $O(c^{-4})$.  We will however not discuss the NLKG here (it is not scale-invariant and so the study of this equation at critical regularities becomes messier).

\subsection{The generalised Korteweg-de Vries equation}

Take $d=1$, and consider the \emph{generalised Korteweg-de Vries} (gKdV) equation\footnote{This family of equations should not be confused with the \emph{Korteweg-de Vries hierarchy} or the \emph{modified Korteweg-de Vries hierarchy}, which are a commuting sequence of completely integrable equations starting from KdV or mKdV which are of increasingly high order (involving more and more spatial derivatives) as one proceeds up the hierarchy.}
\begin{equation}\label{gkdv}
u_t + u_{xxx} = \mu (u^p)_x
\end{equation}
where $u: \R \times \R \to \R$ is a real scalar field, $p \geq 2$ is an integer, and $\mu = \pm 1$ is a sign.  When $p$ is even the sign of $\mu$ is irrelevant (as one can remove it via the change of variables $u \mapsto -u$); but when $p$ is odd we make a distinction between the defocusing case $\mu=+1$ and the focusing case $\mu=-1$.  The case $p=2$ is known as the \emph{Korteweg-de Vries} (KdV) equation, while the case $p=3$ is the \emph{modified Korteweg-de Vries} (mKdV) equation, which
are both well-known examples of completely integrable systems.  The higher values of $p$ are not completely integrable.  These equations can arise as dispersive models for the evolution of one-dimensional water waves in shallow canals.

The gKdV equations are somewhat similar to the one-dimensional NLS equations with the same values of $\mu$ and $p$ (especially when $p$ is odd).
One evidence of this similarity can be seen the conserved mass and energy for gKdV,
\begin{align*}
M(u) = M(u(t)) &:= \int_\R u(t,x)^2\ dx \\
E(u) = E(u(t)) &:= \int_\R \frac{1}{2} u_x(t,x)^2 + \mu \frac{1}{p+1} |u(t,x)|^{p+1}\ dx
\end{align*}
and the scaling symmetry
$$ u(t,x) \mapsto \frac{1}{\lambda^{2/p-1}} u(\frac{t}{\lambda^3}, \frac{x}{\lambda}).$$
The energy is once again the Hamiltonian for the flow, but now using a slightly different symplectic form, $\{u,v\} := \int_\R u \partial^{-1}_x v\ dx$.  On the other hand, in contrast to NLS, the gKdV equation is not Galilean-invariant, although in the limiting case of very coherent wave trains with almost constant frequency, the envelope of these trains does behave in a Galilean-invariant manner and indeed is asymptotically modeled by NLS; more precisely, if $u: \R \times \R \to \C$ solves NLS with $d=1$ and $p$ an odd integer, then the field $\tilde u_N: \R \times \R \to \R$ defined for a large frequency parameter $N \gg 1$ by
$$ \tilde u_N := \left( \frac{2^{p-1}}{N \binom{p}{(p-1)/2}} \right)^{1/(p-1)} \Re\left( e^{iN x} e^{iN^3 t} \overline{u}(t, \frac{x+3N^2 t}{3^{1/2} N^{1/2}} ) \right)$$
solves gKdV up to errors which are small (or at least ``non-resonant'') in the limit $N \to \infty$; see \cite{cct}, \cite{tao-gkdv} for some applications of this asymptotic embedding of NLS in gKdV.

When $\mu = -1$, the gKdV equation admits traveling wave (soliton) solutions $u(t,x) = Q_v(x-vt)$, where $v > 0$ is a rightward velocity and $Q_v$ solves the ground state equation
$$ \Delta Q_v + |Q_v|^{p-1} Q_v = v Q_v.$$
Once again, we expect these solitons to mark the transition between linear and nonlinear behavior, and to be involved in the mechanism for blowup, and we have a certain number of results in these directions, especially concerning small perturbations of the ground state (or vacuum state).

The energy for gKdV is always supercritical.  The mass is subcritical for $p < 5$, critical for $p=5$, and supercritical for $p>5$.  One complication in this equation compared to the NLS is that there is no exact Galilean invariance, and no conserved momentum; nevertheless, one still has the same type of failure of compactness that one would normally associate with this invariance.  On the other hand, this equation has a useful decoupling property, in that radiative components of the solution tend to propagate to the left, while soliton-type components of the solution tend to propagate to the right.  The derivative in the nonlinear term in \eqref{gkdv} causes some difficulty, though these are largely compensated for by the strong dispersive and local smoothing properties of the linear counterpart of the gKdV equation, namely the \emph{Airy equation} $u_t + u_{xxx} = 0$.

The KdV equation (with the normalisation $\mu=3, p=2$) and the defocusing mKdV equation (with the normalisation $\mu=2, p=3$) are connected by the remarkable \emph{Miura transform}: if $u$ solves mKdV, then $u_x + u^2$ solves KdV.  This transform is almost a bijection between $H^s$ and $H^{s-1}$ for various values of $s$, which has allowed one to derive analytical results for one equation via analytical results (at one higher or lower derivative of regularity) for the other.  We will however not discuss these types of results here, focusing instead on the scale-invariant theory (which for a number of reasons is not currently available either for KdV or for mKdV).

\subsection{The wave maps equation}

We now move from the scalar field models to the geometric model nonlinear wave equations, which we shall also refer to as \emph{systems} to emphasise their non-scalar nature.  These systems are often significantly more nonlinear in nature, but to compensate for this they have an extremely geometric structure which can be exploited (e.g. via gauge symmetries) to renormalise the equation.

Let $d \geq 1$, let $M = (M,g)$ be an $m$-dimensional Riemannian manifold with Levi-Civita connection $\nabla$, which acts on smooth sections of the tangent bundle $TM$.  If $\phi: \R^{1+d} \to M$ is a smooth map, then we obtain the pullback $\phi^* \nabla$, which acts on smooth sections of the pullback bundle $\phi^*(TM)$.  We say that $\phi$ is a \emph{wave map} if we have
$$ (\phi^* \nabla)^\alpha \partial_\alpha \phi = 0$$
where we again use the usual raising and lowering conventions; this is the Euler-Lagrange equation for the functional
$$ \int_{\R^{1+d}} \langle \partial^\alpha \phi(t,x), \partial_\alpha \phi(t,x) \rangle_g\ dx dt$$
and is thus the natural Lagrangian generalisation of the free wave equation to fields that take values in Riemannian manifolds.  This equation is also the natural hyperbolic generalisation of harmonic maps (or of the parabolic counterpart, the harmonic map heat flow), and also is a simplified model for studying certain symmetric cases of the Einstein equations of general relativity.

If we parameterise $M$ by local coordinates, thus $\phi = \phi^i$
for $i=1,\ldots,m$, then we can recast the wave maps equation as a nonlinear wave equation
$$ \Box \phi^i = - \Gamma(\phi)^i_{jk} \partial^\alpha \phi^j \partial_\alpha \phi^k$$ 
where $\Gamma$ is the Christoffel symbol.  If $M$ is the unit sphere $S^m \subset \R^{m+1}$, so that $\phi$ can be viewed as taking values in the Euclidean space $\R^{m+1}$ subject to the constraint $\langle \phi, \phi \rangle_{\R^{m+1}} = 1$, then the wave maps equation becomes
$$ \Box \phi = - \phi \langle \partial^\alpha \phi, \partial_\alpha \phi \rangle_{\R^{m+1}}$$
which can be viewed as a ``defocusing'' case of the wave maps equation,
whereas if $M$ is the hyperbolic space $H^m \subset \R^{1+m}$, which can be thought of as the upper unit sphere
$H^m = \{ (t,x) \in \R^{1+m}: t = \sqrt{1+|x|^2} \}$ of Minkowski space $\R^{1+m}$, then the wave maps equation becomes
$$ \Box \phi = \phi \langle \partial^\alpha \phi, \partial_\alpha \phi \rangle_{\R^{1+m}}$$
which can be viewed as a ``focusing'' case of the equation.  Note in all cases the wave maps equation takes the schematic form
$$ \Box \phi = O( F(\phi) \partial \phi \partial \phi )$$
for some specific function $F()$.  In particular, the nonlinearity contains first derivatives of $\phi$, which creates significant new technical difficulties (not present in simpler models such as NLW) when trying to control the nonlinear terms by perturbative methods.

Now we set $c=1$.  The wave maps equation has a scale invariance
$$ \phi(t,x) \mapsto \phi(\frac{t}{\lambda}, \frac{x}{\lambda})$$
and so the natural scale-invariant norm to analyse this data would be the homogeneous Sobolev norm
$$ \| \phi(t) \|_{\dot H^{d/2}_x(\R^d)} + \| \phi_t(t) \|_{\dot H^{d/2-1}_x(\R^d)} $$
ignoring for now the delicate issue of how to properly define this norm for fields taking values in a manifold $M$.  Comparing this against
the conserved energy
$$ E(\phi) = E(\phi[t]) = \int_{\R^d} \frac{1}{2} |\partial_t \phi(t,x)|_g^2 + \frac{1}{2} |\nabla \phi(t,x)|_g^2 \ dx$$
of the equation, we see that the energy is subcritical in one dimension $d=1$, critical in two dimensions $d=2$, and supercritical in higher dimensions.  Unlike NLW, the distinction between focusing and defocusing wave maps is not immediately apparent from the energy density, but can be seen from a number of more subtle considerations, such as the embedding of NLW in WM discussed below.

The current tools used to analyse solutions of nonlinear PDE, such as the Fourier transform, are well adapted to scalar fields but are not as suitable for more complicated fields, such as the field $\phi$, as they are sensitive to the choice of co-ordinates used.  Indeed, selecting good coordinates on $M$ (or on the pullback tangent bundle $\phi^* TM$) is a key step in obtaining a satisfactory critical-regularity analysis.

The analogue of solitons for the WM equation are the harmonic maps (and their Lorentz boosts).  One reason why the negative curvature case is considered defocusing (and thus easier to study) is because such target manifolds cannot support any non-trivial finite energy harmonic maps (thanks to the Bochner identity); heuristically, this should thus prevent the wave map equation from blowing up in finite time, though it turns out that in the supercritical case $d > 2$ that blowup can still occur.
In the focusing case, harmonic maps played a key role in the recent establishment of blowup in the critical case $d=2$. In contrast, in the defocusing case it is conjectured (and widely believed) that no blowup occurs.  It seems that harmonic maps in fact play a decisive role in the blowup and asymptotics of the wave map equation, but the situation is certainly far from understood at present (except when one imposes strong symmetry assumptions on the initial data).

There is a connection between $U(1)$-equivariant energy-critical wave maps, and (spherically symmetric) energy-critical NLW.  For instance, if $M$ is the surface $\{ (s,\alpha): \R^+ \times \R/2\pi \Z: 1 + \frac{\mu}{2} s^2 > 0 \}$ with the metric
$ds^2 + (s^2 + \frac{\mu}{2} s^4) d\alpha^2$, and $\phi: \R^{1+2} \to M$ is an \emph{equivariant} map in the sense that
$$ \phi( t, r \cos \theta, r \sin \theta ) = (r u(t, r), \theta )$$
for all $r \geq 0$, $t \in \R$, and $\theta \in \R$, and some $u: \R \times \R^+ \to \R$ then one can verify (assuming that $\phi$ avoids the singularity $1 + \frac{\mu}{2} s^2 = 0$, which only occurs in the focusing case $\mu=-1$) that the spherically symmetric field $u: \R^{1+4} \to \R$ defined by $u(t,x) := u(t,|x|)$ solves the energy-critical NLW \eqref{nlw} with $d=4$ and $p=3$. Note that $M$ has negative curvature when $\mu=+1$ and positive curvature when $\mu=-1$, thus reinforcing the analogy between negative (resp. positive) curvature and defocusing (resp. focusing) nonlinear equations.

\subsection{Schr\"odinger maps}

Schr\"odinger maps are the analogue of wave maps, but where the linear operator underlying the evolution is the Schr\"odinger operator $i\partial_t + \Delta$ rather than the d'Lambertian $\Box$.  (Similarly, harmonic maps and the harmonic map heat flow have the Laplacian $\Delta$ and the heat operator $\partial_t + \Delta$ respectively as the underlying linear operator.)  The geometric setup is the same as that for wave maps, except that the domain is now Galilean spacetime $\R \times \R^d$ instead of Minkowski spacetime $\R^{1+d}$ and that the manifold $M$ is not just a Riemannian manifold, but is in fact a K\"ahler manifold.
In particular, the tangent bundle $TM$ has a complex structure $z \mapsto iz$.  A map $\phi: \R \times \R^d \to M$ is then said to be a \emph{Schr\"odinger map} (SM) if it obeys the equation
$$ i \partial_t \phi + (\phi^* \nabla)_j \partial_j \phi = 0.$$
In coordinates, the SM equation takes the schematic form
$$ i \partial_t \phi + \Delta \phi = O( F(\phi) \partial \phi \partial \phi )$$
for some function $F(\phi)$ depending on the manifold $M$ (and the coordinate system chosen).  While very similar in form to the
wave maps equation, the derivatives in the nonlinearity are significantly harder to handle here, because the linear operator
$i \partial_t + \Delta$, being only first order in time, has more difficulty compensating for (or ``recovering'') the loss of derivative in the nonlinearity than the linear operator $\Box = -\partial_t^2 + \Delta$, which is second order in time.  Thus while the geometry and algebraic structure of the SM equation is very similar to that of the WM equation, the analysis is significantly more technical.

For simplicity let us restrict attention to the case when the target manifold $M$ is the Riemann sphere $S^2$; this has positive curvature and should thus be viewed as a ``focusing'' case.  If we embed $S^2$ in the Euclidean space $\R^3$, thus viewing $\phi$ as a map from $\R^{1+d}$ to $\R^3$ with $\langle \phi, \phi \rangle_{\R^3} = 1$, then the equation becomes
$$ \partial_t \phi = \phi \times \Delta \phi$$
where $\times$ is the cross product on $\R^3$.  This is not obviously a nonlinear Schr\"odinger equation.  If however we place complex coordinates on the sphere, for instance by using the stereographic projection
$$ \left( \frac{2 \Re(z)}{1+|z|^2}, \frac{2\Im(z)}{1+|z|^2}, \frac{1-|z|^2}{1+|z|^2} \right) \mapsto z$$
(ignoring for now the issue of the singularity at the north pole $(0,0,1)$) to identify $S^2$ with the complex plane $\C$ with
the metric $\frac{4}{(1+|z|^2)^2} |dz|^2$, then the equation becomes
$$ i \partial_t z - \Delta z = \frac{2 \overline{z}}{1 + |z|^2} \partial_j z \partial_j z.$$

The Schr\"odinger maps equation has the scale invariance
$$ \phi(t,x) \mapsto \phi(\frac{t}{\lambda^2}, \frac{x}{\lambda})$$
and so the natural scale-invariant norm to analyse this data would be the homogeneous Sobolev norm $\dot H^{d/2}_x(\R^d)$.
Comparing this against the conserved energy
$$ E(\phi) = E(\phi[t]) = \int_{\R^d} \frac{1}{2} |\nabla \phi(t,x)|_g^2 \ dx$$
we see (as with WM) that the energy is subcritical in one dimension $d=1$, critical in two dimensions $d=2$, and supercritical in higher dimensions.

As with wave maps, harmonic maps are the natural analogue of the soliton solutions for the SM equation.  However, at present we have virtually no understanding of the role these stationary solutions play in the evolution.  Nevertheless, there has been some extremely recent progress towards a global critical theory for these equations, and while the results here lag somewhat the analogous results for wave maps, it seems reasonable to expect parity in these theories in the long term.

\subsection{The Maxwell-Klein-Gordon system}

After the wave maps equation, the next most complicated field equation is the \emph{Maxwell-Klein-Gordon} (MKG) system, which is a coupled system of a section $\phi$ of a complex line bundle on $\R^{1+d}$ and a $U(1)$ connection $D$ on this bundle, being the Euler-Lagrange equation for the Lagrangian
$$ \int_{\R^{1+d}} \frac{1}{2} \langle D^\alpha \phi, D_\alpha \phi \rangle + \frac{1}{4} \langle F^{\alpha \beta}, F_{\alpha \beta} \rangle\ dx dt$$
where $F^{\alpha \beta} = [D^\alpha,D^\beta]$ is the curvature of the connection.  Physically, $\phi$ represents a charged particle field, while $D$ represents the electromagnetic field which is both generated by and drives the particle field.  If one removes the particle field $\phi$, one obtains the (linear) Maxwell equations, while if one instead removes the electromagnetic field $D$ then one obtains the free wave equation.  The nonlinear effects of the MKG system thus arise solely from interactions between the two fields.

We can recast the MKG system in coordinates by choosing a trivialisation $\R^{1+d} \times \C$ of the complex line bundle, thus $\phi: \R^{1+d} \to \C$ now is interpreted as a complex scalar field, and $D_\alpha = \partial_\alpha + i A_\alpha$ for some real one-form $A_\alpha: \R^{1+d} \to \R$.  We then have $F_{\alpha \beta} = i (\partial_\alpha A_\beta - \partial_\beta A_\alpha)$, and the Maxwell-Klein-Gordon system can be written as
\begin{align*}
\partial^\beta F_{\alpha \beta} &= i \Im( \phi \overline{D_\alpha \phi} )\\
D_\alpha D^\alpha \phi = 0.
\end{align*}
The second equation can be regarded as a nonlinear equation for $\phi$, which schematically has the form
$$ \Box \phi = O( A \partial \phi ) + O( \partial A \phi ) + O( A^2 \phi ).$$
The first equation can be viewed as partially describing an evolution for the connection $A$, but it is underdetermined (roughly speaking, it only specifies the curl of $A$ but not the divergence).  This is ultimately due to the fact that there are many possible trivialisations of the complex line bundle, each leading to essentially the same field, and that the evolution should really be quotiented out by the action of the \emph{gauge symmetry}
$$ (\phi, A_\alpha) \mapsto (e^{i\chi} \phi, A_\alpha - \partial_\alpha \chi )$$
for any smooth gauge function $\chi: \R^{1+d} \to \R$.  Ideally, all of the analytical tools used to study this equation should be invariant under this gauge invariance.  This turns out however to be impractical (at least with current technology), and instead one selects a gauge for this equation in order to make the evolution determined, and also as ``linear'' as possible, in order to maximise the effectiveness of the analytical tools.  A particularly popular gauge for this equation is the \emph{Coulomb gauge} $\div A = 0$.  This turns the equation for $A$ into something schematically resembling
$$ \Box A = O( \phi \partial \phi ) + O( A \phi^2 ).$$
Thus we see that we obtain a system of nonlinear wave equations, containing derivatives in the nonlinearity.

We again set $c=1$.  The Maxwell-Klein-Gordon system enjoys the scaling symmetry
$$ (\phi(t,x), A_\alpha(t,x)) \mapsto (\frac{1}{\lambda} \phi(\frac{t}{\lambda}, \frac{x}{\lambda}), 
\frac{1}{\lambda} A_\alpha(\frac{t}{\lambda}, \frac{x}{\lambda})) $$
and the conserved energy
$$ E(\phi,A) = E(\phi[t],A[t]) := \int_{\R^d} \frac{1}{2} |F_{0i}(t,x)|^2 + \frac{1}{2} |F_{ij}(t,x)|^2 + 
\frac{1}{2} |D_0 \phi(t,x)|^2 + \frac{1}{2} |D_i \phi(t,x)|^2\ dx$$
where the Roman indices $i,j$ are implicitly summed from $1$ to $d$.  One can then easily verify that the equation is energy-subcritical in three and fewer dimensions, energy-critical in four dimensions, and energy-supercritical in five and higher dimensions.

Although not apparent at first glance, the Maxwell-Klein-Gordon equation has many similarities with the wave maps equation, especially if the target manifold of the latter is a Riemann surface.  Then both equations can be rewritten as a $U(1)$-covariant wave equation, where the $U(1)$ connection itself obeys some differential equation.  However, a key difference is that in wave maps the connection obeys (after suitable gauge fixing)
an elliptic equation which makes the connection close to flat, whereas in Maxwell-Klein-Gordon the connection itself evolves by a nonlinear wave equation.  For the critical regularity global theory, one is then forced to develop more ``covariant'' techniques, in which one exploits the dispersive properties of covariant wave equations rather than free wave equations.
Also, the MKG equation is not considered to be either focusing nor defocusing; the nonlinear effects do not have a preferred sign.

\subsection{The Yang-Mills equation}

The (hyperbolic) \emph{Yang-Mills} (YM) equation is the time-dependent analogue of the more well-known elliptic Yang-Mills equation, which plays an important role in physics, geometry, and integrable systems.  Informally, the hyperbolic Yang-Mills equation describes the free evolution of a connection, just as the wave maps equation describes the free evolution of an immersed surface.  It is closely related to the Maxwell-Klein-Gordon equation; it does not have the scalar field $\phi$, but to compensate for this the connection $D$ now acts on a vector bundle with a \emph{nonabelian} gauge group, thus re-introducing nonlinearity back into the system.  (One can simultaneously generalise the NLW, MKG, and YM by considering the \emph{Yang-Mills-Higgs equation}, but we will not discuss this more complicated system here.)  

More formally, given a vector bundle\footnote{One can of course define Yang-Mills connections on other $G$-bundles, such as principal bundles; the theory is essentially the same.} on Minkowski space $\R^{1+d}$ with the orthonormal action of a compact Lie group $G$ (with Lie algebra $\g$), consider (smooth) connections $D$ on this bundle, and form the curvature $F_{\alpha \beta} = [D_\alpha, D_\beta]$ in the usual manner; one can view $F_{\alpha \beta}$ as an equivariant two-form on the bundle taking values in $\g$, and so in particular the Yang-Mills density $\langle F_{\alpha \beta}(t,x), F^{\alpha \beta}(t,x) \rangle$ is well-defined (here the inner product is the Hilbert-Schmidt inner product).  One then defines $D$ to be a \emph{Yang-Mills connection} if it is a critical point for the Yang-Mills functional
$$ \int_{\R^{1+d}} \langle F_{\alpha \beta}(t,x), F^{\alpha \beta}(t,x) \rangle\ dx dt.$$
In co-ordinates (choosing a trivialisation $\R^{1+d} \times \R^m$ of the vector bundle, and identifying $G$ with a subgroup of the orthogonal group $O(m)$), the connection $D$ (when acting on the original vector bundle) takes the form $D_\alpha = \partial_\alpha + A_\alpha$, where $A$ is a $\g$-valued one-form, and the connection $F_{\alpha \beta}$ is now the $\g$-valued two-form
$$ F_{\alpha \beta} = \partial_\alpha A_\beta - \partial_\beta A_\alpha + [A_\alpha, A_\beta].$$
The Yang-Mills equation is then
$$ D^\alpha F_{\alpha \beta} = 0$$
where the connection $D_\alpha$ acts on $\g$-valued forms $\omega$ by the formula
$$ D_\alpha \omega = \partial_\alpha \omega + [A_\alpha, \omega]$$
and is raised and lowered via the Minkowski metric in the usual manner.  We remark that the curvature $F_{\alpha \beta}$, by definition, also automatically satisfies the \emph{Bianchi identity}
$$ D_{\alpha} F_{\beta \gamma} + D_\beta F_{\gamma \alpha} + D_\gamma F_{\alpha \beta} = 0,$$
thus in some sense the curvatures of Yang-Mills connections are simultaneously ``divergence-free'' and ``curl-free''.

As with the Maxwell-Klein-Gordon equation, the Yang-Mills equation has a gauge symmetry due to the fact that bundles have multiple trivialisations. Indeed, given any smooth map $U: \R^{1+d} \to G$, we have the gauge invariance
$$ A_\alpha \mapsto U A_{\alpha} U^{-1} - (\partial_\alpha U) U^{-1}; \quad D_\alpha \mapsto U D_\alpha U^{-1};
\quad F_{\alpha \beta} \mapsto U F_{\alpha \beta} U^{-1}.$$
Thus we need to fix the gauge (at least partially) before the Yang-Mills system is well-posed.  One possible choice is the \emph{Lorenz gauge} $\partial^\alpha A_\alpha = 0$, which would convert the Yang-Mills equation into a nonlinear wave equation, schematically of the form
$$ \Box A = O( A \partial A ) + O( A^3 ).$$
As it turns out, however, this is not the ideal formulation for this system, and a slight variant of this gauge (the Coulomb gauge) is preferred instead.  Nevertheless, one should still think of the Yang-Mills equations as a type of nonlinear wave equation, whose nonlinearity is similar in strength to that of the Maxwell-Klein-Gordon system.

Now we set $c=1$.  The Yang-Mills equation enjoys the scaling symmetry
$$ A_\alpha(t,x) \mapsto \frac{1}{\lambda} A_\alpha( \frac{t}{\lambda}, \frac{x}{\lambda} ); \quad
F_{\alpha \beta}(t,x) \mapsto \frac{1}{\lambda^2} F_{\alpha \beta}( \frac{t}{\lambda}, \frac{x}{\lambda} )$$
(thus $A$ scales like a first-order derivative, while $F$ scales like a second-order derivative) and also has the conserved energy
$$ E(A) = E(A[t]) := \int_{\R^d} \frac{1}{2} |F_{0i}(t,x)|^2 + \frac{1}{2} |F_{ij}(t,x)|^2\ dx$$
where we sum Roman indices $i,j$ from $1$ to $d$, and the magnitude of $F$ is taken in the Hilbert-Schmidt sense.  As with MKG, the equation is energy-subcritical in three and fewer spatial dimensions, energy-critical in four spatial dimensions, and energy-supercritical in five and higher dimensions.

Progress on the Maxwell-Klein-Gordon and Yang-Mills systems have proceeded more or less in tandem, with the Yang-Mills equations
considered slightly more difficult due to the non-abelian gauge group and due to the less decoupled nature of the nonlinear
interactions (in MKG, the connection $A$ evolves in a nearly linear manner, while the nonlinear effects on the particle field $\phi$ are caused entirely by $A$).  In the most recent progress on these systems, in which gauge theory has played a more prominent role, the non-abelian nature of the gauge group has caused some highly nontrivial technical difficulties for YM that were not present for MKG.  Nevertheless, these two systems of equations are still considered very similar (for instance, they are closer to each other than they are to WM).  

As with MKG, the YM equations are not considered to be either focusing or defocusing.  Nevertheless, they have an important family of stationary solutions, the \emph{instantons} (finite-energy global smooth solutions to the elliptic Yang-Mills equations), which are analogous to the soliton solutions for other models such as NLW, NLS, and gKdV.  Based on this analogy one would expect the instantons to play a role in the large data global theory of YM, but the theory here is virtually non-existent (except for numerics), due to the significant analytical difficulties encountered in trying to obtain a critical theory for the Yang-Mills equation.

\section{The scaling heuristic}

In this section we try to informally motivate the importance of the criticality, sub-criticality, or super-criticality of the conserved quantities in determining whether the evolution is ultimately linear or nonlinear; in the next section we discuss how to make these heuristics rigorous.  To illustrate the principle, we shall work with one of the simplest models, namely the NLS \eqref{nls}, and with a simple conserved quantity, namely the mass.

By restricting the class of initial data $u(0)$ appropriately, one may assume that this initial data is smooth and rapidly decreasing, and thus bounded in all norms.  However, as the evolution progresses, the solution may well grow in many of these norms.  The only norms which we know for certain to be bounded uniformly in time are those given by conserved quantities (or variants of conserved quantities, such as monotone quantities or quantities which are conserved up to lower order errors).  If we know or suspect that the linear behaviour will be dominant for all time, then we also expect to control the solution in all the norms for which we know the linear solution to be bounded.  This type of result can often be established for small data by perturbative and boostrap techniques, and (with much more effort) for large data when the nonlinearity is defocusing.  However, in many cases we cannot assume \emph{a priori} that the linear behaviour is dominant, and so we can only rely on the control on the solution given by the conserved
quantities\footnote{One could also hope to exploit the heuristics of thermodynamics, which predict that for sufficiently complex systems, the evolution should be distributed ``uniformly'' across all areas of phase space which are consistent with the conservation laws, the initial data, and other structures of the equation.  Such uniform distribution results could significantly augment the control on the solution given by the conservation laws alone.  However, for deterministic PDE such as the ones studied here, there have been no rigorous results in this direction with the current level of technology.}.  This naturally leads to the following question: if all we know about the initial data is that its conserved quantities are all bounded, is this enough to determine whether the linear behaviour of the solution dominates the nonlinear behaviour or not?

Of course, we have not rigorously defined what it means for the linear behaviour to ``dominate'' the nonlinear behaviour.  Let us experiment by using a very crude test for this domination.  Write $u_0(x) := u(0,x)$ for the initial data. Rewrite the NLS equation \eqref{nls} at time $t=0$ as
$$ u_t(0,x) = i \Delta u_0(x) - i \mu |u_0(x)|^{p-1} u_0(x),$$
thus the initial time variation $u_t(0,x)$ of the solution has a linear component $i\Delta u_0(x)$ and a nonlinear component $i \mu |u_0(x)|^{p-1} u_0(x)$.  
We shall naively decide that the linear evolution dominates if the initial magnitude $|i \Delta u_0(x)|$ of the linear component exceeds that of the initial nonlinear component $i \mu |u_0(x)|^{p-1} u_0(x)$, or in other words that
$$ |\Delta u_0(x)| \gg |u_0(x)|^p.$$
Of course, if the reverse inequality holds then we shall decide that the nonlinear evolution will dominate.  Note that this crude test is insensitive to the sign $\mu$ of the nonlinearity, as we are ignoring whether the linear and nonlinear components are interfering constructively or destructively.  Also, this test is only inspecting the behaviour at the initial time $t=0$; at late times the solution may be so different from the initial data that the initial comparison is no longer relevant.  As this is only a heuristic discussion, we will not try to address these objections here.

Now suppose we know that the mass of the initial data is equal to some value $M$, thus
$$ \int_{\R^2} |u_0(x)|^2\ dx = M.$$
There are of course infinitely many such data which obey this mass bound.  But let us make some guesses as to which data should provide the ``worst'' or ``most nonlinear'' behaviour.  Typically, the nonlinear effects tend to be strongest when the solution is concentrated all in one place (so that its amplitude is maximised), rather than when it is dispersed in multiple places.  One model for depicting such a concentration is by assuming that $u_0(x)$ is a \emph{rescaled bump function}
$$ u_0(x) := M^{1/2} N^{d/2} \varphi( N x )$$
where $\varphi \in C^\infty_0(\R^d)$ is a bump function, which we normalise to have total mass 
$\int_{\R^d} |\varphi(x)|^2\ dx = 1$.  The factor $M^{1/2} N^{d/2}$ is needed to ensure that the mass of $u_0$ remains at $M$.
Informally, $u_0$ has magnitude $\sim M^{1/2} N^{d/2}$ on a ball of radius $\sim 1/N$; the parameter $N$ then represents the main frequency magnitude of this data, while the inverse parameter $1/N$ represents the spatial scale.  Thus large $N$ corresponds to high frequencies and fine scales, while small $N$ corresponds to low frequencies and coarse scales.

In this rescaled bump function example, the initial linear component magnitude $|\Delta u_0(x)|$ has magnitude $\sim M^{1/2} N^{d/2} N^2$ on a ball of radius $\sim 1/N$, while the initial nonlinear component magnitude $|u_0(x)|^p$ has magnitude $\sim (M^{1/2} N^{d/2})^p$ on the same ball.  Thus we expect the linear behaviour to dominate when
$$ M^{1/2} N^{d/2} N^2 \gg (M^{1/2} N^{d/2})^p$$
which can be rearranged as
\begin{equation}\label{nheu}
 N^{p - (1 + \frac{4}{d})} \ll M^{(p-1)/d}.
\end{equation}
Thus, in the \emph{mass-subcritical} case, when $p - (1 + \frac{4}{d})$ is negative, we thus expect the linear behaviour to dominate for high frequencies $N \gg 1$, but not for low frequencies $N \ll 1$.  However, in the latter case we see that the components $i\Delta u_0$ and $-i\mu |u_0|^{p-1} u_0$ to the time variation $\partial_t u_0$ are both small compared to $u_0$ itself.  Informally, this suggests that while the low-frequency behaviour is nonlinear, this nonlinear behaviour will not manifest itself for some time.  Thus for short times we expect linear behaviour at both low and high frequencies, but for long times we expect nonlinear behaviour at low frequencies; in practice, this is reflected by the phenomenon that local existence is typically easy to establish at subcritical regularities, but that control of long-time asymptotics is very difficult unless one also has a critical or supercritical conservation law which prevents mass or energy from flowing completely to low frequencies.  If the mass $M$ increases, the time for which linear behaviour is expected will shrink, in some inverse polynomial relationship to the mass (which can also be deduced from dimensional analysis considerations).

Now we turn to the \emph{mass-supercritical} case, when $p - (1 + \frac{4}{d})$ is positive, it is the high frequencies which one expects to behave nonlinearly.  Furthermore, in this case $i\Delta u_0$ and $-i\mu |u_0|^{p-1} u_0$ are both large compared to $u_0$, so one expects the nonlinear behaviour to manifest itself very quickly.  Thus we expect supercritical equations to behave very badly; unless there is another property of the equation, such as energy conservation, which prevents mass from moving to high frequencies, it might happen that the mass concentrates at finer and finer scales, leading to blowup in finite time even from very smooth initial data.  Note that shrinking the mass $M$ may delay the time in which blowup occurs, but from scaling considerations we see that such shrinking cannot prohibit blowup entirely unless the mass is zero.  Thus, in the absence of any control of higher regularities on the time interval of interest, we expect the solution to be very unstable, and the Cauchy problem to either be illposed or to exhibit some form of blowup.  When the initial data is smooth in a supercritical equation, then one still expects local existence (because the high frequencies are initially quite small) but once the mass and energy flows into fine scales (e.g. by self-similar concentration, or by some sort of turbulence effect) it is not known in general what happens to the evolution.  (The notorious global regularity problem for the Navier-Stokes equations falls into this category, as all the known conserved or monotone quantities are supercritical.)

Now we turn to the \emph{critical} case, which for the mass in NLS occurs when $p = 1 + \frac{4}{d}$.  Now we see from
\eqref{nheu} that when the mass $M$ is small, we expect the linear behaviour to dominate the nonlinear behaviour at every scale; however, when the mass is large, it is possible at any given frequency scale $N$ for the nonlinear behaviour to dominate
the linear behaviour.  In such a case, one can check that $i\Delta u_0$ and $-i\mu |u_0|^{p-1} u_0$ have size roughly comparable to $N^2 u_0$, so that we expect the solution to stay close to the initial data $u_0$ only for time $O(1/N^2)$.
Thus we expect global existence, regularity, and scattering to a linear solution when the mass is small, but when the mass is large one only expects to the linear approximation to the solution to be valid for a time $T \sim 1/N^2$ depending on the natural frequency scale $N$ of the data (which can be arbitrary).  Beyond this time scale, one must account for nonlinear effects in order to determine the future behaviour of the evolution.  It is usually here that the sign of the nonlinearity (focusing, defocusing, or neither) is decisive.

The above heuristics can be remarkably accurate, but they are implicitly assuming that the rescaled bump functions are the ``worst'' type of initial data in a certain class (e.g. data with a certain prescribed mass), where by ``worst'' one means that the ratio between the nonlinear and linear components of the equation is strongest.  This is often the case, but when other symmetries than the scaling symmetry are present (particularly symmetries arising from a non-compact group) then one sometimes has to consider other types of data instead.  For instance, because of the Galilean invariance of NLS, one might expect \emph{frequency-modulated} bump functions such as $M^{1/2} e^{i \xi_0 \cdot x} \varphi(x)$ to be a competitor for
the title of worst initial data; more typically, hybrid examples such as rescaled frequency-modulated bumps $M^{1/2} N^{d/2} e^{i \xi_0 \cdot x} \varphi(Nx)$, whose Fourier transform is concentrated on some ball of radius $N$ centred at a frequency $\xi_0$, tend to play an important role.  In wave equations, \emph{Lorentz-transformed} bump functions (related to the \emph{Knapp example} in restriction theory) are also often of importance, when the Lorentz invariance is somehow ``stronger'' or ``higher-regularity'' than the scale-invariance.  See e.g. \cite{cct} for some discussion of the relative strengths of these symmetries for various classes of equations.

\section{Perturbation theory}

In the previous section we made some extremely informal computations regarding the ``ratio'' between the nonlinear and linear components of an equation for certain initial data, to then deduce predictions as to what the evolution should look like.  Now we formalise this intuition in the case where the linear behaviour is expected to dominate; in subcritical cases this corresponds to restricting time to a small interval depending on the norm of the initial data, while in critical cases this corresponds to either global solutions with small norm, or local solutions with large norm (and with time of existence depending on the initial data itself and not just on the norm).   

To achieve this formalisation, it is plausible that one should view the nonlinear equation as a perturbation of the linear equation, so that the nonlinearity is a kind of error term.  It turns out that one of the most effective ways to accomplish this is by converting the differential equation into an \emph{integral} (or Duhamel) equation, via the fundamental solution of the linear operator; this is basically because integral operators are far more likely to be bounded on various function spaces than differential operators.

To illustrate the method, we once again take the NLS \eqref{nls}, with initial data $u(0) = u_0$ in some data class, and solutions $u: I \times \R^d \to \C$ restricted to some time interval $I$.  (For second-order-in-time equations such as nonlinear wave equations, some slight modifications to the scheme below are needed to account for the initial velocity as well as initial position.) Typically one selects a Sobolev space such as $H^s_x(\R^d) = W^{2,s}_x(\R^d)$; these $L^2$-based spaces are preserved by the linear propagator $e^{it\Delta}$ (as can be seen from Plancherel's theorem) and thus have at least some chance of being stable under the nonlinear evolution as well.  The differential equation \eqref{nls} is then equivalent\footnote{This equivalence requires some mild regularity and decay assumptions on the solution; for instance, it will suffice that $u$ and $F(u)$ are both tempered distributions of spacetime which have some continuity in time.  In practice it is not difficult to justify these formal computations for the classes of solution that one is interested in, and we will not dwell on these technical issues here.} by Duhamel's formula
\begin{align*}
u(t) &= e^{it\Delta} u(0) + \int_0^t \frac{d}{dt'} [e^{i(t-t')\Delta} u(t')]\ dt' \\
&= e^{it\Delta} u(0) - i \int_0^t e^{i(t-t')\Delta} (iu_t + \Delta u)\ dt' 
\end{align*}
to the integral equation\footnote{In some cases it is convenient to apply a smooth time cutoff which equals $1$ on $I$ and vanishes outside of a neighbourhood of $I$, but this is a minor technical issue which we will not discuss here.}
\begin{equation}\label{duhamel}
 u(t) = e^{it\Delta} u_0 + (i\partial_t + \Delta)^{-1}(F(u)) (t)
\end{equation}
where $F$ is the nonlinearity function $F(z) := \mu |z|^{p-1} z$, $e^{it\Delta}$ is the propagator associated to the free Schr\"odinger equation $iu_t + \Delta u = 0$, or equivalently is defined via the Fourier inversion formula
$$ f(x) = \int_{\R^d} \hat f(\xi) e^{ix\cdot \xi}\ dx$$
as
$$ e^{it\Delta} f(x) = \int_{\R^d} \hat f(\xi) e^{-it|\xi|^2} e^{ix\cdot \xi}\ dx,$$
and $(i \partial_t + \Delta)^{-1}$ is the \emph{Duhamel operator}, defined by the formula
$$ (i\partial_t + \Delta)^{-1} f(t) :=  := -i\int_0^t e^{i(t-t')\Delta} f(t')\ dt'.$$
The first term on the right-hand side of \eqref{duhamel} if the nonlinearity $F()$ was absent, or in other words if one evolved purely by the linear evolution.  Thus the Duhamel formulation splits the nonlinear solution $u(t)$ as the sum of the linear solution $u_\lin(t) := e^{it\Delta} u_0$, and the cumulative effect $(i\partial_t + \Delta)^{-1}(F(u))(t)$ of the nonlinearity.  Thus we can view solutions $u$ of \eqref{nls} as fixed points of the map
\begin{equation}\label{ucos}
u \mapsto u_\lin + (i\partial_t + \Delta)^{-1} (F(u)).
\end{equation}
Note that $F$ is the only source of nonlinearity in this equation, while the initial data $u_0$ only intervenes via its linear development $u_\lin$.
To find fixed points of \eqref{ucos}, one surprisingly effective method (for semilinear evolution equations of the type discussed here) is the \emph{Duhamel iteration method} (also known as the \emph{contraction mapping method} or \emph{inverse function theorem method}), which is a variant of the classical \emph{Picard iteration method} and is one of the fundamental perturbative methods in the subject.  This method proceeds by establishing \emph{iterates} $u^{(j)}: I \times \R^d \to \C$ for $j=-1,0,1,\ldots$ recursively by setting $u^{(-1)} := 0$ and then setting 
\begin{equation}\label{ujj}
u^{(j)} := u_\lin + (i\partial_t + \Delta)^{-1}(F(u^{(j-1)}))
\end{equation}
for $j=0,1,\ldots$.  Thus for instance $u^{(0)}$ is just the linear solution $u_\lin$, while the \emph{first nontrivial iterate} $u^{(1)} = u_\lin + (i\partial_t + \Delta)^{-1}(F(u_\lin))$ is formed by combining the linear solution with the cumulative forcing term generated by that solution.  Further iterates become significantly more complicated to express non-recursively\footnote{In the case where $p$ is an odd integer, then the nonlinearity $F(z)$ is a polynomial of $z$ and $\overline{z}$, and the iterates can be expressed as a certain sum over $p$-ary trees with bounded size.  While this explicit expansion does clarify a few things, in particular the connection between the iteration method and the method of power series, it is unwieldy to work with in practice.}.  The strategy of the iteration method is then to conclude that the iterates $u^{(j)}$ converge (in suitable topologies) to a limit $u$; taking limits in \eqref{ujj} one should then obtain a fixed point of \eqref{ucos}, provided that $D$ and $F$ are continuous in appropriate topologies. 

In order to obtain this desired convergence, the standard approach is to show that the map \eqref{ucos} is not only continuous in some topology, but is in fact a \emph{Lipschitz} map from some complete metric space (typically a closed ball in a Banach space) to itself, with Lipschitz constant less than $\frac{1}{2}$ (say).  Then the existence of a fixed point follows from the contraction mapping theorem.  Furthermore, one automatically gains uniqueness of the fixed point (at least in the metric space used), as well as some stability properties relative to the linear solution $u_\lin$ (and hence on the initial data $u_0$).  If the nonlinearity $F$ is real analytic, then the solution map $u_0 \mapsto u_\lin$ will be also.  A basic way to achieve this Lipschitz behaviour is to design a Banach space $\S$ of functions on the spacetime slab $I \times \R^d$ to hold the solution $u$, and a Banach space $\N$ of functions on the same slab to hold the nonlinearity $F(u)$.  If one has the linear estimate
\begin{equation}\label{uduh} \| (i\partial_t + \Delta)^{-1} f \|_{\S} \leq C_0 \| f \|_{\N}
\end{equation}
and the nonlinear estimate
\begin{equation}\label{unil}
 \| F(u) \|_{\N} \leq C_1 \| u \|_{\S} \hbox{ whenever } \|u\|_{\S} \leq R
 \end{equation}
and more generally
\begin{equation}\label{ufuh}
\| F(u) - F(v) \|_{\N} \leq C_1 \| u - v\|_{\S} \hbox{ whenever } \|u\|_{\S}, \|v\|_{\S} \leq R
\end{equation}
for some $C_0, C_1, R > 0$
then we easily verify that the map \eqref{ucos} is a contraction on the complete metric space $\{ u \in \S: \|u\|_\S \leq R \}$
with Lipschitz constant at most $\frac{1}{2}$ whenever 
\begin{equation}\label{ulines}
\|u_\lin\|_{\S} \leq \frac{R}{2}
\end{equation}
and $C_0 C_1 \leq \frac{1}{2}$, thus generating a unique fixed point of \eqref{ucos} in this space.  (The quantity $C_0 C_1$ is a rigorous analogue of the informal concept of the ``ratio between the nonlinear and linear parts of the equation'' from the preceding section.)  Notice that this type of perturbative argument is insensitive to the sign $\mu$ of the nonlinearity, and so cannot be used to detect phenomena which are only present in the focusing case but not the defocusing case, or vice versa.

The task now reduces to one in harmonic analysis, namely to come up with spaces $\S, \N$ which obey the estimates
\eqref{ufuh}, \eqref{uduh}, \eqref{ulines} for suitable constants $C_0, C_1, R$.  In order to generate the smallness condition $C_0 C_1 \leq \frac{1}{2}$, one typically either has to make the initial data $u_0$ small (in order to allow $R$ and hence $C_1$ to be small, see \eqref{ulines}, \eqref{ufuh}) or to make the interval $I$ small (in order to make $C_0$ small, see \eqref{uduh} and the definition of $D$), or some combination of both (e.g. to make the size of $I$ small depending in some inverse manner on the norm of the initial data).  When the initial data lies in a scale-invariant space, one can use scaling considerations to see that without loss of generality we must take the spaces $\S$ and $\N$ to also be scale-invariant (note however that the nonlinearity $F(u)$ scales slightly differently frmo the solution $u$ itself).  This reduces the number of spaces and estimates available, which makes the harmonic analysis component of the argument slightly trickier, though as compensation the arguments are then insensitive to the exact length of the time interval involved and so can extend more readily to global control of solutions as opposed to merely local control.

As a simple example of the iteration strategy, the classical \emph{energy method} (or \emph{semigroup method}) for generating local solutions from initial data $u_0$ in a high regularity (and definitely subcritical) Sobolev space $H^s_x(\R^d)$ with $s > d/2$ proceeds by taking\footnote{We use $C^0_t H^s_x(I \times \R^d)$ to denote the Banach space of bounded continuous functions from $I$ to $H^s_x(\R^d)$ with the uniform norm.  This should be contrasted with the Frechet space $C^0_{t,\loc} H^s_x(I \times \R^d)$, which are the space of merely continuous (and thus \emph{locally} bounded) functions from $I$ to $H^s_x(\R^d)$.} $\S = \N = C^0_t H^s_x(I \times \R^d)$.  The linear estimate \eqref{uduh} is then true with $C_0 = |I|$ from Minkowski's inequality and the observation that the linear propagator $e^{it\Delta}$ preserves the $H^s_x(\R^d)$ norm.  The estimate \eqref{ulines} is similarly true so long as the initial data $u_0$ has $H^s_x$ norm less than $R/2$.  Finally, Schauder estimates combined with the hypothesis $s > d/2$ (which allows the $H^s_x$ norm to control boundedness and even H\"older continuity of the solution) imply (at least in the case when $p$ is an odd integer) that \eqref{ufuh} holds with $C_1 = C_{p,d} R^{p-1}$ for some constant $C_{p,d}$ depending only on $p$ and $d$.  Putting all this together, one obtains a local existence result for initial data in $H^s_x(\R^d)$ for an interval $I$ of length $|I| \approx \| u_0 \|_{H^s_x(\R^d)}^{-1/(p-1)}$.  It is instructive to compare this result against what one might expect from the scaling heuristics of the previous section.

While the energy method does give local existence and uniqueness for smooth solutions, it is unsatisfactory in a number 
of ways.  Firstly, it does not work at low regularities; in particular the energy class $H^1_x(\R^d)$ and the mass class $L^2_x(\R^d)$ are often out of reach of the energy method.  Secondly, and perhaps more importantly (from the perspective of smooth solutions), the time of existence given by this argument depends on a high-regularity norm $\|u_0\|_{H^s_x(\R^d)}$
rather than a lower regularity norm such as the energy norm.  This can cause difficulty when considering the long-time evolution of the equation, because low regularity norms are often easier to control (for instance via a conservation law) than higher regularity ones.  In some cases one can use \emph{ad hoc} methods, for instance using the Duhamel formula \eqref{duhamel} combined with harmonic analysis estimates and tools such as Gronwall's inequality or a bootstrap argument, to convert low regularity control (and high regularity control of the initial data) to high regularity control of the entire solution, thus allowing one to continue the solution globally.  However, it turns out that one can often obtain even more precise control on the solution by reworking the local existence argument so that it relies on less regularity on the initial data.  To do this, one must use finer properties of the linear equation $i u_t + \Delta u = 0$ (as represented both in the linear solution $u_\lin$ and in the Duhamel operator $(i\partial_t + \Delta)^{-1}$, and in particular in the \emph{dispersive} properties of this equation.  Informally, the dispersive property (which is the analogue of the \emph{elliptic regularity} effect for elliptic equations, or \emph{parabolic smoothing} effect for parabolic equations) asserts that solutions to this linear equation cannot concentrate significant amounts of mass or energy in small regions of space for extended periods of time; indeed, once a solution concentrates at one point in space and time, then at all later (or earlier) points in time, that component of the solution must disperse away from that point and towards spatial infinity.  There are many ways to capture this dispersive effect.  One basic and useful one is via the \emph{Strichartz inequalities}, which are the dispersive analogue of the well-known (and extremely fundamental) \emph{Sobolev inequalities} in elliptic theory, and control the boundedness of the propagators $e^{it\Delta}$ and $(i\partial_t + \Delta)^{-1}$ in various Sobolev and Lebesgue spaces.  There are many such Strichartz inequalities; a typical one is the estimate
\begin{equation}\label{strich}
\| (i \partial_t + \Delta)^{-1} f \|_{L^2_t L^{2d/(d-2)}_x(\R \times \R^d)} \leq C_d \| f \|_{L^2_t L^{2d/(d+2)}_x(\R \times \R^d)} 
\end{equation}
for all $d \geq 3$ and all spacetime test functions $f$ (see \cite{Keel-Tao}); compare this with the Sobolev inequality
$$ \| \Delta^{-1} f \|_{L^{2d/(d-2)}_x(\R^d)} \leq C_d \|f\|_{L^{2d/(d+2)}_x(\R^d)},$$
which is in fact a special case of the above Strichartz inequality, specialised to the limiting case of time-invariant functions.

Strichartz inequalities have been intensively studied; they ultimately arise from the $L^\infty_x(\R^d)$ decay properties in time of the fundamental solution $\frac{1}{(4\pi i t)^{d/2}} e^{i|x|^2/4t}$ of the propagator $e^{it\Delta}$.  Using these inequalities, one can develop a very satisfactory local (and in some cases global) well-posedness theory for NLS and NLW (excluding some technical cases of very low regularity or very rough nonlinearities) at the subcritical and critical regularities\footnote{Scaling arguments can be used to show that iteration methods must fail for supercritical regularities, and examples are known (especially in focusing cases) where the equation is either extremely unstable or for which blowup occurs instantaneously at these regularities. Our understanding of evolution in supercritical spaces, where the nonlinearity is significantly stronger than the linear part of the equation, is still extremely poor, and further progress may well require a radically different way to construct and control solutions.}.  For instance, the theory for NLS in the energy space $H^1_x(\R^d)$ for $d \geq 3$ in the energy-subcritical ($p < 1 + \frac{4}{d-2}$) case is as follows.

\begin{theorem}[LWP for energy-subcritical NLS]\label{nls-subcrit-lwp} Let $d \geq 3$, $p < 1 + \frac{4}{d-2}$, $\mu= \pm 1$, and $u_0 \in H^1_x(\R^d)$.  Then there exists a unique maximal Cauchy development $u \in C^0_{t,\loc} H^1_x(I \times \R^d)$, where $I \subset \R$ is an open time interval (possibly half-infinite or infinite) containing zero, which solves \eqref{nls} in the sense that \eqref{duhamel} holds.  Furthermore:
\begin{itemize}
\item (Lifespan estimate) We have $I \supset [-T,T]$ for some time $T \geq c_{d,p} \|u_0\|_{H^1_x(\R^d)}^{-C_{d,p}}$ and some constants $c_{d,p}, C_{d,p} > 0$ depending only on $d,p$.  Furthermore, if $p \geq 1 + \frac{4}{d}$ (i.e. the equation is not mass-supercritical\footnote{In the mass-supercritical case we in fact have global existence for arbitrary finite energy, or even finite mass, initial data, but this relies on the mass conservation law and so we do not include that result in this section, which is devoted to purely perturbative methods.}) and $\|u_0\|_{H^1_x(\R^d)} \leq \epsilon_{d,p}$ for some sufficiently small $\epsilon_{d,p} > 0$, then $I = \R$ (thus we have global existence for small energy data).
\item (Blowup criterion) If $T_*$ is a finite endpoint of $I$ then $\lim_{t \to T_*} \|u(t)\|_{H^1_x(\R^d)} = +\infty$.  (This follows easily from the lifespan estimate.)
\item (Persistence of regularity) If $u_0$ is Schwartz (resp. in $H^s_x(\R^d)$ for some $s \geq 0$) and $p$ is an odd integer, then $u$ will be smooth in space and Schwartz in time (resp. in $C^0_{t,\loc} H^s_x(I \times \R^d)$).
\item (Scattering criterion) Suppose $p \geq 1 + \frac{4}{d}$ (i.e. the equation is not mass-supercritical).  If $I$ contains $[0,+\infty)$ and $\|u\|_{L^{(p-1)(d+2)/2}_{t,x}([0,+\infty) \times \R^d))} < \infty$, then there exists a unique $u_+ \in  H^1_x(\R^d)$ such that $\lim_{t \to +\infty} \| u(t) - e^{it\Delta} u_+ \|_{H^1_x(\R^d)} = 0$.  Furthermore, if $u \in H^s_x(\R^d)$) for some $s \geq 0$ and $p$ is an odd integer, then $u_+$ is also in $H^s_x(\R^d)$ and $\lim_{t \to +\infty} \| u(t) - e^{it\Delta} u_+ \|_{H^s_x(\R^d)} = 0$.  Similarly if $I$ contains $(-\infty,0]$.
\item (Continuous dependence on the data) If $u_0^{(n)}$ is a sequence which converges in $H^1_x(\R^d)$ norm to $u_0$, and $J$ is a compact subinterval of $I$ containing zero, then for sufficiently large $n$ there exists solutions $u^{(n)}$ to \eqref{nls} (or \eqref{duhamel}) with initial data $u_0^{(n)}$ which converge to $u_0$ in $C^0_t H^1_x(J \times \R^d)$ norm.
\item (Energy and mass conservation) We have $E(u(t)) = E(u_0)$ and $M(u(t)) = M(u_0)$ for all $t \in I$.
\end{itemize}
\end{theorem}

\begin{remark}
The various components of this theorem are obtained by several variations on the iteration scheme discussed above, using various Sobolev and Lebesgue spaces to control the solution and nonlinearity, and using Sobolev and Strichartz estimates (together with such mundane tools as the Leibnitz rule and H\"older's inequality) to establish the required linear and nonlinear estimates.  See e.g. \cite{cazbooknew}, \cite{taobook}.  The energy and mass conservation laws are obtained by the usual \emph{density method}, namely by first establishing these results for smooth solutions (where everything can be easily justified rigorously) and then taking limits using the continuous dependence and persistence of regularity theory.  (When $p$ is not an odd integer, one sometimes also needs to smooth out the nonlinearity $F$ slightly; see \cite{cazbooknew}.)  There are more technical estimates one can obtain here, which roughly speaking assert that the solution $u$ obeys all the same estimates (up to a factor of two or so) as the linear solution $u_\lin$ on the interval $[-T,T]$ identified above, but we will not explicitly state those estimates here.  The hypothesis that $p$ be an odd integer is a technical one and is only needed when considering very high regularity solutions (e.g. in $H^s_x(\R^d)$ where $s > p$).  The spacetime norm $L^{(p-1)(d+2)/2}_{t,x}$ in the scattering criterion may seem arbitrary, but it is the unique pure Lebesgue spacetime norm which is invariant under the scaling of the equation.  It arises naturally when trying to stretch the iteration argument to noncompact time intervals such as $[T,+\infty)$ for large $T$ (which is what one needs to do to obtain the scattering result), as one can not afford to lose any power of the length of the time interval from H\"older's inequality when running such an argument.  Actually, one could replace this norm by several other scale-invariant norms, and often control of one such scale-invariant norm automatically implies control of many other scale-invariant norms.  We remark that energy class scattering for mass-supercritical data is unknown even if the norm is assumed to be small (the problem is somewhat similar to that of establishing local existence in supercritical norms), although in some cases one can still recover scattering results if additional decay conditions are placed on the data (e.g. $x u_0 \in L^2_x(\R^d)$).
\end{remark}

As $p$ approaches the energy-critical limit $p = 1 + \frac{4}{d-2}$, the exponent $C_{d,p}$ in the above theorem goes to infinity (as can be seen from scaling heuristics), and we obtain a slightly different local existence theorem:

\begin{theorem}[LWP for energy-critical NLS]\label{nls-crit-lwp} Let $d \geq 3$, $p = 1 + \frac{4}{d-2}$, $\mu= \pm 1$, and $u_0 \in \dot H^1_x(\R^d)$.  Then there exists a unique maximal Cauchy development $u \in C^0_{t,\loc} \dot H^1_x(I \times \R^d)$, where $I \subset \R$ is an open time interval (possibly half-infinite or infinite) containing zero, which solves \eqref{nls} in the sense that \eqref{duhamel} holds.  Furthermore:
\begin{itemize}
\item (Lifespan estimate) We have $I \supset [-T_-,T_+]$, where $T_-, T_+ > 0$ are any times for which
$\| u_\lin \|_{L^{2(d+2)/(d-2)}_{t,x}([-T_-,T_+] \times \R^d)} \leq \epsilon_d$, where $\epsilon_d > 0$ is a small constant depending only on $d$.  Furthermore, if $\|u_0\|_{\dot H^1_x(\R^d)} \leq \epsilon_d$, then $I = \R$ (thus we have global existence for small energy data).
\item (Blowup criterion) If $J$ is any subinterval of $I$ containing a finite endpoint of $I$ then $\|u\|_{L^{2(d+2)/(d-2)}_{t,x}(J \times \R^d)} = +\infty$.  
\item (Persistence of regularity) If $u_0$ is Schwartz (resp. in $\dot H^s_x(\R^d)$ for some $s \geq 0$) and $p$ is an odd integer, then $u$ will be smooth in space and Schwartz in time (resp. in $C^0_{t,\loc} \dot H^s_x(I \times \R^d)$).
\item (Scattering criterion) If $I$ contains $[0,+\infty)$ and $\|u\|_{L^{2(d+2)/(d-2)}_{t,x}([0,+\infty) \times \R^d))} < \infty$, then there exists a unique $u_+ \in \dot H^1_x(\R^d)$ such that $\lim_{t \to +\infty} \| u(t) - e^{it\Delta} u_+ \|_{\dot H^1_x(\R^d)} = 0$.  Furthermore, if $u \in \dot H^s_x(\R^d)$) for some $s \geq 0$ and $p$ is an odd integer, then $u_+$ is also in $\dot H^s_x(\R^d)$ and $\lim_{t \to +\infty} \| u(t) - e^{it\Delta} u_+ \|_{\dot H^s_x(\R^d)} = 0$.  Similarly if $I$ contains $(-\infty,0]$.
\item (Continuous dependence on the data) If $u_0^{(n)}$ is a sequence which converges in $\dot H^1_x(\R^d)$ norm to $u_0$, and $J$ is a compact subinterval of $I$ containing zero, then for sufficiently large $n$ there exists solutions $u^{(n)}$ to \eqref{nls} (or \eqref{duhamel}) with initial data $u_0^{(n)}$ which converge to $u_0$ in $C^0_t \dot H^1_x(J \times \R^d)$ norm.
\item (Energy and mass conservation) We have $E(u(t)) = E(u_0)$ and (if $u_0 \in L^2_x(\R^d)$) $M(u(t)) = M(u_0)$ for all $t \in I$.
\end{itemize}
\end{theorem}

Here, we see that the spacetime scale-invariant norm $L^{2(d+2)/(d-2)}_{t,x}$ plays a governing role in the existence of the solution.  Very roughly speaking, when this norm is small, the solution behaves linearly; when the norm is large but finite, the solution behaves nonlinearly but does not blow up, and even scatters to a free solution at $t=\pm \infty$; and when the norm is infinite, then the solution of course blows up.  The above results are achieved by pure perturbative analysis, relying only on variants of the iteration method and on harmonic analysis estimates such as Strichartz and Sobolev inequalities; see \cite{cwI}, \cite{cazbooknew}, \cite{tao-visan}. 

We have seen how perturbative analysis allows one to demonstrate existence, uniqueness, regularity, and spacetime bounds on solutions.  Another important application of perturbation theory is in showing that equations such as \eqref{nls} are \emph{stable}, in the sense that one can add or remove small additional forcing terms to the right-hand side (or to the initial data) without significantly affecting the evolution.  Thus for instance if $v$ \emph{approximately} solves \eqref{nls} in the sense that
\begin{equation}\label{uvt}
 i v_t + \Delta v = F(v) + e 
 \end{equation}
for some small $e$, and $v(0)$ is close to $u_0$ in some suitable norm, then we expect $v$ to be close to the \emph{exact} solution $u$ to \eqref{nls} with initial data $u_0$,
\begin{equation}\label{uut}
i u_t + \Delta u = F(u); \quad u(0) = u_0
\end{equation}
for short times at least  This type of stability result has a number of uses.  Firstly, it can permit one to use the model equation (in this case, NLS) to approximate more complicated equations from which the model was derived (by dropping various ``small'' terms).  Related to this, one can use stability results to rigorously justify the convergence of various numerical schemes to the exact equation, thus allowing for rigorous numerical results for this equation.  Finally, it gives a powerful method to construct exact solutions to the equation, namely by first constructing a sufficiently accurate \emph{approximate solution} to the equation (for instance, by some asymptotic expansion, or by suppressing some nonlinear interactions from the equation), and then using the stability theory to perturb the approximate solution to a nearby exact solution.

There are many stability results in the literature.  The basic idea is to express $v$ as a perturbation of $u$ or vice versa, and solve for the difference.  For instance, if we write $u = v+w$, then $w$ is small at time zero and solves the \emph{difference equation}
$$ i w_t + \Delta w = F(v+w) - F(v) + e.$$
One can then use iterative methods (or other perturbative methods, such as the energy method and Gronwall's inequality)
to control $w$, at least for short and medium times.  A typical stability result, for the energy-critical NLS discussed above, is as follows.

\begin{theorem}[Long-time perturbations]\label{long-time theorem}\cite{tao-visan} Let $d \geq 3$, $p = 1 + \frac{4}{d-2}$, and $\mu = \pm 1$.
Let $I$ be a time interval containing $0$ and let $v \in C^0_t \dot H^1_x(I \times \R^d)$ solve \eqref{uvt}
with the bounds
\begin{align*}
\|v\|_{L_{t,x}^{\frac{2(d+2)}{d-2}}(I \times \R^d)}&\leq M\\
\|v\|_{C_t^0\dot{H}^1_x(I \times \R^d)}&\leq M \\
\|\nabla e\|_{L^2_t L^{\frac{2d}{d+2}}(I \times \R^d)} &\leq \eps
\end{align*}
for some $M > 0$, $\eps > 0$.  Suppose also that $u_0 \in \dot H^1_x(\R^d)$ is such
that $\| v(0) - u_0 \|_{\dot H^1_x} \leq \eps$.  Then if $\eps$ is sufficiently small depending on $d,M$,
Then there exists a solution $u$ to \eqref{uut} such that
$$ \|u-v\|_{L_{t,x}^{\frac{2(d+2)}{d-2}}(I \times \R^d)} + \|u-v\|_{C^0_t \dot H^1_x(I \times \R^d)} \leq C(M,d) 
\bigl(\eps+\eps^{\frac{7}{(d-2)^2}}\bigr) 
$$
for some $C(M,d) < \infty$ depending only on $\eps$ and $d$.
\end{theorem}

The exponent $\frac{7}{(d-2)^2}$ is a technicality arising from the low regularity of the nonlinearity $F()$ in higher dimensions and should be ignored.  The stability result in \cite{tao-visan} is in fact slightly stronger than stated here but we have given a simplified version for sake of exposition.  The argument is purely perturbative; the key idea is to first subdivide the interval $I$ so that the $L^{\frac{2(d+2)}{d-2}}_{t,x}$ norm of $v$ is small rather than merely finite, and then to apply perturbative arguments of the type sketched above to each subinterval separately.  This type of stability result turns out 
to play a crucial role in the large data theory for critical equations, as it is usefully encapsulates a large portion of the perturbative theory.

\subsection{Other function spaces}

The above considerations for NLS in the energy class have analogues for the other equations listed previously, at various levels of regularity.  For the NLS and NLW equations, which have no derivatives in the nonlinearity, the Strichartz estimates are sufficient to establish a satisfactory theory.  However, for the more complicated models which contain derivatives, the need to establish (the analogue of) the estimate \eqref{unil} will force the nonlinearity space $\N$ to be at least one derivative rougher in regularity than the solution space $\S$.  Inspecting \eqref{uduh}, we thus see that the task then falls to the Duhamel operator (such as $(i \partial_t + \Delta)^{-1}$, $\Box^{-1}$, or $(\partial_t + \partial_{xxx})^{-1}$) to ``recover'' this loss of derivative.  This is often not possible to establish with Strichartz estimates alone (except sometimes when the linear part is second-order in time, which is the case with nonlinear wave equation models), and so more advanced spaces have been developed for this recovery of derivatives.  In the case of highly dispersive models such as the gKdV equations,
it turns out that \emph{local smoothing estimates} (coupled with the more technical \emph{maximal function estimates} that give some complementary local control on the solution) are a useful tool.  A typical local smoothing estimate (first observed by Kato) is as follows: if $u \in C^0_t L^2_x(\R \times \R \to \R)$ solves the Airy equation $u_t + u_{xxx} = 0$, then we have
$$ \int_{0}^1 \int_{-1}^1 u_x(t,x)^2\ dx dt \leq C \int_\R u(0,x)^2\ dx$$
for some absolute constant $C$.  Note the gain of one degree of regularity on the left-hand side.  This particular estimate
can be proven by a direct integration by parts argument, using firstly the conservation of the $L^2$ mass $\int_\R u(t,x)^2\ dx$ and secondly the monotonicity of a weighted $L^2$ mass such as $\int_\R \tanh^{-1}(x) u(t,x)^2\ dx$; we omit the details.  More refined local smoothing estimates can be proven by harmonic analysis techniques, in particular invoking the Fourier transform,
which can then be be used to give local wellposedness results for the gKdV equation which are largely sharp; see \cite{kpv:gkdv}.  

When approaching critical regularities, it seems that even local smoothing and maximal function estimates are not sufficient.  For slightly subcritical regularities, a very useful tool has been the development of the \emph{Fourier restriction norm spaces} $X^{s,b}$ (also called $H^{s,b}$) developed by Bourgain \cite{borg:xsb} for nonlinear dispersive equations and by Klainerman and Machedon \cite{klainerman:xsb} for nonlinear wave equations\footnote{These spaces also appeared in earlier work on propagation of singularities in \cite{beals:xsb}, \cite{rauch.reed}.}.  These spaces are to dispersive and wave equations as Sobolev spaces
are to elliptic equations.  For sake of discussion let us work with the $X^{s,b}$ spaces associated with the Schr\"odinger operator $(i \partial_t + \Delta)$.  Just as a Sobolev space $H^s_x(\R^d)$ is essentially given for $s \in \R$ by the norm
$$ \| u \|_{H^s_x(\R^d)} \approx \| \langle \nabla \rangle^s u \|_{L^2_x(\R^d)},$$
where $\langle x \rangle := (1 + |x|^2)^{1/2}$ is the Japanese bracket, interpreted appropriately for operators such as $\nabla$ using a functional calculus, the $X^{s,b}(\R \times \R^d)$ spaces are essentially given for $s,b \in \R$ by the norm
$$ \| u \|_{X^{s,b}(\R^d)} \approx \| \langle \nabla \rangle^s \langle i \partial_t + \Delta \rangle^b u \|_{L^2_{t,x}(\R \times \R^d)}.$$
To formalise this properly one needs the spacetime Fourier transform, and there are also some technical adjustments needed to localise this norm to a compact time interval.  For details see \cite{ginibre-survey}.

The indices $s$ and $b$ measure the ``elliptic'' and ``dispersive'' regularity of the solution respectively.  The power of 
these spaces lies in the fact that they fully capture the smoothing effect of the Duhamel operator $(i \partial_t + \Delta)^{-1}$; indeed, to oversimplify substantially, this operator is essentially an isometry from $X^{s,b-1}$ to $X^{s,b}$ for all $s$ and ``reasonable'' values of $b$.  Strichartz estimates can be reinterpreted as ``dispersive Sobolev embedding theorems'' from the $X^{s,b}$ spaces to other Lebesgue spaces.  The task of establishing nonlinear estimates such as \eqref{unil} in these spaces requires a certain amount of multilinear harmonic analysis but the techniques for doing so are now rather well understood; see e.g. \cite{tao:multi}.

At the critical regularity, even the $X^{s,b}$ spaces begin to break down.  The problem is similar to that faced in Sobolev spaces, when the fundamental Sobolev embedding $H^s_x(\R^d) \subset L^\infty_x(\R^d)$ breaks down at the endpoint $s = d/2$.  However, critical substitutes for these $X^{s,b}$ spaces are known, thanks to the work of Tataru \cite{tataru:box5}, \cite{tataru:wave1}, \cite{tataru:wave2}.  These substitutes are rather technical and messy to describe, but roughly speaking they combine Besov-space variants of the $X^{s,b}$ spaces with certain spacetime frequency-localised versions of Strichartz spaces; the idea is to use $X^{s,b}$ type control in the ``non-resonant'' region where the symbol of the linear operator is large, and Strichartz type control in the ``resonant'' region when the symbol is small.  In low dimensions, when the standard Strichartz estimates are weak, one must also sometimes introduce more exotic Strichartz estimates, for instance adapted to null frames.  This is in particular the case for two-dimensional wave map equations; see \cite{tataru:wave2}.

The need to use Besov spaces at the critical level means that perturbation theory often hits a natural limit at the scale-invariant \emph{Besov space} $\dot B^{s,1}_2(\R^d)$ rather than the scale-invariant Sobolev space $\dot H^s(\R^d)$.  To break this barrier for wave maps (and more recently for Schr\"odinger maps) has required the additional technique of gauge transformations; see Section \ref{gauge-sec}.

At present it seems that our collection of function spaces and estimates are sufficient for the subcritical and critical perturbative theory for most of the standard model equations, although some of the spaces are rather messy and one can hope for further simplification in the future.  There are a variety of results and heuristics which indicate that the supercritical theory is out of reach of perturbation theory, no matter how refined the spaces and estimates one uses.  Firstly, there is the problem that perturbation theory does not seem able to exploit the defocusing sign in a nonlinearity, which appears to be essential in the supercritical theory since focusing equations often blow up instantaneously at supercritical regularities.  Secondly, there are a number of instability results \cite{cct}, \cite{lebeau} for supercritical equations which are inconsistent with the type of control that perturbative techniques naturally give.  Finally, basic dimensional analysis shows that it is not possible to simultaneously have all three estimates \eqref{uduh}, \eqref{ufuh}, \eqref{ulines} for any supercritical data class.  Thus the establishment of a good existence theory for supercritical data classes\footnote{If the degree of supercriticality is only logarithmic, then it turns out that one can sometimes augment the perturbative method with nonlinear \emph{a priori} estimates to continue to control the solution; see \cite{tao:supercrit}.} will have to rely on some sort of non-perturbative method which fully exploits the defocusing nature of the nonlinearity.

\subsection{Alternatives to perturbative methods}\label{alternative-sec}

To close this section, we should emphasise that perturbative techniques, while very effective in the regime where the linear behaviour dominates the nonlinear behaviour, are not the only way to construct solutions; we mention two key ones here.

An important non-perturbative method to construct solutions is the \emph{weak compactness method}, in which penalisation, 
viscosity, discretisation or other approximation methods (generally based on suppressing fine-scale behaviour) are used to construct a family of \emph{approximate} solutions to the equation, obtain uniform bounds on such solutions (typically using conservation laws) and then weak limits extracted to obtain a limiting object which solves the equation in some weak sense.  This method is very robust and can work even for large data in supercritical equations provided that one has a sufficiently positive-definite conservation law.  However, the solution obtained is typically of low regularity (e.g. the energy class) even when the initial data is smooth, and \emph{a priori} is only known to solve the equation in a weak (distributional) sense.  This has some non-trivial consequences regarding the justification of various formal computations regarding such solutions; for instance, a quantity which is conserved for smooth solutions may merely be non-increasing for weak solutions (due to the inequality in Fatou's lemma, for instance).  Substantial additional work is often needed to upgrade the solution to be regular, unique, or to enjoy conservation laws.  To give one example, the existence of global weak solutions for the Navier-Stokes equations from smooth initial data has been known for over seventy years, thanks to the work of Leray, but to this date there has been little progress in upgrading these weak solutions to a globally smooth solution (except when the initial data is small, or some other bound is assumed on the solution).  The basic ``enemy'' in the weak solution method, namely the cascade of energy from coarse scales to fine scales, is ultimately the same as the one encountered in perturbation theory when trying to extend local existence of smooth solutions to global existence, and so it appears that working with weak solutions does not allow one to automatically evade this fundamental obstruction to global regularity. On the other hand, a close relative of the 
weak compactness method, the \emph{concentration compactness method}, has recently proven to be very useful in analysing 
global behaviour or blowup behaviour of these equations, by isolating the key ``blowup profiles'' of the evolution; see Section \ref{ccsec}.

Another major development has been to extend the reach of both perturbative and non-perturbative methods by various nonlinear transformations, most notably \emph{normal forms} and \emph{gauge transforms}, in order to reduce the strength of the nonlinear component of the equation.  (The Miura transform connecting KdV and mKdV also falls into this category.)  Normal form transformations are often motivated from considerations in Hamiltonian dynamics or symplectic geometry, and seek to transform either the equation or the Hamiltonian (often by a symplectic transformation which is a perturbation of the identity map) in order to remove or attenuate the ``non-resonant'' portions of the nonlinearity, possibly replacing them with higher order terms.  While these techniques are important in many problems in this field, they have so far not made much impact on the critical-regularity theory and so we shall not discuss them here.  Gauge transforms, on the other hand, tend to arise from considerations in differential geometry, and can be effective in reducing the strength of nonlinearities which contain first-order derivatives of the solution.  We discuss these in Section \ref{gauge-sec}.

With the important exception of the completely integrable equations, the number of demonstrably effective methods to construct reasonable\footnote{What ``reasonable'' means is of course somewhat subjective, but at a bare minimum, solutions should have some existence and uniqueness theory, be compatible with more classical concepts of a solution, and basic physical properties of these solutions such as conservation laws should be rigorously justifiable.} solutions to nonlinear dispersive and wave equations from general data still remains unacceptably low compared to other areas of PDE.  There are some variants of the basic Duhamel iteration method, such as the Nash-Moser iteration scheme, but while this scheme is rather useful for quasilinear equations, it does not seem to be more effective than Duhamel iteration for semilinear equations.  The classical method of power series expansions (as used for instance to prove the Cauchy-Kowalevski theorem) is useful for real-analytic classes of initial data, but for non-analytic data it seems to be essentially equivalent in strength to (and messier to use than) the Duhamel iteration method.  The lack of anything resembling a maximum principle or comparison principle prevents comparison methods from being effective (except in demonstrating blowup for scalar wave equations), in sharp contrast to elliptic and parabolic PDE.  Similarly, the extreme non-convexity (and non-Palais-Smale nature) of the Lagrangian functionals for these equations has so far prevented the use of variational methods (though see Section \ref{energy-sec}).  Kinetic formulations (for instance, transforming Schr\"odinger equations via the FBI or Wigner transforms) have so far also failed to noticeably improve the existence theory for these equations.  There are also essentially no known topological, dynamical, symplectic, or stochastic methods to construct solutions to these PDE, with the possible exception of some isolated work in constructing invariant measures.  Any new method to construct solutions for such PDE along these or other lines may well represent a significant breakthrough in the field.

\section{Conservation laws}

Having discussed the perturbative theory in the previous section, we now turn to the topic of \emph{non-perturbative} methods for analysing nonlinear dispersive equations, which are valid even for large data or long times (in other words, in regimes where the nonlinear component of the evolution is not insignificant).  For equations which are not completely integrable, one relies primarily on three types of non-perturbative tools: conservation laws, monotonicity formulae, and transformations (such as gauge transformations).  This is admittedly a small list of techniques, and it would be of great interest to develop additional typs of non-perturbative methods.

In this section we discuss conservation laws and how they are used.  One can approach conservation laws either from a algebraic perspective (multiplying the equation against various well-chosen multipliers and then integrating by parts), from a Fourier analytic perspective (studying which multilinear Fourier multipliers of the solution are preserved by the flow), from a Hamiltonian perspective (connecting conserved quantities to symmetries of the equation or Hamiltonian, via Noether's theorem), or from a Lagrangian perspective (viewing conserved quantities in terms of symmetries of the Lagrangian).  All four perspectives are important; for sake of exposition we shall focus here on just one approach, based on the Lagrangian perspective.  (See \cite{taobook} for some discussion of the other approaches.)  This approach is especially well suited to geometric equations, such as the nonlinear wave equations on Minkowski space, as one can take advantage of the diffeomorphism invariance of such equations to obtain a stress-energy tensor which is \emph{pointwise} conserved.  This is in contrast to the Hamiltonian approach, in which finite-dimensional symmetries are used to generate finitely many conserved \emph{integrals}; the infinite-dimensional diffeomorphism symmetry is significantly more powerful than finite-dimensional sub-symmetries (such as translation or rotation symmetry), and the pointwise control will be essential for establishing the monotonicity formulae of the next section.

For sake of discussion, let us consider the nonlinear wave equation NLW, normalised so that $c=1$, although the approach here is very general and applies to any geometric equation associated to a Lagrangian.  We shall work formally for now, ignoring issues such as integrability or regularity; once the form of the conservation laws are obtained, they can be justified rigorously by a number of means.

We view this equation as the Euler-Lagrange equation for the action
\begin{equation}\label{sudef}
S(u, g) := \int_{\R^{1+d}} \frac{1}{2} g^{\alpha \beta} \partial_\alpha u \partial_\beta u
+ \frac{\mu}{p+1} |u|^{p+1} \ dg = \int_{\R^{1+d}} L(u,g)\sqrt{-\det{g}}\ dx dt
\end{equation}
where $g$ is the Minkowski metric $g_{\alpha \beta} x^\alpha x^\beta = - t^2 + x_1^2 + \ldots + x_d^2$, 
$dg = \sqrt{-\det(g)}\ dx dt$ is the associated volume form, and $L(u,g)$ is the Lagrangian density
$$ L(u, g) := \frac{1}{2} g^{\alpha \beta} \partial_\alpha u \partial_\beta u
+ \frac{\mu}{p+1} |u|^{p+1}.$$
Thus if $u$ solves \eqref{nlw}, then $u$ is a critical point for $S(u,g)$ with $g$ fixed:
\begin{equation}\label{lig}
 \frac{\delta S}{\delta u}(u,g) = 0.
\end{equation}

On the other hand, the action $S(u,g)$ is clearly invariant under diffeomorphisms $\phi: \R^{1+d} \to \R^{1+d}$ of the underlying spacetime manifold $\R^{1+d}$:
$$ S( u \circ \phi, \phi_* g ) = S(u, g).$$
In particular, if we consider infinitesimal diffeomorphisms $e^{\eps X}$ associated to an arbitrary (smooth) vector field $X: \R^{1+d} \to T \R^{1+d}$ we have
$$ \frac{d}{d\eps} S( u \circ e^{-\eps X}, (e^{\eps X})_* g )|_{\eps = 0} = 0.$$
From the chain rule, the left-hand side is
$$ - \frac{\delta S}{\delta u}(u,g)[X^\alpha \partial_\alpha u] + \frac{\delta S}{\delta g}(u,g)[{\mathcal L}_X g]$$
where ${\mathcal L}_X g$ is the Lie derivative of $g$ along the vector field $X$.  Applying \eqref{lig} we conclude that
$$ \frac{\delta S}{\delta g}(u,g)[{\mathcal L}_X g] = 0$$
for arbitrary smooth vector fields $X$.  From differential geometry we recall the formula $({\mathcal L}_X g)_{\alpha \beta} = \pi_{\alpha \beta}$, where $\pi$ is the \emph{deformation tensor}
\begin{equation}\label{pidef}
\pi_{\alpha \beta} = \nabla_\alpha X_\beta + \nabla_\beta X_\alpha
\end{equation}
where $\nabla$ is the Levi-Civita connection with respect to the metric $g$ (in the case of the Minkowski metric, this is the same as the ordinary partial derivative $\partial$).  Applying \eqref{sudef}, we can then write
$$  \frac{\delta S}{\delta g}(u,g)[{\mathcal L}_X g] =
\int_{\R^{1+d}} [\frac{\partial L}{\partial g^{\alpha \beta}}(u,g) \pi^{\alpha \beta} - \frac{1}{2} L(u,g) g_{\alpha \beta} \pi^{\alpha \beta}] \sqrt{-\det{g}}\ dx dt.$$
If we then define the \emph{stress energy tensor}
$$ \stress_{\alpha \beta} := \frac{\partial L}{\partial g^{\alpha \beta}}(u,g) - \frac{1}{2} L(u,g) g_{\alpha \beta}$$
we conclude that
$$ \int_{\R^{1+d}} \stress_{\alpha \beta} \pi^{\alpha \beta}\ dg = 0$$
for all smooth vector fields $X$.  Using \eqref{pidef} and the symmetry of $\stress$ we conclude that
$$ \int_{\R^{1+d}} \stress_{\alpha \beta} \nabla^\alpha X^\beta\ dg = 0$$
for arbitrary $X$; integrating by parts and using duality we then conclude the pointwise \emph{conservation of stress-energy}
\begin{equation}\label{conserv}
\nabla_\alpha \stress^{\alpha \beta} = 0.
\end{equation}
In co-ordinates, we thus have
\begin{equation}\label{conserv-coord}
\partial_t \stress^{00} + \partial_j \stress^{0j} = 0; \quad \partial_t \stress^{k0} + \partial_j \stress^{kj} = 0.
\end{equation}
The above computations can be performed for an arbitrary geometric wave equation, though the precise form of $L$ (and hence $\stress$) of course varies from equation to equation.  In the specific case of the NLW, we have
\begin{align*}
\stress^{\alpha \beta} &= (\partial^\alpha u) (\partial^\beta u) - g^{\alpha \beta} (\frac{1}{2} \partial^\gamma u \partial_\gamma u + \frac{2\mu}{p+1} |u|^{p+1}) \\
&= (\partial^\alpha u) (\partial^\beta u) - g^{\alpha \beta} [\frac{1}{4} \Box(|u|^2) - \frac{\mu(p-1)}{p+1} |u|^{p+1}]
\end{align*}
or in coordinates
\begin{align*}
\stress^{00} &= \frac{1}{2} |u_t|^2 + \frac{1}{2} |\nabla u|^2 + \frac{\mu}{p+1} |u|^{p+1} \\
\stress^{0j} = \stress^{j0} &= - \Re( \overline{u_t} u_j) \\
\stress^{jk} &= \Re(\overline{u_j} u_k) - \delta_{jk} (\frac{1}{2} |\nabla u|^2 - \frac{1}{2} |u_t|^2 + \frac{\mu}{p+1} |u|^{p+1})
\end{align*}
where $\delta_{jk}$ is the Kronecker delta.  The density $\stress^{00}$ is known as the \emph{energy density}, while the vector $\stress^{0j}$ is the \emph{energy current} or \emph{momentum density}.  The tensor $\stress^{jk}$ is the \emph{momentum current} or the \emph{stress tensor}.

The pointwise conservation law \eqref{conserv} (or \eqref{conserv-coord}) has many uses.  One of the simplest is obtained simply by integrating \eqref{conserv-coord} in space and using Stokes' theorem, to obtain (formally, at least) 
$$ \partial_t \int_{\R^d} \stress^{00}(t,x)\ dx = \partial_t \int_{\R^d} \stress^{k0}(t,x)\ dx = 0.$$
Thus the \emph{total energy}
$$ E(u[t]) := \int_{\R^d} \stress^{00}(t,x)\ dx = \int_{\R^d} \frac{1}{2} |u_t|^2 + \frac{1}{2} |\nabla u|^2 + \frac{\mu}{p+1} |u|^{p+1}\ dx$$
and the \emph{total momentum}
$$ p^k(u[t]) := \int_{\R^d} \stress^{k0}(t,x)\ dx = -\int_{\R^d}  \Re( \overline{u_t} u_j) 
\ dx$$
are conserved quantities.  We will see further consequences of the conservation laws in the next section.

Recall from Section \ref{nls-eq} that the NLS can be embedded in the NLW of one higher dimension.  Thus the stress-energy conservation law for NLW must have some analogue for NLS.  If one performs the algebraic computations (using the null coordinate frame $\partial_t \pm \partial_{d+1}, \partial_1, \ldots, \partial_d$), one sees that the $d+1$-dimensional stress-energy conservation law for NLW decouples into a $d$-dimensional stress-energy conservation law for NLS
\begin{equation}\label{conserv-coord-again}
\partial_t \stress^{00} + \partial_j \stress^{0j} = 0; \quad \partial_t \stress^{k0} + \partial_j \stress^{kj} = 0.
\end{equation}
where the \emph{pseudo-stress-energy tensor} $\stress^{\alpha \beta}$ is defined by
\begin{align*}
\stress^{00} &:= |u|^2 \\
\stress^{0j} = \stress^{j0} &:= 2\Im( \overline{u} u_j) \\
\stress^{jk} &:= 
4\Re(\overline{u_j} u_k) - \delta_{jk} \Delta(|u|^2) + \frac{4\mu(p-1)}{p+1} \delta_{jk} |u|^{p+1}
\end{align*}
and an additional (scalar) energy conservation law
\begin{equation}\label{energy-conserv}
\partial_t e^0 + \partial_j e^j = 0
\end{equation}
where the \emph{energy density} $e^0$ and \emph{energy current} $e^j$ are defined as
$$ e^0 := \frac{1}{2} |\nabla u|^2 + \frac{\mu}{p+1} |u|^{p+1}; \quad e^j := \Im( \overline{u_{jk}} u_k ) + \mu |u|^{p-1} \Im( \overline{u} u_j ).$$
We thus obtain three important conserved quantities, namely the \emph{total mass}
$$ M(u(t)) := \int_{\R^d} \stress^{00}(t,x)\ dx = \int_{\R^d} |u(t,x)|^2\ dx$$
the \emph{total momentum}
\begin{equation}\label{total-momentum}
 p^k(u(t)) := \int_{\R^d} \stress^{k0}(t,x)\ dx = 2 \int_{\R^d} \Im(\overline{u} u_k)\ dx
\end{equation}
and the \emph{total energy}
$$ E(u(t)) := \int_{\R^d} e^0(t,x)\ dx = \int_{\R^d} \frac{1}{2} |\nabla u|^2 + \frac{\mu}{p+1} |u|^{p+1}\ dx.$$

Similar conservation laws can also be deduced for the other equations (gKdV, SM, WM, MKG, YM) discussed earlier, although for certain equations (notably gKdV and SM) the Lagrangian formulation is not as convenient as the Hamiltonian formulation for locating the conserved quantities.  In the case of the equations with a covariant wave or Schr\"odinger equation (e.g. MKG, WM, SM, NLW, NLS) there are also ``charge conservation laws'' arising from the gauge group, but these have limited usefulness for the analysis of these equations, as neither the charge density nor the charge current enjoy any positivity properties in general\footnote{An exception is NLS, in which the conserved charge density arising from the phase rotation symmetry $u \mapsto e^{i\theta} u$ is in fact the same as the conserved mass density $|u|^2$.  This is because the embedding of NLS into NLW identifies phase rotation with translation in a spacetime null direction, and the mass density is nothing more than the component of the NLW stress-energy tensor in that direction.}.

In the theory of ODE, a conservation law (such as energy conservation) restricts the dynamics to a lower-dimensional subset of phase space, such as the \emph{energy surface} where the energy is constant.  If the conservation law is sufficiently \emph{coercive} (so that the conserved quantity goes to infinity at phase space infinity), then this subset will be bounded.  ODE existence theorems such as the Picard existence theorem thus ensure global existence for the evolution.

In the theory of PDE, which can be viewed as an infinite-dimensional analogue of ODE, the situation is more complicated because there are many inequivalent norms with which to measure the ``boundedness'' of a subset of phase space, and a conserved quantity can give control in one norm whereas the criterion needed for the local existence theory to prevent blowup may require another norm.  A related issue is that even when the energy surfaces are bounded, they are usually quite non-compact.  However, when the conservation laws and the local existence theory are both sufficiently strong, one can combine the two to still obtain global existence.  Typically, this compatibility between the conservation laws and the local existence theory only occurs when a key conserved quantities is subcritical; a large part of recent developments have centred on extending this compatibility to the case when the key conserved quantity is critical.

Let us illustrate the above discussion with the defocusing NLS $\mu=+1$ with subcritical or critical energy (thus we have $d \leq 3$ or $p \leq 1 + \frac{4}{d-2}$), we see from Sobolev embedding that
$$ c_d \|u(t)\|_{H^1_x(\R^d)}^2 \leq M(u(t)) + E(u(t)) \leq C_d ( \|u(t)\|_{H^1_x(\R^d)}^2 + \| u(t)\|_{H^1_x(\R^d)}^{p+1} )$$
for some constants $c_d, C_d > 0$ depending only on the dimension $d$.  Since $M(u)$ and $E(u)$ are conserved in these cases (by 
Theorems \ref{nls-subcrit-lwp}, \ref{nls-crit-lwp}, we thus see that if the solution is initially in $H^1_x(\R^d)$, then it will be bounded in $H^1_x(\R^d)$ throughout the entire lifespan of the solution.  In the subcritical case, the blowup criterion in Theorem \ref{nls-subcrit-lwp} then immediately shows that the solution is in fact global.  In the focusing case $\mu=-1$, the above argument does not quite work directly because $E$ contains a negative component.  However, it turns out that in the mass-subcritical case $p < 1 + \frac{4}{d}$ (or the mass-critical case with small mass) one can use the Gagliardo-Nirenberg inequality to show that the positive (linear) component of the energy $E$ still dominates the negative part when the solution is large, and so one can continue to obtain global existence in this case.

The critical case is however much more delicate, because the blowup condition given by Theorem \ref{nls-crit-lwp} is not precluded by the boundedness of the $H^1_x(\R^d)$ norm.  Global existence in the energy class is indeed known in the defocusing energy-critical setting, but this is an extremely recent and difficult result.  To illustrate the difficulty, let us consider the mass-critical focusing NLS ($p = 1 + \frac{4}{d}$, $\mu = -1$).  For this equation, it is known that there is local existence from $L^2_x(\R^d)$ initial data, and even global existence if the mass is small.  However, in the large mass case
the time of existence depends on the data itself and not just on the mass.  In particular, conservation of mass, while true, is not sufficient by itself to prevent the time of existence shrinking to zero, thus creating finite time blowup.  Indeed, if one
considers a soliton solution $u(t,x) = Q(x) e^{i\omega t}$ and applies a pseudoconformal transformation \eqref{pc-conf} followed by time translation, one obtains the explicit solution
$$ \frac{1}{|t-1|^{d/2}} e^{i|x|^2/4(t-1)} e^{i\omega/(t-1)} Q( \frac{x}{t-1} ) $$
to this NLS.  This solution is smooth with finite mass at $t=0$, and remains smooth with conserved mass for $0 \leq t < 1$, but nevertheless develops a singularity at $t=1$ because the mass has concentrated to a point.  Basically, the scale-invariance of the equation has created a non-compactness in the phase space into which the dynamics can escape into in finite time.  To prevent this type of blowup one must thus exclude this type of \emph{mass concentration} or \emph{energy concentration} where the mass or energy is scaling itself into higher and higher frequencies in finite time.  To do this, conservation laws alone are not enough; one needs the additional tool of \emph{monotonicity formulae}, which we turn to in the next section.

We remark that it is possible, and very useful, to modify conserved quantities by inserting either spatial weights (e.g. cutoff functions) or frequency weights (e.g. derivatives or Littlewood-Paley projections) to create a much larger class of \emph{almost conserved quantities}, whose derivative is not quite zero, but is still somehow ``lower order'' than what one might naively expect.  To give a simple example, in the KdV equation\footnote{For this particular equation, which is completely integrable, one can find a quantity similar to $E_2$ which is exactly conserved.  However, the approach here is more robust, and in particular applies to variants of the KdV equation, such as the difference equation governing the difference to two solutions to KdV; because of this, the ``energy method'' we give here can be used, with some additional arguments, to give a simple local existence theorem in $H^2_x(\R)$ for KdV.  See \cite{bona}, \cite{kato-kdv}, \cite{kato-kdv-2}.}, the ``second energy''
$$ E_2(t) := \int_\R u(t,x)^2 + u_x(t,x)^2 + u_{xx}(t,x)^2\ dx,$$
which is essentially the standard energy weighted by a single derivative (or the mass weighted by two derivatives), is not a conserved quantity.  However some routine integration by parts (and Sobolev embedding) eventually yields the differential inequality
$$ \partial_t E_2(t) = O( E_2(t)^{3/2} )$$
which gives an elementary \emph{a priori} local estimate for the growth of the $H^2_x(\R)$ norm.  The point here is that the right-hand side only involves second derivatives of $u$ at most, whereas a naive inspection of the KdV equation might have suggested instead that as many as five derivatives of $u$ would have to be involved.

These almost conserved quantities can serve as more flexible substitutes for the usual conservation laws, being adaptable to situations where one only has local control of the mass or energy, or for which one is in a rougher or smoother Sobolev space than the mass or energy class.  For instance, the ``$I$-method'' for extending subcritical global existence results from the energy regularity to slightly rougher regularities, as employed for instance in \cite{ckstt:2}, is of this type.  These types of ``local'' almost conservation laws are important in both subcritical and critical equations in controlling how much mass and energy flows low frequencies to high, and from nearby locations to distant ones, or vice versa; see e.g. \cite{taobook} for some examples of this. For reasons of space, however, we will not discuss these techniques in detail here.

\section{Monotonicity formulae}

All the model equations here are examples of Hamiltonian PDE, and in particular are all time reversible.  Thus, in contrast to parabolic equations (such as the heat equation), there is no preferred direction of time.  Thus we do not expect behaviour such as the existence of compact attractors.  In the case of Hamiltonian \emph{ODE}, one has some additional results (e.g. Liouville's theorem on preservation of symplectic volume, Gromov's nonsqueezing theorem, or the Poincar\'e recurrence theorem) which further strengthen this intuition that a Hamiltonian flow cannot ``compress'' the dynamics of arbitrary data into that of a smaller set.

However, the situation can be remarkably different in the case of Hamiltonian \emph{PDE}, especially those on non-compact domains such as Euclidean space $\R^d$.  Here one encounters a phenomenon that while quantities such as energy and mass are conserved, they often radiate away to spatial infinity, so that the \emph{local} mass and energy in a compact region goes to zero
both as $t \to +\infty$ and as $t \to -\infty$.  This mechanism of \emph{dispersion} can serve as a weak substitute for the \emph{dissipation} mechanism for parabolic equations\footnote{Indeed, a useful (though not entirely accurate) rule of thumb is that dispersive models such as the ones studied here are, generally speaking, expected to have similar global existence and blowup properties to their parabolic counterparts; for instance, the theory for wave and Schr\"odinger maps should be roughly analogous to that of the harmonic map heat flow, the theory for NLS and NLW should be analaogous to that of the nonlinear heat equation, etc.  Indeed, the parabolic equations face many of the same key distinctions as the dispersive models, such as subcritical vs. supercritical energies, focusing vs. defocusing, etc.  Nevertheless the actual \emph{proof} of global existence or blowup tends to be quite different in the two settings.}; roughly speaking, the dispersive effect is expected to cause most of the infinite degrees of freedom in the PDE to radiate harmlessly away to spatial infinity, following the linear evolution, leaving only an ``essentially compact'' core of the phase space to evolve in a genuinely nonlinear manner.

Our understanding of this dispersive effect, especially as it pertains to large data over long periods of time, is not well understood in the focusing case, where there are portions of phase space which do not disperse, but instead lead to solitons or to blowup solutions.  However, in defocusing cases we now have a reasonably satisfactory mechanism to rigorously establish dispersion, by modifying the conservation laws of the preceding section to produce quantities which are monotone decreasing or increasing in time, rather than being constant in time.  The reason this can be used to establish dispersion is due to a simple fact (from the fundamental theorem of calculus): if a quantity is both monotone and bounded, then its derivative is absolutely integrable, and in particular decays (at least on average) as time goes to infinity.  This decay can then be combined with the Duhamel formula \eqref{duhamel} and perturbation techniques (e.g. Strichartz estimates) to obtain good control on the solution at infinity (basically, that the linear behaviour dominates the nonlinear behaviour for sufficiently large times).

Now we turn to the algebraic manipulations which create these monotonicity formulae.  For simplicity let us ignore all issues of smoothness and regularity that would be needed to justify the manipulations below; in practice, the rigorous justification can be achieved by standard regularisation or limiting arguments and will not be discussed here.

Monotonicity formulae are close cousins of conservation laws, and so it is not surprising that the stress-energy tensor $T_{\alpha \beta}$ is a rich source of such formulae.  Indeed, if $T_{\alpha \beta}$ is any rank-two tensor obeying the conservation laws \eqref{conserv-coord}, then on multiplying these laws against an arbitrary scalar weight\footnote{Here we are taking a ``spatial'' perspective, in which we decouple the roles of space and time; this is particularly useful for NLS.  For nonlinear wave equations it is more profitable to take a ``spacetime'' approach which we discuss shortly.  On the other hand, it is sometimes useful for NLS to consider weights $a$ which depend on time as well as space, see e.g. \cite{nak-2d}.}
 $a(x)$ or a vector weight $a_k(x)$ and integrating by parts, we obtain (formally, at least) the integral identities
\begin{align}
\partial_t \int_{\R^d} \stress^{00}(t,x) a(x)\ dx &= \int_{\R^d} \stress^{0j}(t,x) \partial_j a(x)\ dx \label{var0}\\
\partial_t \int_{\R^d} \stress^{k0}(t,x) a_k(x)\ dx &= \int_{\R^d} \stress^{kj}(t,x) \partial_j a_k(x)\ dx.\label{var1}
\end{align}
The first identity \eqref{var0} is thus a first variation formula for integrals of the energy or mass density $\stress^{00}$, and is particularly useful for understanding the local flux of such densities.  The second identity \eqref{var1} (which is a first variation formula for the momentum density) turns out to be particularly useful when the stress-energy tensor is symmetric, and $a_k= \partial_k a$ is a gradient vector field, in which case it becomes a second variation formula for the above integrals:
\begin{align}
\partial_{tt} \int_{\R^d} \stress^{00}(t,x) a(x)\ dx &= \partial_t 
\int_{\R^d} \stress^{0j}(t,x) \partial_j a(x)\ dx \label{second-1}\\
&= \int_{\R^d} \stress^{jk}(t,x) \partial_{jk} a(x)\ dx.\label{second-2}
\end{align}
We have complete freedom to choose the weight $a$.  It turns out that if this weight is sufficiently ``convex'' (so that $\partial_{jk} a$ is positive definite), the quantity \eqref{second-2} can be non-negative, thus leading to a monotonicity formula for the weighted momentum\footnote{Note that it is only the weighted momentum which has a chance to enjoy a monotonicity formula.  A weighted mass or weighted energy cannot be monotone in time as this would be inconsistent with time reversal symmetry; on the other hand, reversing time also reverses the momentum and so does not contradict a momentum monotonicity formula.  However, we see from \eqref{second-1}, \eqref{second-2} that a weighted mass or energy can be \emph{convex} in time.  These convexity formulae are known as \emph{virial identities} and play an important role in both focusing and defocusing equations.}
$$M_a(t) := \int_{\R^d} \stress^{0j}(t,x) \partial_j a(x)\ dx = \partial_t \int_{\R^d} \stress^{00}(t,x) a(x)\ dx.$$
If for instance we specialise to the NLS, then we have
$$ M_a(t) = 2\int_{\R^d} \Im( \overline{u} u_j)(t,x) \partial_j a(x)\ dx = \partial_t \int_{\R^d} |u(t,x)|^2 a(x)\ dx$$
and (after one last integration by parts)
\begin{align*}
 \partial_t M_a(t) &= \int_{\R^d} 4\Re(\overline{u_j} u_k)(t,x) \partial_{jk} a(x)\ dx \\
&\quad - \int_{\R^d} |u(t,x)|^2 \Delta\Delta a(x)\ dx \\
&\quad + \frac{4 \mu(p-1)}{p+1} \int_{\R^d} |u(t,x)|^{p+1} \Delta a(x)\ dx.
\end{align*}
Now suppose that $a$ is (non-strictly) convex, so that $\partial_{jk} a(x)$ is positive-definite; then the first term on the right-hand side (which is the top order term) is non-negative.  If
we also have $\mu \geq 0$ (so we either have a defocusing NLS, or the linear Schr\"odinger equation), and we also have the sub-biharmonic condition\footnote{This term arises from ``quantum corrections'' to the classical analogue of this formula, which asserts that if a particle $t \mapsto x(t)$ evolves by Newton's first law $\partial_{tt} x(t) = 0$, then the weighted momentum $M_a(t) := \partial_t x_j(t) \partial_j a( x(t) ) 
= \partial_t a(x(t))$ evolves by the formula $\partial_t M_a(t) = \partial_t x_j(t) \partial_t x_k(t) \partial_{jk} a(x(t))$.  In general, while any monotonicity formula for the Schr\"odinger equation must necessarily imply a classical monotonicity formula for Newtonian particle motion (by taking the \emph{semiclassical limit} $\hbar \to 0$), the converse is not always true, unless one is only interested in top order terms.}
 $-\Delta \Delta \geq 0$, then the two lower order terms are also non-negative, and so we have a genuine monotonicity formula.  This is the ideal situation; however, even when some of the lower order terms have no preferred sign, one can often still extract nontrivial control on the solution as long as the top order term is mostly positive.  Thus it is really the convexity of $a$ which leads to important formulae.
 
Let us give some basic examples of this formula in action.  Setting $a := 1$ simply gives conservation of mass.  Setting
$a := x$ (or $a \equiv x_k$ for $k=1,\ldots,d$) simply gives conservation of the total momentum \eqref{total-momentum}, and also reveals that the (un-normalised) centre-of-mass $\int_{\R^d} x |u(t,x)|^2\ dx$ varies linearly in time, with rate of change equal to the total momentum.  
Setting $a := |x|^2$ gives rise to \emph{Glassey's virial identity}\cite{glassey-blow}
\begin{equation}\label{virial}
\begin{split}
\partial_{tt} \int_{\R^d} |x|^2 |u(t,x)|^2\ dx &= 2\int_{\R^d} x_j \Im( \overline{u} u_j)(t,x)\ dx \\
&= 8 \int_{\R^d} |\nabla u(t,x)|^2\ dx + \frac{8 \mu (p-1) d}{p+1} |u(t,x)|^{p+1}\ dx.
\end{split}
\end{equation}
For simplicity let us consider the pseudoconformal case $p =1 + \frac{4}{d}$, in which the virial identity takes the particularly appealing form
$$ \partial_{tt} \int_{\R^d} |x|^2 |u(t,x)|^2\ dx = 16 E(u),$$
thus the second variation of the (un-normalised) mass variance $\int_{\R^d} |x|^2 |u(t,x)|^2\ dx$ is essentially equal to the
conserved energy.  This variance can be viewed as a measure of how close the mass clusters to the origin; thus when the energy is positive, we expect the mass to be repelled from the origin, while when the energy is negative (which can happen in the focusing case $\mu=-1$), the mass should be attracted to the origin.  (For stationary solitons in the pseudoconformal NLS, the energy is precisely zero; this is a special case of the \emph{Pohozaev identity}.)  One consequence of this is that when the energy is negative (and assuming suitable decay and regularity conditions on the initial data), then the solution to NLS must blow up in finite time (both for positive and negative times); see \cite{glassey-blow}.  Similar results hold for higher (mass-supercritical) powers.  In the defocusing case $\mu=+1$ and with arbitrary power $p$ (or in the focusing case $\mu=-1$ and mass-subcritical power), we obtain the inequality
$$ \partial_{tt} \int_{\R^d} |x|^2 |u(t,x)|^2\ dx \geq c_{p,d} E(u)$$
for some positive constant $c_{p,d} > 0$ depending only on $p,d$.  If the energy is strictly positive, this implies that $\int_{\R^2} |x|^2 |u(t,x)|^2\ dx$ goes to infinity as $t \to \pm \infty$; thus the solution cannot stay strongly localised near the origin indefinitely.  This statement is not always directly useful, because it requires a lot of decay on the solution $u$ (in particular, that $xu$ is square integrable) but in practice one can modify the above argument by smoothly truncating the weight $a(x) = |x|^2$ smoothly at infinity and dealing somehow with the error terms.  There are many instances of this trick in the literature; see e.g. \cite{merlekenig} for a very recent one.

An alternative monotonicity formula, which is especially useful in the defocusing case is the \emph{Morawetz inequality} of Lin and Strauss \cite{linstrauss}, which is obtained by setting $a(x) = |x|$.  It is geometrically obvious that $a$ is non-strictly convex.  For sake of discussion let us specialise to three dimensions $d=3$, to defocusing nonlinearities $\mu=+1$, to finite energy and mass solutions, and to energy critical or sub-critical nonlinearities $p \leq 5$ (in order to be able to use Theorem \ref{nls-subcrit-lwp} or Theorem \ref{nls-crit-lwp}).  In this setting we have $-\Delta\Delta a = 8 \pi \delta$ in the sense of distributions, where $\delta$ is the Dirac mass; thus (after doing some standard arguments to handle the singularity of $a$ and its derivatives at the origin) we obtain
$$ M_a(t) = 2\int_{\R^d} \Im( \overline{u} u_j)(t,x) \frac{x_j}{|x|}\ dx$$
and
\begin{align*}
 \partial_t M_a(t) &= \int_{\R^d} 4 \frac{|\nabb u(t,x)|^2}{|x|}\ dx \\
&\quad + 8 \pi |u(t,0)|^2\\
&\quad + \frac{8 (p-1)}{p+1} \int_{\R^d} \frac{|u(t,x)|^{p+1}}{|x|}\ dx
\end{align*}
where $\nabb u$ is the angular component of the gradient, thus
$$ |\nabb u(t,x)|^2 = |\nabla u(t,x)|^2 - |\frac{x}{|x|} \cdot \nabla u(t,x)|^2.$$
In particular, in the defocusing case $\mu=+1$ in three dimensions, we have the monotonicity formula\footnote{Physically, $M_a(t)$ represents the radially outward momentum; the portion of the momentum which is radiating away of the origin.  As time progresses, inward momentum gets converted into outward momentum, but not vice versa, thus explaning the monotonicity.  The nonlinear factor $\int_{\R^d} \frac{|u(t,x)|^{p+1}}{|x|}\ dx$ represents the fact that the defocusing nonlinearity also converts inward momentum to outward momentum, but not vice versa.}
$$ \partial_t M_a(t) \geq c_p \int_{\R^d} \frac{|u(t,x)|^{p+1}}{|x|}\ dx \geq 0$$
for some absolute constant $c_p > 0$.  On the other hand, from the Cauchy-Schwarz inequality and conservation of mass and energy we have
the upper bound
$$ |M_a(t)| \leq 2 M(u)^{1/4} E(u)^{1/4}.$$
From the fundamental theorem of calculus we thus obtain the global spacetime bound
\begin{equation}\label{m1}
\int_I \int_{\R^d}\frac{|u(t,x)|^{p+1}}{|x|}\ dx dt \leq C_p M(u)^{1/4} E(u)^{1/4}
\end{equation}
where $I$ is the maximal interval of existence.  Note that the right-hand side does not depend on the size of $I$ (which in fact turns out to be infinite), and also does not require any decay on the solution other than finite mass and energy.  If $I$ is infinite, then this estimate shows that the quantity $\frac{|u(t,x)|^{p+1}}{|x|}$ is globally integrable in spacetime, and in particular decays in some suitable norm as $t \to \pm \infty$.  

One drawback of the above Morawetz estimate is the presence of the $\frac{1}{|x|}$ weight, which means that the estimate is strong near the spatial origin $x=0$ and weak away from this origin.  For the class of spherically symmetric solutions, one can use the radial Sobolev inequality
$$ |f(x)| \leq C \min( \frac{1}{|x|}, \frac{1}{|x|^{1/2}} ) \| f \|_{H^1_x(\R^3)} \hbox{ for all } x \in \R^3 \backslash \{0\},$$
which localises the finite-energy function $u(t)$ to near the origin, to effectively exploit the Morawetz estimate \eqref{m1}.  However, for solutions in translation-invariant classes such as the energy class without any symmetry assumption, the estimate \eqref{m1} can be arbitrarily weak and thus will not be able by itself to establish translation-invariant control on such solutions.  There is however an interesting ``doubling'' trick that can get around this difficulty, by introducing \emph{two} spatial variables $x,y \in \R^d$ instead of one.  Indeed, a routine modification of the second variation formula \eqref{second-1}, \eqref{second-2} yields the ``two-particle'' variant
\begin{align}
\partial_{tt} \int_{\R^d} \stress^{00}(t,x)& \stress^{00}(t,y) a(x-y)\ dx dy = 2\partial_t 
\int_{\R^d} \stress^{0j}(t,x) \stress^{00}(t,y) \partial_j a(x-y)\ dx dy \label{third-1}\\
&= 2 \int_{\R^d} [\stress^{jk}(t,x) \stress^{00}(t,y) - \stress^{j0}(t,x) \stress^{k0}(t,y)] \partial_{jk} a(x-y)\ dx dy\label{third-2}
\end{align}
whenever $a$ is an even function.  In the case of the NLS, we obtain the identity
\begin{align*}
 \partial_t M_{a,2}(t) 
 &= 2 \int_{\R^d} \Re(\overline{p_{j}(t,x,y)} p_{k}(t,x,y)) \partial_{jk} a(x-y)\ dx \\
 &\quad 2 \int_{\R^d} \Re(\overline{q_{j}(t,x,y)} q_{k}(t,x,y)) \partial_{jk} a(x-y)\ dx \\
&\quad - 2\int_{\R^d} |u(t,x)|^2 |u(t,y)|^2 \Delta\Delta a(x-y)\ dx dy \\
&\quad + \frac{8 \mu(p-1)}{p+1} \int_{\R^d} |u(t,x)|^{p+1} |u(t,y)|^2 \Delta a(x-y)\ dx
\end{align*}
where 
$$ M_{a,2}(t) := 4\int_{\R^d} \Im( \overline{u} u_j)(t,x)  |u(t,y)|^2 \partial_j a(x-y)\ dx dy = \partial_t \int_{\R^d} |u(t,x)|^2 |u(t,y)|^2 a(x-y)\ dx dy$$
and
\begin{align*}
p_{j}(t,x,y) &:= u(t,x) u_j(t,y) - u_j(t,x) u(t,y) \\
q_{j}(t,x,y) &:= u(t,x) \overline{u_j}(t,y) + u_j(t,x) \overline{u}(t,y).
\end{align*}
In particular, if $a$ is non-strictly convex then
$$ \partial_t M_{a,2}(t) \geq - 2\int_{\R^d} |u(t,x)|^2 |u(t,y)|^2 \Delta\Delta a(x-y)\ dx dy.$$
Setting $a(x) := |x|$, $\mu=1$ and $d=3$ as before, we conclude that
$$ \partial_t M_{a,2}(t) \geq c \int_{\R^3} |u(t,x)|^4\ dx$$
for some absolute constant $c > 0$, which eventually leads to the spacetime bound
\begin{equation}\label{iu}
 \int_I \int_{\R^3}|u(t,x)|^4\ dx dt \leq C M(u)^{3/2} E(u)^{1/2}
\end{equation}
for energy-class solutions to the NLS with $\mu=1, d=3, p \leq 5$ and $I$ the maximal interval of existence.  This ``interaction'' or ``two-particle'' Morawetz inequality is similar to the ``one-particle'' Morawetz inequality
\eqref{m1}, but now does not have the weight $\frac{1}{|x|}$ and is now better suited for translation-invariant situations.
These spacetime bounds can be inserted (possibly after combining them with other spacetime bounds, such as those arising from mass and energy conservation or from the Duhamel formula \eqref{duhamel}) into the scattering criterion in theorems such as Theorem \ref{nls-subcrit-lwp}, for instance giving a fairly quick proof of scattering in the energy class in the regime $\mu=1, d=3, 3 < p < 5$ (a result first obtained in \cite{gv:scatter}); see \cite{ckstt:french}.

Analogous monotonicity formulae exist for nonlinear wave equations (although finding good analogues of the interaction Morawetz inequality for such equations has proven surprisingly elusive, except in one dimension when they correspond to the classical Glimm interaction estimates).  Here it is more natural geometrically (and physically) to treat spacetime as a unified object (Minkowski spacetime).  Again we work formally, ignoring issues of regularity or integrability.  From the conservation law \eqref{conserv} we have the divergence identity
$$ \partial_\alpha( T^{\alpha \beta} X_\beta ) = \frac{1}{2} T^{\alpha \beta} \pi_{\alpha \beta}$$
for any vector field $X = X^\alpha$, where $\pi = \pi_{\alpha \beta}$ is the \emph{deformation tensor}
$$ \pi^{\alpha \beta} := \nabla^\alpha X^\beta + \nabla^\beta X^\alpha$$
(or $\pi^{\alpha \beta} = \partial^\alpha X^\beta + \partial^\beta X^\alpha$ in the usual Minkowski coordinate system).  This identity is particularly simple when $X$ is a \emph{Killing vector field} (i.e. an infinitesimal isometry of Minkowski space), since in this case the deformation tensor vanishes, and we obtain a conserved current $T^{\alpha \beta} X_\beta$.  However, the number of linearly independent Killing vector fields is very small (basically one only obtains the conservation of energy, momentum, and energy momentum this way).  One can often also extract conserved (or almost conserved) currents from \emph{conformal} Killing vector fields (such as the scaling vector field $t \partial_t + x_j \partial_j$ or the \emph{Morawetz vector field} $(t^2+x^2) \partial_t + 2 t x_j \partial_j$), in which the deformation tensor $\pi_{\alpha \beta}$ is a scalar multiple of the metric $g_{\alpha \beta}$, basically because the trace $T^{\alpha \beta} g_{\alpha \beta}$ of the stress-energy tensor is often either zero, or is itself the divergence of another vector field.  For instance, using the scaling vector field
$t \partial_t + x_j \partial_j$ in the energy-critical defocusing case $\mu=+1$, $d=3$, $p=5$ and Stokes' theorem, combined with some additional arguments, one can obtain a non-concentration property for the potential energy density:
\begin{equation}\label{6decay}
\lim_{t \to 0} \int_{|x| \leq t} |u(t,x)|^6\ dx = 0;
\end{equation}
see e.g. \cite{shatah-struwe}.  When combined with finite speed of propagation and perturbative analysis (based on Strichartz estimates), one can use \eqref{6decay} to establish global regularity (or well-posedness in the energy class) for this equation; the point is that \eqref{6decay} shows that even large energy data will behave like small (potential) energy locally in spacetime, at which point the perturbative theory can be used to show that no blowup can occur.  

It is also useful to consider other types of vector fields than conformal Killing vector fields.  As was the case with NLS, it is profitable to consider vector fields which are gradients of some scalar potential $a$, thus $X^\alpha = \partial^\alpha a$, and we obtain
$$ \partial_\alpha( T^{\alpha \beta} \partial_\beta a ) = T^{\alpha \beta} \partial_{\alpha \beta} a;$$
in the specific case of NLW, the right-hand side becomes
$$ (\partial^\alpha u) (\partial^\beta u) \partial_{\alpha \beta} a
- \frac{1}{4} (\Box a) \Box(|u|^2) + \frac{\mu(p-1)}{p+1} |u|^{p+1} \Box a.$$
Once again, one can often obtain a useful monotonicity formula from the case when $a$ is non-strictly convex.  For instance,
with the same equation $\mu=+1, d=3, p=5$ as before, one can use the weight $a(t,x) = |x|$ to obtain the Morawetz inequality
$$ \int_I \int_{\R^3} \frac{|u(t,x)|^6}{|x|}\ dx dt \leq C E(u)$$
for some absolute constant $C$ (compare with \eqref{m1}).  This can be used as a substitute for \eqref{6decay} for the purposes of establishing global regularity or scattering.

A variety of monotonicity formulae are also known for wave maps, especially in the presence of symmetry; see \cite{shatah-struwe}, \cite{tao:forges}, \cite{taobook}.  Generally speaking, these formulae assert that as one approaches a potential singularity of a wave map, that the (rescaled) wave map converges (in some weak sense) to a harmonic map.  It would be of interest to make this phenomenon more quantitative, as this would undoubtedly be useful in the (still incomplete) theory of large energy wave maps.  For the Maxwell-Klein-Gordon and Yang-Mills equations, no nontrivial monotonocity formulae appear to be known.

It would be of interest to obtain further monotonicity formulae which are not so dependent on the stress-energy tensor.  One tentative step in this direction is in \cite{tao-gkdv}, in which the mass and energy conservation laws for (defocusing) gKdV are played off against each other to obtain a dispersion estimate.

\section{Induction on energy}\label{energy-sec}

Historically, the first large data global regularity result for a critical nonlinear dispersive or wave equation was that for the defocusing energy-critical NLW in three dimensions ($\mu=+1$, $d=3$, $p=5$); see \cite{struwe}, \cite{g_waveI}, \cite{grillakis.semilinear}, \cite{struweshatah.energy}, \cite{shatah-struwe}.  The approach (which was inspired by some similar arguments in nonlinear elliptic and parabolic equations) was based upon two basic ingredients:

\begin{itemize}
\item (Small energy implies regularity) If the energy is sufficiently small, then no singularities can form; this follows from perturbative analysis.  In practice, one needs stronger versions of this statement, in which only the potential energy is assumed to be (locally) small.
\item (Nonconcentration of energy) The (potential) energy is shown to locally decay as one approaches any given point in spacetime.  This is non-perturbative and is achieved by a monotonicity formula approach (e.g. Morawetz estimates).
\end{itemize}

This two-step approach then formed the model for a number of other critical global regularity results, such as those for radially or equivariantly symmetric critical wave maps \cite{christ.spherical.wave}, \cite{christ2}, \cite{shatah.shadi.eqvt}
or Yang-Mills-Higgs \cite{keel}.  A crucial feature of these equations was that the quantity which was shown to decay by a monotonicity formula was critical (scale-invariant); otherwise, there was no chance that smallness of this quantity would be at all helpful for establishing regularity.  For nonlinear wave equations, this type of scale-invariance was achievable, ultimately because the momentum density (which was the source of monotonicity formulae) had the same scaling as the energy density (which was already assumed to be critical).

The energy-critical defocusing NLS (e.g. $\mu=+1$, $d=3$, $p=5$) thus presented a new difficulty, because the momentum and energy no longer had the same scaling, and so the known monotonicity formulae (such as \eqref{m1}) did not establish decay of any useful critical quantity near a potential singularity.  This difficulty was resolved by Bourgain \cite{bourg.critical} and Grillakis \cite{grillakis:scatter} in the case of spherical symmetry, and later by Colliander-Keel-Staffilani-Takaoka-Tao \cite{gopher} in the general case, based on a number of additional observations.  The first was that a non-critical monotonicity formula such as \eqref{m1} could be localised via cutoff functions to obtain a critical estimate, albeit one which now depended on the scale of the cutoff.  For instance, by smoothly truncating the weight $a(x) = |x|$ to a ball centred at the origin, one can modify \eqref{m1} to the estimate
\begin{equation}\label{moraw}
\frac{1}{|J|^{1/2}} \int_J \int_{|x| \leq K |J|^{1/2}} \frac{|u(t,x)|^6}{|x|}\ dx dt \leq C_K E(u)
\end{equation}
for all intervals $J$ inside the maximal interval of existence and all $K > 0$, where the constant $C_K$ depends on $K$; see \cite{bourg.critical}, \cite{grillakis:scatter}.  The point is that the right-hand side only involves the critical energy $E(u)$ and not the supercritical mass $M(u)$; indeed, both sides of this inequality are scale-invariant.  The drawback to this estimate was the unusual nature of the left-hand side, in particular the presence of the weight $\frac{1}{|J|^{1/2}}$.  This made it difficult to convert this type of scale-invariant control to an estimate which could be used as input for the perturbation theory (which would require a critical unweighted spacetime norm, such as the $L^{10}_{t,x}$ norm of the solution).  The basic problem is that any two given norms on the solution need not be comparable, even after insisting that both norms are critical; there is a serious \emph{lack of compactness} in the space of solutions that is not resolved simply by quotienting out by symmetries such as scale invariance.

A key breakthrough\footnote{This method is not strictly necessary for the energy-critical NLS; see \cite{grillakis:scatter}, \cite{tao-nls} for some alternate approaches.  However, the induction on energy philosophy seems to provide a powerful and unified tool to approach many other critical problems, decreasing the need to rely on more ad hoc methods.}
was made by Bourgain \cite{bourg.critical}, who introduced an \emph{induction on energy} method which ``compactified'' the dynamics of solution sufficiently that one could begin comparing different (but critical) norms on the solution.  The method is closely related, though not identical, to the \emph{concentration compactness} method of Lions; we compare the two methods in the next section.  

The induction on energy method is the analogue of the energy minimisation method used to construct solutions of elliptic equations (for instance, minimising the Dirichlet energy to solve the Dirichlet problem).  The fact that one works with \emph{minimisers} of a functional, rather than merely \emph{critical points}, can allow one to restrict the solution\footnote{This is of essentially the famous \emph{Palais-Smale condition} for variational functionals.} to a compact set (perhaps after quotienting out by the symmetries of the problem); in practical terms, this means that the minimiser ``behaves like a bump function'' in the sense that it is localised in both space and frequency.

We illustrate this technique with the energy-critical defocusing three-dimensional NLS (so $\mu=+1$, $d=3$, $p=5$), though the method is quite general and has been extended to several other equations.  
The main result here is the following \emph{a priori} estimate:

\begin{theorem}\label{l10-bound}\cite{bourg.critical}, \cite{gopher}  Let $u \in C^0_t \dot H^1_x(I \times \R^3)$ be an energy-class solution to NLS with $\mu=+1, d=3, p=5$ on a compact time interval $I$ with energy $E(u) \leq E$.  Then we have the bound
$$ \int_I \int_{\R^3} |u(t,x)|^{10}\ dx dt \leq A(E)$$
for some finite quantity $A(E)$ depending only on $E$.
\end{theorem}

From this theorem and Theorem \ref{nls-crit-lwp} one obtains

\begin{corollary} Let $u_0 \in \dot H^1_x(\R^3)$.  Then there is a unique global energy-class solution $u \in C^0_t \dot H^1_x(\R \to \R^3)$ to NLS with $\mu=+1, d=3, p=5$ and $u(0) = u_0$, which also lies in the space $L^{10}_{t,x}(\R \times \R^3)$.  Also, $u$ scatters to a linear solution $e^{it\Delta} u_\pm$ as $t \to \pm \infty$ for some $u_\pm \in H^1_x(\R^3)$, and if $u_0$ is Schwartz then $u$ will be Schwartz in space and smooth in time.
\end{corollary}

Thus the main task is to establish Theorem \ref{l10-bound}.  We introduce the function $A: \R \to [0,+\infty]$ by
$$ A(E) := \sup \{ \int_I \int_{\R^3} |u(t,x)|^{10}\ dx dt : E(u) \leq E \}$$
where the supremum ranges over all energy-class solutions $u$ to NLS of energy at most $E$, with the convention that $A(E) = 0$ when $E$ is negative.  The task is to establish that $A(E)$ is finite for all $E$; note that Theorem \ref{nls-crit-lwp} already gives this for small $E$.

The basic induction-on-energy strategy of Bourgain \cite{bourg.critical} is to establish this finiteness by estimating $A(E)$ in terms of $A(E')$ for various explicit smaller energies $E' \leq E$.  In particular, when restricting to spherically symmetric solutions (thus decreasing $A(E)$), the recursive inequality
\begin{equation}\label{ae}
 A(E) \leq C \exp( C \eta^{-C} A( E - \eta^4 )^C )
\end{equation}
was proven for all energies $E \geq C^{-1}$, where $C$ is an absolute constant and $\eta$ was a small quantity depending on $E$ (one can take $\eta := 1 / (C E^C)$).  Very briefly, this type of inequality was obtained by first performing some lengthy analysis (both perturbative and non-perturbative)
to argue that if a solution with energy $E$ had very large $L^{10}_{t,x}$ norm, then at some time the solution must decouple into an isolated ``bubble'' (of energy comparable to some power of $\eta$), together with a remainder component of energy at most $E - \eta^4$.  By inductive hypothesis, the remainder would evolve with an $L^{10}_{t,x}$ norm controlled by $A(E-\eta^4)$.  One then applies stability theory (such as Theorem \ref{long-time theorem}), combined with the isolation property, to then control the $L^{10}_{t,x}$ norm of the original solution.

From iterating \eqref{ae} it is not difficult to show that $A(E)$ is finite for all $E$, although the upper bound obtained in $A()$ is rather poor (it is a tower of exponentials of height polynomial in $E$).

In \cite{gopher} the induction-on-energy strategy was reinterpreted as an analysis of minimal-energy blowup solutions, in analogy to the method of mathematical induction can often be reinterpreted in the contrapositive as the method of descent.  
It is not hard to show (using Theorem \ref{nls-crit-lwp}) that $A$ is monotone non-decreasing, left-continuous, and finite for small $E$.  From this we obtain a dichotomy: either $A(E)$ is finite for all $E$, or else there exists a \emph{critical energy} $0 < E_{\crit} < \infty$ with the property
that $A(E) < \infty$ for all $E < E_\crit$ and $A(E) = +\infty$ for all $E \geq E_\crit$.  Thus $E_\crit$ is the minimal energy required for the solution to blow up in the sense that the $L^{10}_{t,x}$ norm becomes infinite (which is a natural criterion for blow up, in light of Theorem \ref{nls-crit-lwp}).  Thus to show that $A(E)$ is finite for all $E$, we may assume for contradiction that a finite critical energy $E_\crit$ exists, and then obtain a contradiction.  Note for instance that a bound such as \eqref{ae} can achieve this.  One advantage of this formulation is that it allows one to exploit a \emph{system} of inequalities connecting $A()$ with various quantities such as $\eta$, as opposed to just a single inequality; such systems often require a multiple induction if one wanted to apply them directly.  Conversely, if one uses the minimal-energy blowup formulation it is quite difficult to establish any explicit bounds on $A(E)$ other than that it is finite.
For instance, if one established the conditional inequality
$$ A(E) \leq C_0( E, \eta_0, \eta_1, A(E-\eta_0), A(E-\eta_1) ) \hbox{ whenever } 1/\eta_1 \leq C_1( E, \eta_0, A(E-\eta_0) ) \hbox{ and } \eta_0 \leq 1/(CE^C)$$
where $C_0(), C_1()$ denote various explicit functions,
then it is easy to see that this is inconsistent with the existence of a finite critical energy $E_\crit$, although to establish the finiteness of $A(E)$ directly from this inequality requires a double induction.  

Suppose that the critical energy $E_\crit$ was finite.  Then we can find solutions $u: I \times \R^3 \to \C$ of energy $E(u) \leq E_\crit$ whose $L^{10}_{t,x}$ norm is arbitrary large.  In fact, it turns out (by the concentration compactness arguments below, see \cite{keraani}) that we can find a maximal-lifespan solution $u: I \times \R^3 \to \C$ of energy exactly $E(u) = E_\crit$ whose $L^{10}_{t,x}$ norm is \emph{infinite}; in fact with a little refinement (see \cite{tvz-compact}) we can ensure $L^{10}_{t,x}$ blowup in both directions, thus the $L^{10}_{t,x}$ norm is infinite on both $(I \cap [t_0,+\infty)) \times \R^3$
and $(I \cap (-\infty,t_0]) \times \R^3$ for any $t_0 \in I$.
we refer to such solutions as \emph{minimal-energy blowup solutions}.  For the purposes of the induction-on-energy argument, it is not strictly necessary to work with minimal-energy blowup solutions, and one can instead work with \emph{almost-blowing-up} solutions of nearly the minimal energy, in which the $L^{10}_{t,x}$ norm is very large rather than infinite (see e.g. \cite{gopher}), but we shall use exactly minimal-energy blowup solutions as they are conceptually and technically simpler to deal with.

\begin{remark} In the focusing case $\mu = -1$, there is a smooth non-negative stationary solution $u(t,x) = W(x)$, where $\Delta W = - W^5$ (in fact we have the explicit formula $W(x) := 1 / (1 + |x|^2/3)^{1/2}$).  This solution exists globally, but blows up in the sense that its $L^{10}_{t,x}$ norm is infinite.  Thus in the focusing case, the analogue of the critical energy is at most $E(W)$.  It is conjectured in the focusing case that the critical energy is in fact \emph{exactly} $E(W)$, thus any solution with energy (and $\dot H^1$ norm) less than that of the stationary solution should exist globally with finite $L^{10}_{t,x}$ norm; then the stationary solution $u(t,x) = W(x)$ would become a minimal-energy blowup solution.  This conjecture has recently been verified in the spherically symmetric case \cite{merlekenig}.
\end{remark}

The key advantage of working with minimal-energy blowup solutions, as opposed to more general solutions, lies by exploiting the following informal principle:

\begin{centering}
Minimal-energy blowup solutions are \emph{irreducible} and hence \emph{localised}.  In fact they are \emph{almost periodic} modulo symmetries.
\end{centering}

Readers who are familiar with elliptic variational theory may see an analogy here between minimal-energy blowup solutions and energy-minimisers of various elliptic functionals, such as the Dirichlet energy functional.  Thus the induction-on-energy method can be viewed as an analogue of the variational method for evolution equations.

Let us now explain some of the terms in the above principle more precisely.  By \emph{irreducible}, we mean that a minimal-energy blowup solution cannot ever decompose into the sum of two weakly interacting components of non-trivial energy.  For, if this were the case, each of the components would have strictly smaller energy than the critical energy $E_\crit$, and hence they each evolve separately by NLS with bounded $L^{10}_{t,x}$ norm.  Because the NLS equation is not completely linear, the superposition (sum) of these two evolutions is not quite a solution to NLS. However, if the interaction between the two components is sufficiently weak, then this superposition will \emph{approximately} solve the NLS equation, with an accuracy which is sufficient for the stability theory (Theorem \ref{long-time theorem}) to be applicable.  This allows us to establish an $L^{10}_{t,x}$ bound on the original solution, contradicting the blowup hypothesis (i.e. that the $L^{10}_{t,x}$ norm is infinite).  We illustrate this informal strategy by sketching a proof of \emph{frequency irreducibility} from \cite[Proposition 4.3]{gopher}:

\begin{proposition}[Minimal-energy blowup solutions are frequency-irreducible]\label{freq-irreducible}\cite{gopher}  Let $u: I \times \R^3 \to \C$ be a solution to NLS (with $\mu = +1$, $p=5$, $d=3$).  Suppose that we have a time $t_0 \in I$, a frequency $N > 0$, and $\eta, K > 0$ such that we have the frequency separation property
$$ \| P_{\leq N} u(t_0) \|_{\dot H^1} \leq \eta$$
and
$$ \| P_{\geq K N} u(t_0) \|_{\dot H^1} \leq \eta$$
where $P_{\leq N}$ is a Littlewood-Paley frequency projection\footnote{The exact definition of Littlewood-Paley projection will not be important for this informal discussion.} to low frequencies $\{ \xi: |\xi| \leq N \}$, and $P_{\geq KN}$ is 
a Littlewood-Paley projection to high frequencies $\{ \xi: |\xi| \geq KN \}$.  Then, if $K$ is sufficiently large depending on $\eta$, then $u$ \emph{cannot} be a minimal-energy blowup solution.
\end{proposition}

\begin{proof}(Sketch) Let $u, t_0, N, \eta, K$ be as above; suppose for contradiction that $u$ is a minimal-energy blowup solution, so in particular $E(u) = E_\crit$.  We first invoke a useful pigeonholing trick to locate a suitably ``empty'' region of frequency space in which to split the solution.

Let $0 < \eta' \ll \eta$ be a small quantity, and $K' \gg 1$ be a large quantity.  If $K$ is sufficiently large depending on $\eta', K'$, then by the pigeonhole principle one can find $N'$ between $N$ and $KN/K'$ such that
$$ \| P_{N' < \cdot < K'N'} u(t_0) \|_{\dot H^1} \leq \eta',$$
where $P_{N' < \cdot < K'N'}$ is a Littlewood-Paley projection to frequencies $\{ \xi: N' < |\xi| < K'N' \}$.  We can then split
$$ u(t_0) = u_{\lo}(t_0) + u_{\hi}(t_0) + O_{\dot H^1}(\eta'),$$
where $u_\lo(t_0) := P_{\leq N'} u(t_0)$ and $u_\hi(t_0) := P_{\geq K'N'} u(t_0)$ are the low and high frequency components of $u(t_0)$, and $O_{\dot H^1}(\eta')$ is an error whose $\dot H^1$ norm is $O(\eta')$.  By hypothesis we see that $u_\lo(t_0)$ and $u_\hi(t_0)$ both have an $\dot H^1$ norm of at least $\eta$, and from this it is not too difficult (from orthogonality arguments, assuming $\eta'$ small and $K'$ large) that $u_\lo(t_0)$ and $u_\hi(t_0)$ have energy strictly less than $E_\crit$; more precisely one has
$$ E(u_\lo(t_0)), E(u_\hi(t_0)) \leq E_\crit - c(\eta)$$
for some $c(\eta) > 0$ depending only on $\eta$.  By induction hypothesis, we thus see that we may evolve $u_\lo$ and $u_\hi$ by NLS to create global solutions $u_\lo, u_\hi: \R \times \R^3 \to \C$ with bounded $L^{10}_{t,x}$ norm:
$$ \| u_\lo \|_{L^{10}_{t,x}(\R \times \R^3)},  \| u_\hi \|_{L^{10}_{t,x}(\R \times \R^3)} \leq A( E_\crit - c(\eta) ).$$
In particular, the scalar field $\tilde u := u_\lo + u_\hi$ has bounded $L^{10}_{t,x}$ norm on $\R \times \R^3$.

Now we compare $\tilde u$ and $u$.  At time $t_0$, the two fields only differ in $\dot H^1$ norm by $O(\eta')$, by construction.  Now at later times, the field $\tilde u$ does not quite solve NLS; instead, it solves the equation
$$ (i \partial_t + \Delta) \tilde u = |\tilde u|^4 \tilde u + e$$
where
$$ e := |u_\lo|^4 u_\lo + |u_\hi|^4 u_\hi - (|u_\lo + u_\hi|^4)(u_\lo + u_\hi) = O( |u_\lo| |u_\hi|^4 + |u_\hi| |u_\lo|^4 ).$$
One can show (with some effort) that $e$ is quite small in appropriate norms.  Roughly speaking, the reason is that $u_\hi$ and $u_\lo$ are widely separated in frequency at time $t_0$, and hence (by perturbation theory and the $L^{10}_{t,x}$ bounds) will also be essentially widely separated in frequency at all other times also.  It turns out (due to certain ``bilinear Strichartz estimates'', which ultimately stems from the basic dispersive fact that different frequencies propagate at different velocities) that the interaction of two linear solutions to the Schr\"odinger equation with widely different frequencies will be quite small.  The $L^{10}_{t,x}$ bounds ensure that the solutions $u_\hi$, $u_\lo$ behave somewhat linearly (at least at short times), and it is possible (by choosing $K'$ sufficiently large) to ensure that the interaction is indeed small; for details see \cite{gopher}.  If $\eta'$ is also
sufficiently small, Theorem \ref{long-time theorem} now applies, and we pass from $L^{10}_{t,x}$ control of the approximate solution $\tilde u$ to $L^{10}_{t,x}$ control of the exact solution $u$.  But this implies that $u$ cannot be a minimal-energy blowup solution, and the claim follows.
\end{proof}

By applying the above proposition (in the contrapositive) for all values of $\eta$ at once, it is not difficult to conclude:

\begin{corollary}[Minimal-energy blowup solutions are frequency-localised]\label{freqlocal}  Let $u: I \times \R^3 \to \C$ be a minimal-energy blowup solution to NLS.  Then there exists a function $N: I \to \R^+$, and for every $\eta > 0$ there exists $K(\eta) > 0$ such that
$$ \| P_{\leq N(t_0)/K(\eta)} u(t_0) \|_{\dot H^1} \leq \eta$$
and
$$ \| P_{\geq K(\eta) N(t_0)} u(t_0) \|_{\dot H^1} \leq \eta.$$
\end{corollary}

Indeed one can select $N(t_0)$ to be (say) the median frequency of the $\dot H^1$ energy distribution.  A similar (but more intricate) argument can also be employed to obtain spatial concentration:
 
\begin{proposition}[Minimal-energy blowup solutions are spatially-localised]\label{spatloc}\cite{gopher}  Let $u: I \times \R^3 \to \C$ be a minimal-energy blowup solution to NLS, and let $N: I \to \R^+$ be as above.  Then there exists $x: I \to \R^3$, and for every $\eta > 0$ there exists $K(\eta) > 0$ such that
$$ \int_{|x-x(t_0)| \geq K(\eta) / N(t_0) } |\nabla u(t_0,x)|^2\ dx \leq \eta$$
for all $t_0 \in I$.
\end{proposition}

\begin{proof}(Sketch) The first step is to establish the weaker property of spatial \emph{concentration} of energy, namely that there exists an $x(t_0) \in \R^3$ for each $t_0 \in I$ such that
$$ \int_{|x-x(t_0)| = O(1/N(t_0)) } |\nabla u(t_0,x)|^2\ dx \geq c > 0$$
for some $c > 0$ depending only on the critical energy $E_\crit$.  For if this were not the case for some $t_0 \in I$, one can use some harmonic analysis to show that the free evolution $e^{i(t-t_0)\Delta} u(t_0)$ of $u$ from $t_0$ is dispersed for times $t$ near $t_0$, in the sense that
$$ \| e^{i(t-t_0)\Delta} u(t_0) \|_{L^{10}_{t,x}( [t_0 - C/N(t_0)^2, t_0 + C/N(t_0)^2] \times \R^3 ) } \leq 1/C$$
for $C>0$ which can be arbitrarily large (this $C$ is essentially the reciprocal of the $c$ appearing above).  On the other hand, if the free evolution is \emph{globally} small in $L^{10}_{t,x}$ norm, then perturbative theory (e.g. Theorem \ref{long-time theorem}) easily lets one show that $u$ is globally bounded in $L^{10}_{t,x}$ norm, contradicting the blowup hypothesis.  Thus $e^{i(t-t_0)\Delta} u(t_0)$ must concentrate at some time $t_1$ far away from $t_0$, say at a past time $t_1 < t_0$.  Thus the backward-propagated wave $e^{i(t_1-t_0)\Delta} u(t_0)$ has a large inner product with some highly concentrated ``wavelet'' $f$;
by duality, this means that $u(t_0)$ has a large inner product with a forward-propagated wavelet $e^{i(t_0-t_1)\Delta} f$.  We can
then split\footnote{This splitting argument is based on an earlier argument in \cite{bourg.critical}.} $u(t_0)$ into a small multiple $v(t_0)$ of this propagated wavelet, plus a remainder $w(t_0)$ of strictly smaller energy.  We use the induction on energy hypothesis to propagate $w$ to all of $\R \times \R^3$ by the nonlinear evolution, and $v$ by the linear evolution.
The point is that because $v(t_0)$ was already a wavelet propagated forward by a long period of time, the further linear propagation of $v(t_0)$ will be extremely small to the \emph{future} of $t_0$.  This allows one to apply the perturbative
theory (Theorem \ref{long-time theorem}) on the future interval $[t_0,+\infty)$, and pass from $L^{10}_{t,x}$ control of the solution $w$ to $L^{10}_{t,x}$ control of the solution $u$.  But we are assuming that $u$ blows up both to the future and to the past\footnote{In the ``finitary'' version of this argument, where $u: I \times \R^3 \to \C$ has very large but finite $L^{10}_{t,x}$ norm, what we have to do instead is split $I = I_- \cup I_0 \cup I_+$, where $I_-, I_0, I_+$ are intervals which each capture one third (or more precisely $1/3^{1/10}$) of the $L^{10}_{t,x}$ norm.  The physical space concentration effect
then only is valid on the middle third interval $I_0$; dispersion can occur at one or both of the endpoints $I_-$, $I_+$ (think of a near-soliton which stays coherent for a long time interval $I_0$ but disperses both to the future and past of this interval).  The point is that while dispersion can occur, any energy which has radiated away by dispersion cannot be subsequently reconcentrated, and so one no longer has true critical energy behaviour.}, which is a contradiction.

Once we have physical space concentration, the stronger property of localisation is obtained by a variant of the arguments used to prove Proposition \ref{freq-irreducible}.  Indeed, if localisation failed, so that a significant portion of energy at some time $t_0$ was distributed far away from $x(t_0)$, then by pigeonholing as before we can locate a splitting $u(t_0) = v(t_0) + w(t_0) + \hbox{small}$, where $v(t_0)$ is supported near $x(t_0)$, $w(t_0)$ is supported well away from the support of $v(t_0)$, and the error is very small in energy norm.  Also one can arrange matters so that $v(t_0)$ and $w(t_0)$ have energy strictly smaller than $E_\crit$.  Thus by the induction hypothesis one can propagate $v$ and $w$ by the NLS flow and obtain $L^{10}_{t,x}$ bounds.  To finish the argument one needs to show that the nonlinear interactions $O(|v| |w|^4 + |v|^4 |w|)$ between $v$ and $w$ are suitably small.  For times $t$ near $t_0$ this can be accomplished by exploiting approximate finite speed of propagation phenomena for linear and nonlinear Schr\"odinger flows, which will keep $v$ and $w$ more or less separated in physical space.  For times $t$ far away from $t_0$, the physical supports of $v$ and $w$ can intermingle; however, the physical space localisation of $v$ at time $t_0$, combined with dispersive estimates (such as those arising from pseudoconformal energy identities) will ensure that $v$ will be so small away from these times that the interactions at these times will necessarily be quite weak.  
\end{proof}

\begin{remark} An alternate approach to establishing space and frequency concentration (but not localisation) for arbitrary
solutions with large $L^{10}_{t,x}$ norm appeared in \cite{tao-nls}.  The main point there is that in order for the $L^{10}_{t,x}$ norm to be large, the nonlinear component of the Duhamel formula \eqref{ucos} must dominate.  One then inspects this component using harmonic analysis to deduce concentration, which turns out to be sufficient (in the radial case) to establish global $L^{10}_{t,x}$ bounds.  A somewhat related approach also appears in \cite{grillakis:scatter}.
At present, however, the only known proof of global existence in the non-radial case for this equation requires the full strength of the induction-on-energy machinery (or the closely related concentration compactness machinery of the next section).
Also the reliance on fundamental solution methods (i.e. the Duhamel formula) requires a substantial amount of decay on the fundamental solution, which is typically available only in high dimensions (such as three and higher), whereas the induction on energy approach extends to general dimension.
\end{remark}

Informally, what we have shown is that for a minimal-energy blowup solution $u: I \times \R^3 \to \C$, the solution concentrates at each time $t_0$ essentially all of its energy in a frequency annulus $\{ \xi: |\xi| \sim N(t_0) \}$ and in a dual spatial ball $\{ x: |x-x(t_0)| \lesssim 1/N(t_0) \}$.  A particularly elegant of saying this is that after quotienting out by the scaling and spatial translation symmetries of the NLS equation, the orbit $\{ u(t_0): t_0 \in I \}$ of the minimal-energy blowup solution is \emph{precompact} (its closure is compact).  In the language of dynamical systems, minimal-energy blowup solutions are \emph{almost periodic} modulo the symmetries of the equation.  This phenomenon is in fact very general and can be extended to other model equations in which all the ``defects of compactness'' are caused by symmetries; see \cite{tvz-compact} and the next section.  In the case of spherical symmetry (which eliminates the defect of compactness caused by translation invariance) one can basically set $x(t_0) = 0$; see \cite{bourg.critical}, \cite{tao-nls}, or \cite{tvz-compact}.

Aside from this compact dynamics, the only remaining non-compact degrees of freedom are the frequency $N(t_0)$ and the position $x(t_0)$.  The above perturbative arguments do not provide any significant long-term control on these quantities\footnote{One can however use perturbative theory to show that on time intervals centered at $t_0$ of length $\ll 1/N(t_0)^2$, the frequency $N(t_0)$ does not move by more than a constant multiplicative factor, while the position $x(t_0)$ moves by a displacement of at most $O( 1 / N(t_0) )$.  One can use this to view the solution as being composed of a sequence of ``bubbles'' of energy concentration in spacetime, where each bubble has a spatial width of $1/N$ and lifespan of $1/N^2$ for some $N$.  See \cite{taobook} for further discussion.}.  On the other hand, one can recast spacetime integrals in terms of these degrees of freedom, and thus use tools such as monotonicity formulae to obtain further control.  For instance, the fact that the $L^{10}_{t,x}$ norm of $u$ blows up both forward and backward in time can be shown to be equivalent to the assertion that the improper integral $\int_I N(t)^2\ dt$ also blows up forward and backward in time\footnote{One can view the renormalised time variable $s$ defined infinitesimally by $ds := N(t)^2\ dt$ (as well as the renormalised spatial parameter $y := N(t) (x-x(t))$) as natural scale-invariant spacetime coordinates in which to view the dynamics; see \cite{sulem} for some elaboration of this viewpoint.  This has some advantages for numerical computations, but is difficult to use analytically for a number of reasons, notably the lack of control on derivatives of $N(t)$ and $x(t)$.}  In the radial case (so $x(t) \equiv 0$), the Morawetz estimate \eqref{moraw} can be shown to be equivalent for minimal-energy blowup solutions to the Morrey-Campanato type estimate
\begin{equation}\label{jn}
\frac{1}{|J|^{1/2}} \int_J N(t)\ dt \lesssim 1
\end{equation}
for all $J \subset I$.  This comes close to contradicting the blowup of $\int_I N(t)^2\ dt$, except that the power of $N(t)$ is wrong (this is related, via scale invariance, to the undesirable weight of $\frac{1}{|J|^{1/2}}$ on the left-hand side of \eqref{moraw}).  Nevertheless, the Morawetz estimate does show that the frequency $N(t)$ cannot stay bounded by any given frequency cutoff $N_0$ for times much longer than $1/N_0^2$.  By iterating this fact in an elementary manner (see \cite{bourg.critical}, \cite{tao-nls}) one can show
that $N(t)$ must move from low frequencies to high frequencies in a rapid amount of time; indeed one can show that for any $K > 1$ there exist times $t_0, t_1$ with 
\begin{equation}\label{nkn}
N(t_1) \geq K N(t_0) \hbox{ and } t_1 = t_0 + O( N(t_0)^{-2} ).
\end{equation}
It is important to note here that the implied constant in the $O()$ notation is independent of $K$; this is ultimately due to the convergence of the geometric series $\sum_j N_j^{-2}$ when the $N_j$ are growing exponentially in $j$.

In order to exclude this last remaining blowup scenario (which can be viewed as a kind of ``self-similar'' blowup scenario)
one can exploit local approximate conservation of mass in physical space.  At time $t_0$, the frequency $N(t_0)$ is 
relatively low compared to $N(t_1)$, which (because the energy is fixed) will imply that the mass is relatively large; indeed, the mass in the ball $\{ x = O( 1/N(t_0) ) \}$ at time $t_0$ is at least as large as $c/N(t_0)^2$ for some $c > 0$.  One can then use localised mass conservation laws such as \eqref{var0} (with $a$ being a smooth cutoff to a dilated version of this spatial ball) to show that the mass in the ball $\{ x = O( 1/N(t_0) ) \}$ at time $t_1$ is also at least as large as $c'/N(t_0)^2$.  Some Fourier analysis
then shows that at time $t_1$, a significant portion of the energy must be concentrated near the frequency $N(t_0)$.  But this contradicts Corollary \ref{freqlocal} since $N(t_1) \geq K N(t_0)$ and $K$ can be taken arbitrarily large.  This concludes
the proof of Theorem \ref{l10-bound} in the spherically symmetric case.

An alternate approach, given recently by Kenig and Merle \cite{merlekenig}, uses the viriel identity as a substitute for the (localised) Morawetz inequality \eqref{moraw}.  Indeed, modifying \eqref{virial} with a suitable spatial cutoff we easily verify that
$$ \partial_{tt} \int_{\R^3} |x|^2 |u(t,x)|^2 \varphi(x/R)\ dx \geq c \int_{|x| \leq R} |\nabla u(t,x)|^2 + |u(t,x)|^6\ dx
+ O( \int_{|x| \sim R} |\nabla u(t,x)|^2 + |u(t,x)|^6\ dx )$$
for some $c > 0$, where $\varphi(x/R)$ is a cutoff supported on the ball $|x| \lesssim R$ which equals one when $|x| \leq R$.  
Integrating this on a time interval $J \subset I$ and specialising to minimal energy blowup solutions, one obtains the inequality
$$ |\{ t \in J: N(t) \gg R^{-1} \}| \lesssim R^2 + |\{ t \in J: N(t) \lesssim R^{-1} \}|.$$
This is a weaker version of \eqref{jn}, but has the same key effect, namely it prevents the frequency $N(t)$ from staying near a constant value $R^{-1}$ for periods of time much longer than $R^2$.  In conjunction with the mass conservation argument one can then obtain a bound on $\int_I N(t)^2\ dt$ as before.  The advantage of using the virial identity is that it also works well in the focusing case, even for solutions close in the energy to the stationary state, due to the variational properties of that state; see \cite{merlekenig}.

Now we turn to the non-radial case (so $x(t) \not \equiv 0$), which is significantly more difficult.  The local mass conservation argument extends to this case without difficulty, and establishes the weak continuity bound
\begin{equation}\label{nbn}
N(t_1) \leq C(B) N(t_0) \hbox{ whenever } B > 1 \hbox{ and } |t_1 - t_0| \leq B N(t_0)^{-2}
\end{equation}
where $C(B)$ is some finite quantity depending on $B$.  However, this by itself is certainly not enough to establish a bound on $\int_I N(t)^2\ dt$ (think of the ``pseudosoliton'' case when $N$ is bounded).  The Morawetz estimate \eqref{moraw} is now much weaker; it essentially asserts that
$$ \frac{1}{|J|^{1/2}} \int_J \frac{1}{N(t)^{-1} + |x(t)|}\ dt \lesssim 1 \hbox{ for all } J \subset I.$$
Since $x(t)$ can be arbitrarily far away from the origin, this estimate does not give much control on either $N(t)$ or $x(t)$, other than to say that $x(t)$ cannot linger close to the time axis for times much longer than $N(t)^{-2}$.  One can use translation invariance to generalise this bound slightly to
$$ \frac{1}{|J|^{1/2}} \int_J \frac{1}{N(t)^{-1} + |x(t)-x_0|}\ dt \lesssim 1 \hbox{ for all } J \subset I \hbox{ and } x_0 \in \R^3$$
but this is still quite weak (for instance, it cannot even prevent a ``moving pseudosoliton'' example in which $N(t)$ stays constant and $x(t)$ moves linearly in $t$).  As of this time of writing, the only monotonicity formula which is known to give a usable spacetime bound on $N(t)$ in the non-radial case is (a localised version of) the interaction Morawetz inequality \eqref{iu}.  Unlike the situation with \eqref{moraw}, it turns out that one needs to localise this inequality in \emph{frequency space} rather than in \emph{physical space}.  Indeed one has

\begin{proposition}[Frequency-localised interaction Morawetz estimate]\label{flim}\cite{gopher}  Let $u: I \times \R^3 \to \C$ be a minimal-energy blowup solution, let $\eta > 0$, and suppose that $J \subset I$ is an interval.  Let $N_*$ be such that $N_* \leq c(\eta) N(t)$ for all 
$t \in J$ and some sufficiently small $c(\eta) > 0$.  Then
\begin{equation}\label{jjo}
 \int_J \int_{\R^3} |P_{\geq N_*} u(t,x)|^4\ dx \leq \eta N_*^{-3}.
\end{equation}
\end{proposition}

The proof of this proposition is quite complicated, taking up almost 24 pages in \cite{gopher}!  The idea is to repeat the derivation of \eqref{iu} but with $u$ replaced by the high-frequency component $P_{\geq N_*} u$.  Note that the analogue of the right-hand side of \eqref{iu} can be easily
estimated as $O(\eta N_*^{-3})$.  However, there are now several new ``low-high interaction'' error terms arising from the fact that the high-frequency component does not quite solve NLS by itself.  To control these interaction terms one needs to use some perturbative analysis (and a bootstrap assumption of $L^4_{t,x}$ control on the high frequencies) to establish some preliminary estimates of Strichartz type on the low and high frequency components of $u$.  Here one crucially needs the hypothesis $N_* \leq c(\eta) N(t)$ to ensure that the low frequencies have very small energy and are thus amenable to a treatement by perturbative theory.  This deals with most of the error terms, but even so there are a few very unpleasant ``top order'' error terms which do not fall to the above estimates.  For this one needs to fully exploit the concentration properties of the minimal-energy blowup solution $u$, especially the spatial energy decay away from $x(t)$, and also to play the forward and backward Duhamel formula against each other.  See \cite{gopher} for full details.

The estimate \eqref{jjo} implies an integral bound on $N(t)$, namely
$$ \int_J N(t)^{-1}\ dt \leq C [\inf_{t \in J} N(t)]^{-3}$$
for all $J \subset I$ and some absolute constant $C$ (depending only on $E_\crit$).  This is somewhat similar to 
\eqref{jn} in that it prevents $N(t)$ from lingering near a constant value for extended periods of time.  Unfortunately
this estimate is in some sense ``too far away'' from control of $\int_I N(t)^2\ dt$ to force a rapid frequency cascade as
in \eqref{nbn}.  Instead, all one can conclude at this point is that if $\int_I N(t)^2\ dt$ is finite, then
$\sup_{j \in I} N(t) / \inf_{j \in I} N(t)$ is unbounded.  In particular, given any $K \geq 1$ we can find times $t_0, t_1$ for which
$$ N(t_1) \geq K N(t_0)$$
but for which we have no upper bound on the time difference $|t_1-t_0|$, thus prohibiting us from exploiting short-time
estimates such as \eqref{nbn} (other than to establish lower bounds on $|t_1-t_0|$).  In order to prevent this from happening, we once again must use some sort of localised mass conservation law.  The spatial localisation used previously is no longer effective at long times, but it turns out that \emph{frequency localisation} of the mass conservation law is much more effective (note that for the linear evolution, frequency localisation of data persists for arbitrarily long times, in contrast to spatial localisation).

We briefly sketch some details of the frequency localisation argument (which, while simpler than the derivation of Proposition \ref{flim}), is still non-trivial, occupying about $10$ pages of \cite{gopher}).  With a little additional argument (rescaling and exploiting the compactness of the dynamics modulo symmetries) one can pass to a minimal-energy blowup solution with a slightly stronger property, namely that there is a time $t_0$ for which $N(t_1) \geq N(t_0) = 1$ for all $t_1 \in I$ with $t_1 \geq t_0$ and 
\begin{equation}\label{soup}
\sup_{t_1 \in I; t_1 > t_0} N(t_1) = +\infty.
\end{equation}
This reduction is not absolutely essential for the argument but it does simplify things slightly.  It implies that for some sequence of times approaching the future endpoint $\sup(I)$ of the maximal lifespan $I$, the energy of the solution goes to infinity in frequency space; in particular, the solution converges weakly to zero.  This allows one to obtain a backward Duhamel formula
$$ u(t) = i \int_t^{\sup(I)} e^{i(t-t')\Delta}( |u(t')|^4 u(t'))\ dt'$$
where the improper integral has to be interpreted in a weak conditional sense, using the above-mentioned sequence of times converging to $\sup(I)$.  On the other hand, from \eqref{jjo} we also have $L^4_{t,x}$ estimates on the high frequencies of $u$ to the future of $t_0$; combining the two using Strichartz estimates, one can obtain quite strong estimates on the \emph{low} frequencies of $u$ to the future of $t_0$; in particular one has very strong energy decay as one approaches the frequency origin - much stronger (by about $3/2$ inverse derivatives) than what one obtains just from Corollary \ref{freqlocal}.  See \cite{gopher} for details.  It turns out that this control is now sufficient to establish that the high-frequency components of the solution obey an approximate mass conservation law, indeed for suitably small $\eta > 0$ one can show
$$ \int_{\R^3} |P_{\geq \eta} u(t_1,x)|^2\ dx \geq \frac{1}{2} \int_{\R^3} |P_{\geq \eta} u(t_0,x)|^2\ dx$$
for all $t_1 \geq t_0$.  In terms of the frequency variable $N()$, this implies that $N(t_1) = O(N(t_0))$ for all $t_1 \geq t_0$, contradicting \eqref{soup}.  This eliminates the last outstanding blowup scenario (a kind of ``slow low-to-high frequency cascade'') and establishes Theorem \ref{l10-bound}.

\begin{remark} The above arguments even give an explicit bound on $A(E)$ in the non-radial case, although due to the extremely heavy reliance of the induction on energy hypothesis, the bound is incredibly poor (an eightfold-iteratred exponential tower!).
In the radial case, there are methods avoiding induction on energy (or compactness) which give a more civilised exponential bound \cite{tao-nls}.  In the case of the critical NLW, the situation is better; one has exponential bounds in the non-radial case \cite{nak2}, \cite{tao} and polynomial bounds in the radial case \cite{gsv}.  We do not know at present whether any of these bounds are sharp (although the analysis from \cite{cct} in principle gives some very weak lower 
bounds).  Improving these bounds has application to pushing the critical theory to slightly supercritical regimes; see \cite{tao:supercrit}.
\end{remark}

\begin{remark} The above general scheme has been extended to higher dimensions \cite{rv}, \cite{visan}, to the nonlinear wave and Klein-Gordon equations \cite{nak-2d}, \cite{nak2}, and recently to the mass-critical NLS in high dimensions with spherical symmetry \cite{tvz}.  It is likely that the method extends further, in particular it should have relevance to the large data theory of energy-critical wave maps and mass-critical gKdV (and more ambitiously to the energy-critical MKG and YM equations, once the perturbative theory of these equations is settled).
\end{remark}

\section{Concentration compactness}\label{ccsec}

In the previous section we described a general ``induction on energy'' strategy to deal with large data solutions to a critical energy, which focused attention on the critical threshold energy between linear and nonlinear behaviour.  The arguments here tended to be quite ``quantitative'' or ``hard'' in nature, in that one relied quite heavily on various estimates arising from
either perturbative analysis (e.g. from harmonic analysis estimates on the linear propagator) or on the bounds arising from conserved and monotone quantities.  

In parallel to this, a seemingly rather different ``qualitative'' or ``soft'' strategy to control solutions, based on compactness methods (notably \emph{concentration compactness}), was developed, originally from calculus of variations (see e.g. \cite{lions0}, \cite{lions}) but in recent years now firmly established in nonlinear wave and dispersive equations.  Like the induction on energy method (when viewed contrapositively as an analysis of minimal-energy blowup solutions), the compactness method\footnote{The methods here should be compared with the compactness methods discussed in Section \ref{alternative-sec}.  In both cases one uses sequential compactness to extract solutions with special properties.  In Section \ref{alternative-sec}, the special property is an initial condition $u(t_0)=u_0$; here, the special property might be that a certain spacetime norm is infinite, that a certain energy is minimal, that there is no radiation at infinity, etc.} is somewhat indirect; in order to prove that solutions exhibit some sort of behavior, assume for contradiction that the behavior is violated, and then consider an ``extreme'' example of this violation and deduce a contradiction.  In the induction-on-energy approach, the extreme solution is obtained by minimising an energy (subject to a blowup condition, which is a kind of boundary condition).  In the compactness method, one takes an arbitrary sequence of progressively egregious examples of bad behaviour, and extracts a convergent subsequence in order to find an extreme example which has ``infinitely bad'' behaviour in some sense.  The power of this method lies in the fact that quantities which were merely decaying to zero for solutions in the sequence, would now be \emph{identically zero} for the limiting solution, which often simplifies the subsequent analysis both technically and conceptually.  
Further applications of this limit-of-subsequence idea can be used to erase all ``good'' behaviour (e.g. linear dispersion) from the solution (because dispersive behaviour often converges to zero in some weak sense), leaving a ``pure'' bad solution which is then often very rigid and can be controlled by a variety of methods (perturbation theory, monotonicity formulae, variational principles).  This latter idea has been particularly fruitful in analysing the stability of solitons for the NLS and gKdV equations (see e.g. \cite{merle}, \cite{liouville}, \cite{gkdv0}, \cite{gkdv1}), though recently it has begun to be extended to more general situations.  As it turns out, these methods can be used to reinterpret
the induction-on-energy method in a clean and qualitative context, albeit at the cost of foregoing any hope of explicit quantitative bounds.

In running the compactness method, one runs into the problem that the sequence of solutions for which one wishes to extract a convergent subsequence need not be sequentially compact, except in very weak topologies.  One can of course use the Banach-Alaoglu theorem (or more precisely the Arzela-Ascoli diagonalisation argument) to extract weakly convergent subsequences from any bounded sequence, but the main difficulty with weak convergence is that properties of the elements of the sequence (e.g. regularity, or largeness of certain norms) need not be preserved in the weak limit (although \emph{uniform upper} bounds will in general be preserved, thanks to the weak closure of the unit ball or by Fatou's lemma).  To resolve this, it becomes
necessary to seek ways in which to upgrade weak convergence to stronger notions of convergence.  

Of course, the basic problem here is that the function spaces one works in (e.g. the energy space $\dot H^1(\R^d)$) have infinitely many degrees of freedom, and thus bounded sequences in such spaces are almost certainly not compact in the strong topology.  In subcritical cases one can sometimes exploit compact embeddings (e.g. the Rellich compactness theorem) to recover compactness in slightly coarser (but still strong) topologies, but in critical cases, the presence of non-compact symmetry groups such as scaling and spatial translation show that one cannot hope for compactness in any norm which is preserved by these
symmetries, unless one somehow ``quotients out'' these symmetries first.  When one is close to a ground state, one can often exploit a variational characterisation of that ground state to obtain the desired compactness modulo symmetries, if the variational functional obeys a suitable Palais-Smale type condition.  

For more general classes of data, not close to a ground state, the presence of symmetries combined with the ability to superimpose two disjoint solutions means that the failure of strong compactness cannot be resolved merely by quotienting out by
the symmetry group.  To give a simple example, let $x_n, y_n \in \R^d$ be a sequence of points which diverge from each other in the sense that $\lim_{n \to \infty} |x_n - y_n| = \infty$, and consider the ``two bump'' examples 
$u_n(x) := \psi(x - x_n) + \psi(x-y_n)$ where $\psi$ is a test function.  Then this sequence $u_n$ is bounded in any reasonable translation-invariant norm (e.g. in the Sobolev norms $H^s(\R^d)$ for any $s$) but have no convergent subsequence in any of thse norms, even if one is allowed to translate each $u_n$ by an arbitrary amount; the problem is that one can make one of the bumps stay confined to a compact region of space (and thus have a convergent subsequence), but only at the cost of the other bump escaping to infinity, thus converging weakly to zero but diverging in every strong topology.  One can concoct similar examples with the translation symmetry replaced by other non-compact symmetries, such as scaling symmetry and modulation symmetry, provided of course that all topologies one is studying are invariant with respect to these symmetries.

Fortunately, in many situations this type of example - superpositions of fixed objects - each moved around by a different symmetry of the equation, and with the symmetries becoming ``asymptotically orthogonal'' in the limit $n \to \infty$ - turns out
to be the \emph{only} source of non-compactness for bounded sequences, provided that one is willing to measure errors in a slightly coarse topology, which allows the error to be large in energy or mass so long as it is somehow ``dispersed'' (asymptotically orthogonal to all concentrated objects).  This phenomenon, known as \emph{concentration compactness}, was introduced by Lions for applications to elliptic variational problems, although it has since proven to have many further applications. It is a surprisingly effective substitute for genuine compactness.  Informally, it says that \emph{any} bounded function splits as the ``asymptotically orthogonal'' sum of boundedly many concentrated objects (each of which can be placed into a compact region of space and frequency after applying suitable symmetries), plus a dispersed error.  In many applications the dispersed error is negligible, and the asymptotically orthogonal components become decoupled, and so the analysis reduces to understanding the compact dynamics of an evolution of concentrated fields - just as in the induction-on-energy method.

Let us now briefly outline some details of this theory.  One typically works in a Hilbert space $H$ such as $L^2_x(\R^d)$ or $\dot H^1_x(\R^d)$.  We will capture the symmetries\footnote{One can also replace this group with a more general collection of bounded operators satisfying certain axioms; see \cite{tintarev}.} by introducing a (non-compact) finite-dimensional 
Lie group $G$ of unitary transformations on $H$.  For instance, $G$ might be the group of translations $\tau_{x_0}: f(x) \mapsto f(x-x_0)$, or perhaps
the group of $L^2_x(\R^d)$-unitary dilations $\sigma_\lambda: f(x) \mapsto \lambda^{-d/2} f(\frac{x}{\lambda})$, or the group generated by both translations and dilations.  For us, the relevant properties we need are that (a) $G$ is indeed a finite-dimensional Lie group in the strong operator topology, and (b) $G$ can be compactified in the weak operator topology by adjoining $0$.  More precisely, we need the crucial \emph{dislocation property} that if $g_n$ is a sequence in $G$ which goes to infinity (i.e. it escapes every compact set, as measured in the strong operator topology), then it converges to zero in the weak operator topology.  One can easily verify that the groups discussed above have this property.

The dislocation property has the following important consequence.  Call two sequences $g_n, g'_n \in G$ \emph{asymptotically orthogonal} if $(g'_n)^{-1} g_n$ goes to infinity in $G$.  Then for every $f, f' \in H$ we have $\lim_{n \to \infty} \langle g_n f, g'_n f' \rangle_H = 0$, explaining the terminology ``asymptotically orthogonal''.  

A related consequence is as follows.  Let us say that a bounded sequence $f_n \in H$ \emph{converges weakly to zero with $G$-concentration} if the sequence $g_n f_n$ converges weakly to zero for \emph{any} sequence $g_n \in G$; this is thus intermediate in strength between weak and strong convergence.  For instance, the two-bump example mentioned earlier does not converge weakly to zero modulo the group of translations, because we can translate so that one of the bumps stays near the origin, thus ensuring failure of weak convergence to zero.  Intuitively, sequences which converge weakly to zero with $G$-concentration are ``dispersed'' even if they stay large in the strong norm $\| \|_H$, because they are asymptotically orthogonal to all concentrated functions (fixed functions, moved around by arbitrary group elements). 

\begin{lemma}[Abstract dichotomy between dispersion and concentration]  Let $G,H$ be as above. Let $f_n \in H$ be a bounded sequence which does not converge weakly with $G$-concentration to zero.  Then by passing to a subsequence if necessary, we can find a non-zero $\phi \in H$ and a decomposition $f_n = g_n \phi + f'_n$, where $g_n \in G$, and $g_n^{-1} f'_n$ converges weakly to zero.  In particular $g_n \phi$ and $f'_n$ are asymptotically orthogonal.

Furthermore, if $g'_n$ is any sequence in $G$ such that $(g'_n)^{-1} f_n$ converges weakly to zero, then $g_n$ and $g'_n$ are asymptotically orthogonal.
\end{lemma}

\begin{proof} Since $f_n$ does not converge weakly with $G$-concentration to zero, we can find $g_n$ such that $g_n^{-1} f_n$ does not weakly converge to zero.  By weak compactness, we may then pass to a subsequence for which $g_n^{-1} f_n$ converges weakly to a non-zero $\phi$.  Setting $f'_n := f_n - g_n \phi$ we obtain the first claim.

To prove the second claim, assume for contradiction that we can find $g'_n$ such that $(g'_n)^{-1} f_n$ converges weakly to zero, but that $g_n$ and $g'_n$ are not asymptotically orthogonal.  By the dislocation property, we may thus pass to a subsequence where $g_n^{-1} g'_n$ converges strongly to some fixed group element $g$, and thus $g_n^{-1} f_n = (g_n^{-1} g'_n) (g'_n)^{-1} f_n$ converges weakly to zero.  But this contradicts the fact that $g_n^{-1} f_n$ converges to the non-zero $\phi$.
\end{proof}

Repeated iteration of this lemma eventually leads to

\begin{corollary}[Abstract concentration compactness]\cite{tintarev} Let $G,H$ be as above. Let $f_n \in H$ be a bounded sequence.  Then after passing to a subsequence we have an absolutely convergent decomposition
$$ f_n = \sum_{j=1}^\infty g_n^{(j)} \phi^{(j)} + w_n,$$
where $\phi^{(j)} \in H$ are functions, $g_n^{(j)}$ are sequences of group elements with $g_n^{(j)}$ and $g_n^{(j')}$ asymptotically orthogonal for all $j \neq j'$, and $w_n$ is bounded in $H$ and converges weakly with $G$-concentration to zero.  Furthermore we have the asymptotic Pythagoras theorem
$$ \limsup_{n \to \infty} \| f_n \|_H^2 = \sum_{j=1}^\infty \| \phi^{(j)} \|_H^2 + 
\limsup_{n \to \infty} \| w_n \|_H^2.$$
\end{corollary}

\begin{remark} It turns out that for many applications in nonlinear dispersive and wave equations it is better to use a truncated version of the above decomposition, in which we only sum finitely many of the main terms $g_n^{(j)} \phi^{(j)}$, at the cost of worsening the behaviour of the error $w_n$.  We shall describe such a truncated version shortly.
\end{remark}

In order to use this type of concentration compactness result effectively, one needs to deal with the error $w_n$.  It is here that the choice of group $G$ becomes important (beyond merely obeying the dislocation property), for when $G$ is sufficiently large, one can upgrade weak convergence with $G$-concentration to strong convergence in various Banach space norms $\| \|_X$ which are controlled by $H$.  Roughly speaking, this occurs when the group $G$ captures all the ``defects of compactness'' of the embedding of $H$ into $X$; in more quantitative terms, this means that the $X$ and $H$ norms are only comparable for functions which correlate with a test function, shifted by a group element in $G$.  A prototypical example arises from non-endpoint Sobolev embedding, such as $H^1_x(\R^3) \subset L^3_x(\R^3)$.  When the domain is compact, the well-known Rellich compactness theorem shows that this embedding is compact, in particular weak convergence in bounded subsets of $H^1_x$ implies strong convergence in $L^3_x$.  For unbounded domains such as $L^3_x(\R^3)$, the invariance under the group $G$ of translations shows that the embedding can no longer be compact; nevertheless, we have

\begin{lemma}[Concentration-compact Sobolev embedding]\label{ccse} 
Any bounded sequences in $H^1_x(\R^3)$ which converge weakly \emph{with $G$-concentration} also converges strongly in $L^3_x(\R^3)$.
\end{lemma}

For a proof, see e.g. \cite{lions}.
One can use ``soft'' arguments to show that the above ``qualitative'' statement is in fact equivalent to the following ``quantitative'' assertion:

\begin{lemma}[Inverse Sobolev theorem]  Let $G$ be the group of translations.
For every $\eta > 0$ there exists a finite set ${\mathcal E}_\eta \subset C^\infty_0(\R^3)$ of test functions
with the following property: for every $f \in H^1_x(\R^3)$ such that $\|f\|_{H^1_x(\R^3)} \leq 1$ and $\| f\|_{L^3_x(\R^3)} \geq \eta$, there
exists $\phi \in {\mathcal E}_\eta$ and $g \in G$ such that $|\langle f, g\phi \rangle| \geq 1$.
\end{lemma}

This lemma can in turn be proven by a variety of means, for instance by using Littlewood-Paley theory, or the wavelet characterisation of various Besov and Sobolev function spaces.  Using this fact, one can convert the abstract concentration compactness result into something more quantitative.  For instance, we have

\begin{proposition}[Profile decomposition for $H^1_x(\R^3) \subset L^3_x(\R^3)$]\cite{gerard} Let $G$ be the translation group on $\R^3$. Let $f_n \in H^1(\R^3)$ be a bounded sequence.  Then after passing to a subsequence we have decompositions
$$ f_n = \sum_{j=1}^l g_n^{(j)} \phi^{(j)} + w_{n,l}$$
for all $l \geq 0$, where $\phi^{(j)} \in H^1_x(\R^3)$ are functions, $g_n^{(j)} \in G$ are sequences of group elements with $g_n^{(j)}$ and $g_n^{(j')}$ asymptotically orthogonal for all $j \neq j'$, and $w_{n,l}$ is bounded in $H^1_x(\R^3)$ with
$$ \lim_{n \to \infty} \limsup_{l \to \infty} \|w_{n,l} \|_{L^3_x(\R^3)} = 0.$$
Furthermore we have the asymptotic Pythagoras theorem
$$ \limsup_{n \to \infty} \| f_n \|_{H^1_x(\R^3)}^2 = \sum_{j=1}^l \| \phi^{(j)} \|_{H^1_x(\R^3)}^2 + 
\limsup_{n \to \infty} \| w_{n,l} \|_{H^1_x(\R^3)}^2$$
for all $l \geq 0$.
\end{proposition}

Note that the embedding $H^1_x(\R^3) \subset L^3_x(\R^3)$ is invariant under translations, but not under other symmetries such as scaling or frequency modulation.  This is basically why the translation group $G$ is the natural group that appears for this embedding.  For applications to \emph{critical} (scale-invariant) problems, however, we need to understand the defect of compactness for embeddings which are invariant both under scaling as well as translation.  A good example is the Strichartz embedding
$$ \| e^{it\Delta} f \|_{L^{10}_{t,x}(\R \times \R^3)} \leq C \|f\|_{\dot H^1_x(\R^3)}$$
which we have already seen to play a major role in the theory of the energy-critical NLS.  This estimate is invariant under the group $G''$ generated by translations, $\dot H^1_x(\R^3)$-preserving scalings $f(x) \mapsto \frac{1}{\lambda^{1/2}} f(\frac{x}{\lambda})$, and the linear propagators $e^{it\Delta}$.  This group $G'$ also enjoys the dislocation property, and one can show the analogue of Lemma \ref{ccse}, namely that if $f_n$ is bounded in $\dot H^1_x(\R^3)$ and converges weakly modulo $G'$, then $e^{it\Delta} f_n$ converges in $L^{10}_{t,x}(\R \times \R^3)$.  As a consequence we have a profile decomposition:

\begin{proposition}[Profile decomposition for linear Schr\"odinger waves]\label{profile}\cite{keraani} Let $G''$ be as above.
 Let $f_n \in \dot H^1(\R^3)$ be a bounded sequence.  Then after passing to a subsequence we have decompositions
$$ f_n = \sum_{j=1}^l g_n^{(j)} \phi^{(j)} + w_{n,l}$$
for all $l \geq 0$, where $\phi^{(j)} \in \dot H^1_x(\R^3)$ are functions, $g_n^{(j)} \in G''$ are sequences of group elements with $g_n^{(j)}$ and $g_n^{(j')}$ asymptotically orthogonal for all $j \neq j'$, and $w_{n,l}$ is bounded in $\dot H^1_x(\R^3)$ with
$$ \lim_{n \to \infty} \limsup_{l \to \infty} \|e^{it\Delta} w_{n,l} \|_{L^{10}_{t,x}(\R \times \R^3)} = 0.$$
Furthermore we have the asymptotic Pythagoras theorem
$$ \limsup_{n \to \infty} \| f_n \|_{\dot H^1_x(\R^3)}^2 = \sum_{j=1}^l \| \phi^{(j)} \|_{\dot H^1_x(\R^3)}^2 + 
\limsup_{n \to \infty} \| w_{n,l} \|_{\dot H^1_x(\R^3)}^2$$
for all $l \geq 0$.
\end{proposition}

Similar profile decompositions are known for other equations and regularities, for instance for the wave equation in the energy class see \cite{bg}.

These profile decompositions combine very well with stability theory such as Theorem \ref{long-time theorem}, especially when the underlying group $G$ is also a symmetry group for the equation.  Roughly speaking, they assert that the asymptotic behaviour of any sequence of solutions from initial data $f_n$ decouples into the asymptotically orthogonal superposition of the solutions arising from the data $\phi^{(j)}$, moved around by symmetries of the group, plus a negligible radiation term.
(See \cite{bg} for a precise formulation of this statement, in the context of the energy-critical NLW.)  This type of decoupling has many uses.  For instance, one can analyse the behaviour of a solution near a singularity by continually rescaling around that singularity and then applying the above profile decompositions to the sequence of rescaled solutions; see \cite{mv}
for a very typical instance of this type of argument.  More recently, in \cite{merlekenig} it was observed that this profile decomposition can be used (together with the stability theory) to very quickly imply the localisation results in Corollary \ref{freqlocal} and Proposition \ref{spatloc}.  A key lemma is

\begin{lemma}[Palais-Smale type lemma modulo $G$]\label{palais}\cite{merlekenig}  Let $\mu = +1, d = 3, p = 5$, and suppose that the critical energy $E_\crit$ for NLS is finite.   Let $G'$ be the group of unitary transformations on $\dot H^1(\R^3)$ generated by translations and dilations, and let $f_n$ be a sequence of initial data with energy less than or equal to $E_\crit$ whose
maximal Cauchy developments $u_n: I_n \times \R^3 \to \C$ blow up in $L^{10}_{t,x}$ both forward and backward in time, thus
$$ \| u_n \|_{L^{10}_{t,x}(I_n \cap [0,+\infty) \times \R^3)}, \| u_n \|_{L^{10}_{t,x}(I_n \cap (-\infty,0] \times \R^3)} \to \infty \hbox{ as } n \to \infty.$$
Then after passing to a subsequence, the $f_n$ will be strongly convergent in $\dot H^1_x(\R^d)$ modulo $G'$, thus there exist $g_n \in G'$ such that $g_n^{-1} f_n$ converges strongly in $\dot H^1_x(\R^d)$.
\end{lemma}

\begin{proof}[Sketch] We use an argument from \cite{tvz-compact}.  We apply the profile decomposition from Proposition \ref{profile}, passing to a subsequence if necessary, thus writing $f_n$ in terms of components $\phi^{(j)}$, moved around by group elements $g_n^{(j)} \in G''$ plus negligible errors $w_{n,l}$.

A technical difficulty arises because of the presence of the linear propagators $e^{it\Delta}$ in the group elements $g_n^{(j)}$, because these propagators are not symmetries of NLS.  For now let us simply ignore the linear propagators and assume that $g_n^{(j)}$ consists entirely of translations and dilations, i.e. that $g_n^{(j)}$ lies in $G'$; we briefly comment on what changes have to be made to address the general case at the end of this sketch.

First suppose that all the components $\phi^{(j)}$ have energy strictly less than $E_\crit$.  Then by induction hypothesis, one can find global solutions with initial data $\phi^{(j)}$ with a bounded $L^{10}_{t,x}$ norm.  By the translation and scaling symmetries of NLS, we can achieve a similar statement for $g^{(j)}_n \phi^{(j)}$.  The asymptotic orthogonality of the $g^{(j)}_n$ (and the dispersed nature of the errors $w_{n,l}$) then allows us to superimpose these solutions together and
obtain an $L^{10}_{t,x}$ bound for the $u_n$ for sufficiently large $n$, a contradiction.  

Thus at least one of the components $\phi^{(j)}$ must have energy at least $E_\crit$.  An asymptotic Pythagoras-type theorem for the energy then shows that that component will have energy \emph{exactly} $E_\crit$, while all other components will vanish, and the error $w_{n,l}$ will have asymptotically vanishing energy as $n \to \infty$.  This implies that $f_n$ converges strongly in $\dot H^1_x(\R^d)$ as desired.  (Compare this with the heuristic from the previous section that minimal energy blowup solutions must be ``irreducible''.)

Now we comment on what happens when the $g_n^{(j)}$ contain some linear time propagation, thus
$g_n^{(j)} = \tilde g_n^{(j)} e^{it_n^{(j)} \Delta}$ for some $\tilde g_n^{(j)} \in G'$ and $t_n^{(j)} \in \R$.  For sake of argument let us just work with a single $j$.  If 
the $t_n^{(j)}$ stay bounded then after passing to a subsequence we can make them converge to a finite time as $n \to \infty$,
at which point it is easy to absorb these propagators into the $\phi^{(j)}$ and $w_{j,l}$ and argue as before.  If instead
the $t_n^{(j)}$ go to $-\infty$ (say) then the nonlinear evolution of $e^{it_n^{(j)} \Delta} \phi^{(j)}$ can be approximated by the nonlinear evolution of $\phi^{(j)}_+$, shifted in time by $t_n^{(j)}$, where $\phi^{(j)}_+$ is the forward scattering state of $\phi^{(j)}$ as Theorem \ref{nls-crit-lwp}.  Applying the symmetry associated to $\tilde g_n^{(j)}$ one can then control
the nonlinear evolution of $g^{(j)}_n \phi^{(j)}$ as before.  Continuing the argument, we eventually see that $f_n$ is asymptotically close to $g^{(j)}_n \phi^{(j)}$ in the $\dot H^1_x(\R^3)$ norm.  But from this and the stability theory one can 
easily show that $u_n$ converges to zero forward in time in the $L^{10}_{t,x}$ norm (because the same is true for the linear evolution of $g^{(j)}_n \phi^{(j)}$), a contradiction.  Hence this case cannot occur. A similar argument also works if $t_n^{(j)}$ goes to $+\infty$.  These three cases cover all the possibilities (after passing to a subsequence), and we are done.
\end{proof}

Just as the classical Palais-Smale condition in calculus of variations implies the existence of minimisers,
Lemma \ref{palais} implies the following result, which in turn can be easily shown by simple compactness arguments to imply
Corollary \ref{freqlocal} and Proposition \ref{spatloc}:

\begin{corollary}[Existence of almost periodic minimal energy blowup solutions]\cite{merlekenig}  Let $\mu = +1, d = 3, p = 5$, and suppose that the critical energy $E_\crit$ for NLS is finite.   Let $G'$ be the group of unitary transformations on $\dot H^1(\R^3)$ generated by translations and dilations.
Then there exists a minimal energy blowup solution $u: I \times \R^3 \to \C$ which blows up both forward and backward in time,
and whose orbit $\{ u(t): t \in I \}$ is precompact modulo $G'$ in $\dot H^1_x(\R^3)$, or in other words there exists a compact set $K \subset \dot H^1_x(\R^3)$ and a map $g: I \to G'$ such that $g(t)^{-1} u(t) \in K$ for all $t \in I$.
\end{corollary}

\begin{proof}(Sketch) We again use an argument from \cite{tvz}.  By definition of $E_\crit$ we can find a sequence of initial data $f_n$ of energy at most $E_\crit$ whose maximal Cauchy developments $u_n$ asymptotically blow up in $L^{10}_{t,x}$ norm.  By translating in time appropriately one can easily ensure that these $u_n$ in fact asymptotically blow up both forward and backward in time.  We apply Lemma \ref{palais} to pass to a limit $f$, and from the stability or well-posedness theory it is not hard to see that the maximal Cauchy development $u$ to this data must blow up forward and backward in time.  In particular $u$ must be a minimal energy blowup solution.

Now suppose for contradiction that the orbit of $u$ is not precompact modulo $G'$, then there exists a sequence of times $t_n$ where $g_n^{-1} u(t_n)$ has no convergent subsequence for any $g_n \in G'$.  But then we can apply Lemma \ref{palais} to the initial data $f_n := u(t_n)$ and obtain the desired contradiction.
\end{proof}

Analogues of this result exist for focusing NLS \cite{merlekenig} and for $L^2$-critical NLS \cite{tvz}.  It is likely that this type of result in fact very general and should apply to any equation with a symmetry group which is large enough to cover all the essential defects of compactness in the perturbation theory.

In view of this Corollary, one can reduce Theorem \ref{l10-bound} to the following rigidity result, which is known as a ``Liouville theorem'' in analogy to the classical result of Liouville that any entire function which is bounded must in fact be constant.

\begin{theorem}[Liouville theorem]  Let $\mu=+1, d=3, p=5$, and let $u: I \times \R^3 \to \C$ be a maximal Cauchy development for NLS whose orbit is precompact modulo $G'$.  Then $u$ is identically zero.
\end{theorem}

This theorem can be proven using the localised Morawetz and mass conservation laws of the previous section; in the spherically symmetric case it can be achieved using localised virial identities and mass conservation, see \cite{merlekenig}.  The latter argument has the significant advantage that it also extends to the focusing case, so long as the energy and $\dot H^1_x(\R^3)$ norm of the initial data are strictly less than that of the ground state.  This two-step approach of controlling arbitrary solutions by first using compactness methods to reduce to ``almost periodic'' solutions, and then using additional arguments (typically based on various localisations of conservation laws and monotonicity formulae) to establish Liouville
theorems for such solutions, also underlies a number of other recent breakthroughs in this field, for instance in the stability theory of solitons
for critical gKdV \cite{merle}, \cite{liouville}, \cite{gkdv0}, \cite{gkdv1} and also for the critical theory of NLS at exponents other than the mass or energy \cite{rap0}, \cite{rap}. 

\section{Gauge fixing}\label{gauge-sec}

In the preceding sections we have discussed the small and large data wellposedness theory for various semilinear wave equations (particularly NLS and NLW), in which the nonlinearity did not involve derivatives.  Because of this low-order nature of the nonlinearity, it was relatively easy to
apply perturbation theory to approximate the nonlinear flow by the linear one (assuming that certain key norms are small or at least finite, 
of course).  This then set the stage for further tools, such as conservation laws, monotonicity formulae, and concentration compactness to be applied.

However, once one turns to equations with derivatives in the nonlinearity, such as the WM, MKG, YM equations\footnote{The gKdV equation also has derivatives which cause some analytical difficulty, but it turns out in this case that the high order of dispersion in the linear term $u_{xxx}$ generates enough of a local smoothing effect to compensate for this loss of one degree of regularity in the nonlinearity, and so the gKdV perturbation theory is closer in spirit to the NLS and NLW than to the WM, MKG, and YM equations.  See \cite{kpv:kdv}, \cite{tao-gkdv4}.}, the presence of a derivative in the nonlinearity becomes highly troublesome for the perturbation theory, especially when one seeks a scale-invariant theory (which is needed in order to obtain global-in-time asymptotic control).  In particular, the \emph{sign} of the nonlinearity, which previously played absolutely no role in the perturbative theory, is now often decisive.
We illustrate this with an example of Nirenberg.  Let us first consider solutions $\phi: \R^{1+2} \to \R$ to the wave maps-type equation
$$ -\partial_{tt} \phi + \Delta \phi = |\phi_t|^2 - |\nabla \phi|^2.$$
Formally, one has solution to this equation of the form $\phi = \log u$, where $u: \R^{1+2} \to \R^+$ solves the \emph{linear} wave equation
\begin{equation}\label{utt}
-\partial_{tt} u + \Delta u = 0.
\end{equation}
Of course, the logarithm function has a singularity at zero.  This is not a problem locally in time if the solution is sufficiently regular, since $u = e^\phi$ will stay away from zero at the initial time $t=0$, and hence for a short time after that if $u$ is smooth enough.  However, if the initial position and velocity of $\phi$ and $u$ lie in the energy class $\dot H^1_x(\R^2) \times L^2_x(\R^2)$, which just barely fails to imply
continuity (or even boundedness) on $u$ or $\phi$ due to the logarithmic failure of Sobolev embedding, then it is not difficult to construct
examples of solutions $\phi$ which have bounded or even small energy at time zero, but develop singularities instantaneously afterwards.
In particular the standard perturbative approach to analysing this equation in the energy class must \emph{necessarily} fail no matter how cleverly one chooses the spaces to iterate in.  This can also be seen by analysing the Taylor expansion
$$ \log (1+u) = u - \frac{u^2}{2} + \frac{u^3}{3} - \ldots$$
for $u$ in $\dot H^1_x(\R^2)$.  The first term of this expansion is of course also in the energy class $\dot H^1_x(\R^2)$, but subsequent terms will not, because the space $\dot H^1_x(\R^2)$ is not closed under multiplication (this is again related to the failure of the endpoint Sobolev theorem to embed $\dot H^1_x(\R^2)$ into $L^\infty_x(\R^2)$).

On the other hand, consider the very similar equation
$$ -\partial_{tt} \phi + \Delta \phi = \phi \times ( |\phi_t|^2 - |\nabla \phi|^2 )$$
where $\phi: \R^{1+2} \to S^1$ now takes values on the unit circle $S^1 := \{ z \in \C: |z| = 1 \}$.  The presence of the additional bounded factor $\phi$ should not significantly affect the perturbation theory.  On the other hand, this equation can be solved explicitly by the substitution
$\phi = e^{iu}$ for \emph{real-valued} $u: \R^{1+2} \to \R$, and one quickly sees that $u$ (formally at least) must solve the linear wave equation \eqref{utt}.  Now the nonlinear map $u \mapsto e^{iu}$ is well-behaved on the energy class $\dot H^1_x(\R^2)$ for real-valued $u$, indeed it clearly preserves the $\dot H^1_x(\R^2)$ norm, and with a little additional effort one can even show this map is continuous in $\dot H^1_x(\R^2)$.  This is despite the failure of the power series
$$ e^{iu} = 1 + iu - \frac{u^2}{2!} - \frac{iu^3}{3!} + \ldots$$
to converge or even have its quadratic and higher terms to make sense in the energy class $\dot H^1_x(\R^2)$; the map $u \mapsto e^{iu}$ is continuous in $\dot H^1_x(\R^2)$ but not analytic.  Note that for this map to be well-behaved one has to crucially exploit the simple but nonlinear (and
non-perturbative) observation that $e^{iu}$ is bounded whenever $u$ is real; the map $u \mapsto e^{iu}$ can easily be shown to be very badly behaved in $\dot H^1_x(\R^2)$ when $u$ is no longer assumed to be real.

The above simple examples already show that a simple algebraic transformation can sometimes simplify a nonlinear equation into a linear one.
In the case of the wave maps equation, this type of transformation is available whenever the target manifold is one-dimensional, or (slightly more generally) if the initial data lies on (and moves tangentially to) a geodesic in the target; a nonlinear transformation based on the arclength parameterisation of the geodesic will then convert the wave maps equation to the free wave equation (actually this is geometrically obvious from any intrinsic formulation of the wave maps equation, such as the Lagrangian one, since geodesics are isometric to subsets of $\R$).  One can generalise this slightly to the case of wave maps from $\R \times \R^2$ into a surface of revolution which has an equivariant $U(1)$ rotation symmetry; in this case, the wave maps equation does not collapse all the way down to the free wave equation due to a residual non-flatness in the angular directions, but it does simplify to a semilinear NLW-type equation which can then be handled by existing perturbation theory techniques (e.g. Strichartz estimates)
even at the critical regularity $\dot H^1_x(\R^2) \times L^2_x(\R^2)$; see e.g. \cite{shatah-struwe}.

For general target manifolds, one cannot hope to find such a nonlinear transformation (essentially a selection of coordinates on the target) that achieves such a dramatic reduction in the strength of the nonlinearity; it is akin to hoping for a coordinate system on an arbitrary manifold which flattens most components of the metric.  Of course, the Riemann curvature tensor provides an inherent geometric obstruction to this goal.  It turns out however that if one works not on the manifold directly, but on the \emph{tangent bundle} of that manifold (basically by differentiating the wave maps equation), one obtains a much richer class of ``gauge transformations'' which can be used to weaken the nonlinearity.  

From an algebraic perspective, the advantage of differentiating the equation lies in the fact that the nonlinearity becomes linear in first derivatives instead of quadratic.  Very schematically, if one starts with an equation of the rough form
$$ \Box \phi = O( \partial \phi \partial \phi )$$
and differentiates it, setting $\psi := \partial \phi$, one expects by the product rule to get a (non-scalar, overdetermined) equation of the rough form
$$ \Box \psi = O( \psi \partial \psi )$$
The nonlinearity now is linear in first derivatives and thus has a ``magnetic'', or more generally a ``connection'' flavour.  This will be formalised geometrically later, but let us first argue algebraically.  Consider a magnetic (or ``$U(1)$ covariant'') wave equation of the form
$$ \Box \psi = - 2 i A^\alpha \partial_\alpha \psi$$
where $\psi: \R^{1+d} \to \C$ and $A^\alpha$ are some real-valued coefficients, which one should think of as being ``smooth'' and fixed\footnote{More generally, one can view $\psi$ as living in a vector space $\R^n$ and $i A^\alpha$ taking values in the skew-adjoint operators on such spaces; this is the case of interest for Yang-Mills equations, and for wave maps into targets of dimension higher than two.  However this case is slightly more complicated due to the non-abelian nature of the gauge group and we shall avoid discussing it here.}.  This equation is linear in $\psi$, but the term on the right-hand side (which is analogous to the ``nonlinearity'') involves first order derivatives in $\psi$.  In some cases however, we can transform this equation to eliminate or at least weaken this derivative term.  If we make the \emph{gauge change} $\tilde \psi := e^{i \chi} \psi$ for some arbitrarily chosen field $\chi: \R^{1+d} \to \R$, then we see (formally at least) that $\tilde \psi$ solves the wave equation
$$ \Box \tilde \psi = - 2i\tilde A^\alpha \partial_\alpha \psi - (\partial_\alpha A^\alpha) \tilde \psi$$
where $\tilde A^\alpha := A^\alpha - \partial^\alpha \chi$.  If we can arrange for the transformed connection $\tilde A^\alpha$ to vanish or be otherwise ``negligible'', then we have significantly improved the right-hand side of this equation as the remaining term no longer involves derivatives of $\psi$ or $\tilde \psi$.  (We will consider $A$ as being smoother than $\psi$, so that, all else being equal, a term with derivatives on $A$ is preferable to one with derivatives on $\psi$.)  

In general, we do not expect to be able to make $\tilde A^\alpha$ to vanish completely, as this is asking $A^\alpha$ to be a gradient\footnote{For scalar Schr\"odinger equations in one spatial dimension, the connection only has one component and is thus a gradient by the fundamental theorem of calculus.  This can be used to eliminate magnetic components completely in this special case.  A variant of this trick has proven decisive in the low-regularity theory of the Benjamin-Ono equation, in effect neutralising the effect of the derivative from the nonlinearity; a key observation is that the Benjamin-Ono equation can be recast using Riesz projections as a nonlinear Schrodinger equation with a nonlinearity which is of magnetic type (plus a small non-local error).  See \cite{tao-bo}, \cite{bp}, \cite{ik}.}. The obstruction to this occuring is described (locally, at least) by the \emph{curvature tensor}
$$ F^{\alpha \beta} := \partial^\alpha A^\beta - \partial^\beta A^\alpha;$$
observe that this curvature is unaffected by gauge transforms.  Thus it is necessary for the curvature to vanish in order
for $A^\alpha$ to be transformed to the zero connection; the contractibility of spacetime $\R^{1+d}$ ensures that the converse is also true.  If the curvature is non-zero but small in some sense, then we cannot make $\tilde A^\alpha$ vanish entirely, but we can hope to make it small also by choosing $\chi$ appropriately.  For instance, one can consider the (formal) variational problem of minimising the $L^2$ norm $\int_{\R^d} \sum_{j=1}^d |\tilde A_j(t,x)|^2\ dx$ for each $t$ (note that we are ignoring the $A_0$ component for now).  This leads to the \emph{Coulomb gauge condition}
$$ \div \tilde A = \partial_j \tilde A_j = 0.$$
In terms of the gauge field $\chi$, this becomes the elliptic equation
$$ \Delta \chi = \partial_j A_j$$
which thus has a unique solution (assuming suitable decay and regularity hypotheses on $A, \chi$.  The gauge transformed connection $\tilde A$ can also be read off directly from the curvature via the elliptic equations
$$ \Delta \tilde A_\alpha = \partial_j \partial_j \tilde A_\alpha - \partial_\alpha \partial_j A_j = \partial_j F_{j\alpha}.$$
Thus schematically we have $\tilde A = O( \nabla^{-1} F )$, so that if $F$ is small in suitable norms then $\tilde A$ is also small in a norm of one higher degree of regularity.  Given that $F$ was essentially a derivative of $A$ (and $\tilde A$) in the first place, we see that this should be about the best we can do in minimising the size of the connection.  

The Coulomb gauge was used crucially\footnote{It is however possible to see these null forms also appear in some other gauges, such as the temporal gauge; see \cite{tao:yang}.} in the sub-critical local wellposedness theory of the MKG and YM equations
in \cite{kl-mac}, \cite{klainerman:xsb}, \cite{kl-tar:yang-mills}, in order to generate certain ``null form'' structures
in the nonlinearity which provided enough cancellation for an iteration argument to establish local existence; this should be constrasted with the examples from \cite{lindblad:sharpduke} which showed that wellposedness can fail even for subcritical regularities for (non-geometric) wave equations whose nonlinearities did not obey the null condition.    Even with this gauge, however, the well-posedness (or regularity) theory at the critical regularity (and in particular, the establishment of global solutions for data with small critical norm) had been elusive until very recently.  An initial breakthrough was established by Tataru \cite{tataru:box5}, \cite{tataru:wave1}, \cite{tataru:wave2}, who introduced sophisticated refinements of existing function spaces to essentially push the iteration method to its natural limit, namely a critical-regularity Besov space (basically, this is the minimal strengthening of the critical Sobolev space required to obtain some substitute for false endpoint Sobolev embeddings such as $\dot H^{d/2}_x(\R^d) \not \subset L^\infty_x(\R^d)$).  These spaces resolved a certain technical ``division problem'' which was preventing scale-invariant iteration methods from working, leaving only the interaction between different frequency ranges as the only remaining obstacle to a critical Sobolev space theory.

For wave maps, the key to proceeding further was to recast this equation as
An equation with a gauge symmetry.  We have already sketched how this could be done by differentiating the equation.  A slightly different approach, adopted first in \cite{tao:wavemaps}, \cite{tao:wavemaps2}, performed Littlewood-Paley

projections instead of taking derivatives in order to reveal a connection-type structure.  Later, in \cite{nahmod}, \cite{shsw:wavemap}, a simpler and more geometric perspective was introduced to greatly clarify the situation.  Given any map (not necessarily a wave map) $\phi: \R^{1+d} \to M$, the tangent bundle $TM$ of $M$ pulls back to a vector bundle $\phi^*(TM)$ on $\R^{1+d}$.  The partial derivatives $\partial_\alpha \phi$ are then sections of this bundle.  The Levi-Civita connection $\nabla$ on $TM$ similarly pulls back to a connection $\phi^* \nabla$ on $\phi^*(TM)$, and the wave maps equation becomes
$$ \phi^* \nabla^\alpha \partial_\alpha \phi = 0.$$
This formulation is manifestly geometric, but difficult to analyze due to the lack of a co-ordinate system for the vector bundle $\phi^*(TM)$.  To address this, one can choose an (at present arbitrary) orthonormal frame bundle $e_1,\ldots,e_m$ on $\phi^*(TM)$, where $m$ is the dimension of $M$ (and hence of the vector bundle).  Note that the Riemannian metric on $M$ pulls back to a Hilbert space structure on each fibre of $\phi^*(TM)$, so the notion of an orthonormal frame makes sense; the contractibility of the domain $\R^{1+d}$ also makes it easy to ensure that at least one continuous orthonormal frame exists (at least for smooth $\phi$).  Using this frame, one can rewrite the derivative $\partial_\alpha \phi$ as an $\R^m$-valued field $\psi_\alpha := (\psi^1_\alpha,\ldots,\psi^m_\alpha)$ by the formula 
$$ \psi^i_\alpha := \langle \partial_\alpha \phi, e_i \rangle$$
where one uses the Hilbert space structure on $\phi^*(TM)$.  Similarly, the connection $\phi^* \nabla_\alpha$ can now be rewritten as $D_\alpha := \partial_\alpha + A_\alpha$, where $A_\alpha$ is the skew-adjoint matrix on $\R^m$ with components
$$ A_\alpha^{ij} := \langle \phi^* \nabla_\alpha e_i, e_j \rangle.$$
The wave maps equation now becomes $D^\alpha \psi_\alpha = 0$, while the torsion-free nature of the Levi-Civita equation forces the compatibility condition
$$ D_\alpha \psi_\beta – D_\beta \psi_\alpha = 0.$$
Finally, the curvature of the target $M$ manifests itself as an equation for the curvature $F_{\alpha \beta} := [D_\alpha,D_\beta]$ of the connection.  For instance, if $M$ has constant curvature $\kappa$, then standard differential geometry computations show that
$$ F_{\alpha \beta} = \kappa \psi_\alpha \wedge \psi_\beta.$$
These are now the three equations of motion for the wave maps equation when viewed using the ``differentiated fields'' $A_\alpha$ and $\psi_\alpha$.  On differentiating the wave maps equation we thus see that $\psi_\alpha$ obeys a covariant cubic nonlinear wave equation:
$$ D^\alpha D_\alpha \psi_\beta = \kappa (\psi_\alpha \wedge \psi_\beta) \neg \psi^\alpha.$$
 Because our orthonormal frame was chosen arbitrarily, one has a gauge freedom
$$ A_\alpha \mapsto U A_\alpha U^{-1} – (\partial_\alpha U) U^{-1}; \quad \psi_\alpha \mapsto U \psi_\alpha$$
for an arbitrary rotation matrix-valued gauge field $U$.  As before, one can exploit this gauge freedom to place the connection $A$ in a convenient form.
By using the Coulomb gauge $\div A = 0$, small data global regularity for wave maps at the critical Sobolev regularity was established in four and higher dimensions in \cite{nahmod}, \cite{shsw:wavemap} (with a microlocal Coulomb gauge approach giving a similar result in five and higher dimensions in \cite{kr:wavemap}).  Roughly speaking, the Coulomb gauge places the connection in the form $A = O( \nabla^{-1} F ) = O( \nabla^{-1}( \psi^2 ) )$, so the cubic wave equation now has the schematic form
$$ \Box \psi = O( \nabla^{-1}( \psi^2 ) \nabla \psi) + O( \psi^3 )$$
which turns out to be amenable to relatively simple Strichartz estimate techniques in four and higher dimensions.
In the special case of hyperbolic space targets, this approach was pushed (with
Substantial difficulty) to three and two dimensions in \cite{krieger:3d}, \cite{krieger:2d}, using much more sophisticated function spaces.  However, the Coulomb gauge actually becomes quite problematic to use here, due to the increasingly divergent nature of the inverse derivative operator $\nabla^{-1}$ in low dimensions at low frequencies.  This made it quite difficult to go beyond small data global regularity and obtain other expected features of the critical perturbation theory, such as a large data result, a usable blowup criterion, and a stability and well-posedness theory.  To resolve these issues, a more geometric ``caloric gauge'' was proposed in \cite{tao:forges}.  Re-interpreting an earlier microlocal gauge construction from \cite{tao:wavemaps}, \cite{tao:wavemaps2} by replacing the (discrete, linear) Littlewood-Paley projections with the (continuous, nonlinear) harmonic map heat flow propagator, it was shown in \cite{tao:forges} that the heat flow naturally induced a gauge which was slightly more regular than the Coulomb gauge, replacing the problematic bilinear form $O( \nabla^{-1}( \psi^2 ) )$ with a nonlinear paraproduct in which the inverse derivative was guaranteed to fall on the higher
frequency factor and thus staying relatively small.  Interestingly, the reliance on the heat flow means that the gauge extends to large data (unlike all previous gauges), provided that the heat flow is known to converge asymptotically to zero for this data (which is true, for instance, for surfaces of constant negative curvature, due to a classical result of Eells and Sampson \cite{eells}). It looks likely that this will lead to a satisfactory large data perturbation theory for critical wave maps in two dimensions; this in turn sets the stage for
the tools of preceding sections, such as induction on energy, to be brought to bear on the large data critical regularity problem in two dimensions, which is currently open except in the case of symmetric data.  This is currently work in progress by the author.  For further discussion of all of these issues on wave maps we refer to the recent survey \cite{igor}.

We close with a brief discussion of status of the corresponding critical regularity theory for the Yang-Mills and Klein-Gordon equations.  Here, many of the expected analogous results – for instance, that the four-dimensional Yang-Mills equations enjoy global regularity for any small energy data – are still open.  One of the main difficulties here is that the connections are significantly more curved than in the wave maps case; indeed, even after taking a good gauge such as the Coulomb gauge, the best thing that can be said about a connection is that it itself obeys a nonlinear wave equation.  One consequence of this is that even after selecting the gauge carefully, one cannot hope to dispense with the influence of the connection via an iteration argument.  Instead, one is forced to work with the connection as an integral part of the equation, and begin developing dispersive estimates for the \emph{covariant} wave equation $D^\alpha D_\alpha \phi = 0$.  This is now a problem in variable-coefficient liner equations rather than nonlinear PDE, and as such requires a rather different set of tools to those discussed above, namely the method of parametrices.  Such parametrices were developed in six and higher-dimensions, first for the Maxwell-Klein-Gordon equations in \cite{tao:igor} (which is simpler due to the abelian nature of the gauge group), and then for the non-abelian Yang-Mills equations in \cite{krieger}.  The basic idea is to construct certain ``distorted plane wave'' functions which almost solve the covariant wave equation, and then superimpose these waves together to create a parametrix (approximate solution) for the equation.  In order to ensure that the error terms accrued in this process are manageable, a large number of harmonic analysis preparations (such as Littlewood-Paley projections) have to be carefully performed first.  In the non-abelian case an additional difficulty arises because the distorted plane waves are obtained by solving a nonlinear ODE, and many regularity estimates on the solutions to that ODE must then be obtained.  See \cite{tao:igor}, \cite{krieger} for details.  The lower-dimensional cases, especially the energy-critical four-dimensional case, remain of great interest; it appears that the necessary step here is to develop covariant null form estimates, but there appear to be significant technical obstacles to doing so at present.

\end{document}